\documentclass[11pt,paper=a4]{article}
 \usepackage{indentfirst, latexsym,bm}
 \usepackage{amsmath}
 \usepackage{pifont}
 \usepackage{amsfonts}
 \usepackage{mathrsfs}
 \usepackage{array}
 \usepackage{multirow} 
 \usepackage{graphicx}
 \usepackage{subfigure}
 \usepackage{picinpar}
 \usepackage{setspace}
 \usepackage[textsize=footnotesize]{todonotes}
 \RequirePackage[colorlinks,citecolor=blue,urlcolor=blue,linkcolor=blue]{hyperref}
 \usepackage[top=2cm,bottom=2cm, outer=2cm, inner=2cm]{geometry}
 
 \usepackage{caption}
 \usepackage{subcaption}
 \usepackage[labelfont={bf,small},textfont={small}]{caption}
 
 \RequirePackage[numbers]{natbib}
 \usepackage{enumitem} 
 
 \usepackage{booktabs}
 \usepackage{authblk}
 
 \usepackage{threeparttable}
 \usepackage{mathrsfs,amsfonts,amsmath}
 
 \usepackage{color}
 
 \makeatletter
 \def\namedlabel#1#2{\begingroup
 	#2%
 	\def\@currentlabel{#2}%
 	\phantomsection\label{#1}\endgroup
 }
 \makeatother

 \numberwithin{figure}{section}

 \newcommand\email[1]{\href{mailto:#1}{ \nolinkurl{#1}}}

 \newtheorem{theorem}{Theorem}[section]
 \newtheorem{definition}[theorem]{Definition}
 \newtheorem{lemma}[theorem]{Lemma}
 \newtheorem{corollary}[theorem]{Corollary}
 \newtheorem{proposition}[theorem]{Proposition}
 \newtheorem{remark}[theorem]{Remark}
 \newtheorem{condition}[theorem]{Condition}
 \newtheorem{example}{Example}[section]

 \def\blemma{\begin{lemma}}\def\elemma{\end{lemma}}
 \def\bproposition{\begin{proposition}}\def\eproposition{\end{proposition}}
 \def\ttheorem{\begin{theorem}}\def\etheorem{\end{theorem}}
 \def\bcorollary{\begin{corollary}}\def\ecorollary{\end{corollary}}
 \def\bremark{\begin{remark}}\def\eremark{\end{remark}}
 \def\bcondition{\begin{condition}}\def\econdition{\end{condition}}

 \newtheorem{assumption}[theorem]{Assumption}

 \def\benumerate{\begin{enumerate}}\def\eenumerate{\end{enumerate}}
 \def\bitemize{\begin{itemize}}\def\eitemize{\end{itemize}}

 \def\beqlb{\begin{eqnarray}}\def\eeqlb{\end{eqnarray}}
 \def\beqnn{\begin{eqnarray*}}\def\eeqnn{\end{eqnarray*}}
 \def\nnm{\nonumber}\def\ar{\!\!\!&}
 
 \def\mrm{\mathrm}
 
 \def\proof{\noindent{\it Proof.~~}}\def\qed{\hfill$\Box$\medskip}

\begin{document} 
  \title{\bf  \Large Large Excursions of Reflected L\'evy Processes: Asymptotic Shapes}
 
 \author{Zhi-Hao Cui\footnote{School of Mathematical Sciences, Nankai University, China; email: cuizh.math@gmail.com} \quad  \  Hao Wu\footnote{School of Mathematical Sciences, Nankai University, China; email: wuhao.math@outlook.com} 
 	\quad\ and\quad 
 	Wei Xu\footnote{School of Mathematics and Statistics, Beijing Institute of Technology,  China; email: xuwei.math@gmail.com. Xu gratefully acknowledges financial support from the National Natural Science Foundation of China (No. 11531001) and the National Key R\&D Program of China (No. 2023YFA1010103).} 
 }  
 \maketitle

 \begin{abstract}
  This paper primarily investigates the geometric properties of  excursions of  L\'evy processes reflected at the past infimum with long lifetime or large height. 
  For an oscillating process in the domain of attraction of a stable law, our results state that excursions with a long lifetime need not have a large height. 
  After a suitable scaling, they behave like stable excursions  with lifetime or height greater than one. 
  These extend the related results in Doney and Rivero [{\it Prob. Theory Relat. Fields}, 157(1) (2013) 1-45].
  In contrast, for the negative-drift case we prove that under a heavy-tailed condition, long lifetime and large height are asymptotically equivalent. 
  Conditioned on either event, excursions converge under spatial scaling to a single-jump process with Pareto-distributed jump size and size-biased jump time. Moreover, after a suitable time rescaling, the effect of the negative drift becomes apparent.
 
 	\bigskip
 	
 	
 	\smallskip
 	

 \end{abstract}

  \section{Introduction and main results} 
 \label{Sec.Introduction}
 \setcounter{equation}{0}
 
 
 Consider a real-valued L\'evy process $X=\{ X_t:t\geq 0 \}$ with non-monotone c\`adl\`ag trajectories and characteristics $(b,\sigma,\nu)$.  
 For $x\in \mathbb{R}$, we write $ \tau^-_x$ and $ \tau^+_x$ for the first passage times of $X$ into $(-\infty,x)$ and $(x,\infty)$, respectively.
 Define the \textit{running infimum} and  \textit{supremum processes} by
 \beqnn
  \underline{X}_t:= \inf\big\{X_s: 0\leq s\leq t\big\}
 \quad \mbox{and}\quad 
 \overline{X}_t:= \sup\big\{X_s: 0\leq s\leq t\big\},\quad t\geq0. 
 \eeqnn 
 It is known that the reflected processes $X-\underline{X}$ and $\overline{X}-X$ are non-negative  strong Markov processes. 
 Their excursions away from $0$ form two Poisson point processes on the space $\mathcal{E}$ of excursions, with respective intensities $\underline{n}(d\epsilon)$ and $\overline{n}(d\epsilon)$.
 Here, each excursion $\epsilon \in \mathcal{E}$ is a non-negative c\`adl\`ag function on $\mathbb{R}_+$ with \textit{lifetime} $\zeta=\inf\{t\geq0: \epsilon_s=0, \forall s\geq t \}$ and \textit{height} $\overline{\epsilon}= \sup\{\epsilon_s: 0\leq s\leq \zeta\}$, with the usual convention that $\inf \emptyset =+\infty$.  
 We make convention that $\epsilon_s=0$ whenever $s\geq \zeta$. 
 

  

 This work is mainly concerned with the asymptotic properties of excursions with long lifetime or large height under the measures $\underline{n}(d\epsilon)$ and $\overline{n}(d\epsilon)$. 
 These objects appear in many contexts and have attracted substantial attention. Several representative examples are listed below: 
 \begin{enumerate}
  \item[$\bullet$] \textit{Path-decomposition.} For $t\geq 0$, the last passage time by $X$ at its infimum before $t$ is defined by
  \beqnn 
  \underline{g}_t:= \sup\big\{ s\leq t: X_s = \underline{X}_t
  \mbox{ or } X_{s-}= \underline{X}_t\big\}  .
  \eeqnn
  Conditionally on $ \underline{g}_t=s \in[0,t]$, the process $X$ over the time interval $[0,t]$ can be decomposed at time  $\underline{g}_t$ into the returned pre-$\underline{g}_t$ part and the post-$\underline{g}_t$ part, that are distributed as the laws $\overline{n}(\cdot \,|\, \zeta >s)$ and $\underline{n}(\cdot \,|\, \zeta >t-s)$, respectively; see e.g. \cite{Chaumont2013,Yano2013}. 
  The long-term behavior of $X$ can be analyzed via the asymptotics of long-lived excursions under $\underline{n}(d\epsilon)$ and $\overline{n}(d\epsilon)$, e.g., the past infimum $\underline{X}_t$ is asymptotically distributed as $ \underline{n}(\zeta >t)\cdot \int_0^\infty \overline{n}(\epsilon_s\in dx, \zeta >s)ds$ for large $t$; see \cite{ChaumontMałecki2016}. 
  
 \item[$\bullet$] \textit{Continuous-state branching processes.}  For a spectrally positive L\'evy process $X$, the Lamperti transform associates it with a continuous-state branching process $\mathcal{L}^X$; see e.g. Theorem~12.2 in \cite[p.337]{Kyprianou2014} and Theorem~10.10 in \cite[p.292]{Li2023}.  
 In particular, the total progeny and maximum  of $\mathcal{L}^X$ are identical in law to $(\tau_0^-, \overline{X}_{\tau_0^-})$, whose asymptotic relations can be understood with the help of conditional limit theorems for excursions under $\underline{n}(d\epsilon)$. 
  
  \item[$\bullet$] \textit{Splitting tree.} When $X$ is of bounded variation and has no negative jumps,  Lambert \cite{Lambert2010} established a bijection between excursions in $\mathcal{E}$ and splitting trees whose width process is a binary homogeneous Crump-Mode-Jagers process, which was generalized to the spectrally positve case later in \cite{LambertUribeBravo2018}. 
  Under this correspondence, the extinction time and total progeny of the binary tree coincide in law with the height $\overline{\epsilon}$ and lifetime $\zeta$ under $\underline{n}(d\epsilon)$. 
 \end{enumerate}
 Despite extensive work on their distribution properties and connections to other areas, the asymptotic properties of L\'evy excursions remain largely unexplored. 
 When $X$ is a Brownian motion,  the distributions of $\zeta$ and $\overline{\epsilon}$  are given explicitly by 
 \beqnn
 \underline{n}\big(  \zeta>t \big) = \frac{ t^{-1/2}}{\sqrt{\pi/2}}
 \quad \mbox{and}\quad 
 \underline{n} ( \overline{\epsilon}>x) =\frac{1}{x} , \quad x,t> 0;
 \eeqnn
 see Proposition~2.8 and Exercise~2.10 in \cite[p.484-485]{RevuzYor2005}. 
 By Williams' decomposition,  the joint law of $(\zeta, \overline{\epsilon})$ can be described by the maximum and length of two independent three-dimensional  Bessel processes, though not in closed form. 
 For stable processes, we also have $\underline{n}\big(  \zeta >t \big) =t^{-\varrho}/\Gamma(1-\varrho)$ with $\varrho= \mathbf{P}(X\leq 0)$. 
 However, to the best of our knowledge, the closed form for $\underline{n} ( \overline{\epsilon}>x) $ is unknown unless $X$ is one-sized. 
 For general oscillating L\'evy process in the domain of attraction of a stable law,  Doney and Rivero \cite{DoneyRivero2013}  obtained a precise local description of the distributions of the lifetime $\zeta$ and the final position $\epsilon_{\zeta-}$, and also used it to study the local probability of the first passage time $\tau_0^-$.

 Before stating the main results of this paper, we need to introduce several notation and terminology. 
 For $\kappa\in\mathbb{R}$, we write $f \in \mathrm{RV}^\infty_\kappa$ (resp. $f\in \mathrm{RV}^0_\kappa$) if  for any $x>0$,
 \beqnn
 \frac{f(tx)}{f(t)} \ar\to\ar  x^\kappa,\quad \mbox{as } t\to\infty \mbox{ (resp. as $t\to 0+$)}.
 \eeqnn 
 For a locally compact set $\mathcal{I}$, let $D(\mathcal{I};\mathbb{R})$ be the space of all c\`adl\`ag functions on $\mathcal{I}$ that is equipped with Skorohod's  topology. 
 A sequence of probability laws on $D(\mathcal{I};\mathbb{R})$ is said to   \textit{converge weakly} if it converges weakly on $D(I;\mathbb{R})$ for any compact subset $I\subset \mathcal{I}$; see Chapter 3 in \cite{Billingsley1999}. 
  \medskip
  
 \textit{\textbf {The oscillating case.}} 
 To obtain asymptotic results for excursions under $\underline{n}(d\epsilon)$, we impose the following standard assumption on $X$:
 \begin{assumption}\label{Assumption01}
 The process $X$ belongs to the domain of attraction of a stable law without centering, that is, there exists a positive function $\boldsymbol{c}$ on $\mathbb{R}_+$ such that 
 	\beqnn
 	\frac{X_t}{\boldsymbol{c}(t)} \ar\overset{\rm d}\to\ar Y_1, \quad \mbox{as $t\to\infty$},
 	\eeqnn
  where $Y$ is a strictly stable process with index $\alpha\in (0,2]$ and negativity parameter $\rho= \mathbf{P} \big(Y_1\leq 0\big)\in(0,1)$. 
 \end{assumption}
 
 It is well known that $\boldsymbol{c} \in \mathrm{RV}^\infty_{1/\alpha}$; see \cite[p.345]{BinghamGoldieTeugels1987}. 
 Under this assumption, the bivariate downward ladder process is in the domain of attraction of a bivariate stable law with index $(\rho, \alpha\rho )$. From  Spitzer's formula, it follows that 
 \beqnn
 \underline{n}(\zeta >t) \in \mathrm{RV}^\infty_{-\rho};
 \eeqnn 
 see \cite{DoneyRivero2013} and Section~\ref{Sec.RecurrentCase} for details.
 Quantities introduced to $X$ have analogues for $Y$ and are denoted by a superscript $Y$, e.g., $\underline{n}^Y$ is the excursion measure of  $Y$ reflected at its past infimum. 
  Our first main result states that after a suitable scaling,  excursions of $X-\underline{X}$ can be well-approximated by the corresponding stable excursions. 
 
 
  \begin{theorem}\label{MainThm01}
  The pushforward of $\underline{n}\big(\cdot\,|\,\zeta>t\big)$ by the scaling map $\epsilon\mapsto  \{ \epsilon_{ts}/\boldsymbol{c}(t):s\geq 0 \}$ converges weakly to $\underline{n}^Y\big(\cdot \,|\,\zeta>1\big)$ on $D\big([0,\infty);\mathbb{R}\big)$ as $t\to\infty$. 
 \end{theorem}
 
 The image of $\underline{n}^Y(\cdot\,|\, \zeta>1)$ under the restriction map $\omega \mapsto \{ \omega_s:s\in[0,1] \}$ is a probability law on $D\big([0,1];\mathbb{R}_+\big)$ and known as the law of \textit{$\alpha$-stable meander}. 
 Actually, Theorem~\ref{MainThm01} generalizes the conditional limit theorem established in  \cite{DoneyRivero2013} for the entrance law of excursions reflected at the minimum, which states that  as $t\to\infty$,
  \beqlb\label{eqn. StableMeander}
 \underline{n}\big( \epsilon_t/\boldsymbol{c}(t)\in dx\, \big|\, \zeta>t  \big) \to \underline{n}^Y\big(\epsilon_1\in dx\,\big|\, \zeta >1\big) . 
 \eeqlb
 Under the Spitzer’s condition, analogous results for a random walk $Z=\{ Z_n:n\geq 0 \}$ has been proved in various setting in \cite{AfanasyevGeigerKerstingVatutin2005,Doney1985,Durrett1978}. 
 Specially, for some regularly varying function $b_n\in \mathrm{RV}^\infty_{1/\alpha}$,  the rescaled process $\{Z_{[ns]}/b_n:s\in[0,1]\}$ conditioned on $\tau_0^->n$ converges weakly as $n\to\infty$ to the $\alpha$-stable meander on $D([0,1];\mathbb{R})$.


%
 
 
 The next corollary uses Theorem~\ref{MainThm01} to establish an analogous result for the L\'evy process $X$ started from a positive state and conditioned to stay positive.
 Define a c\`adl\`ag process $\mathcal{C}^{Y,1}:=\{ \mathcal{C}^{Y,1}_s:s\geq 0 \}$ as the concatenation of the stable process $Y$ with an independent stable meander.
 Specifically, the segment $\big\{ \mathcal{C}^{Y,1}_s: s\in [0,1]\big\}$ is distributed as the pushforward of $\underline{n}^Y(\cdot\,|\,\zeta>1)$ by the restriction $\omega\mapsto \{\omega(s): s\in[0,1]\}$. Conditionally on  $\mathcal{C}^{Y,1}_1=y$, the shifted process $\{\mathcal{C}^{Y,1}_{1+s}: s\geq 0\}$ is independent of the first segment and has the same law as $Y$ under $\mathbf{P}_y$.

 \begin{corollary}\label{MainCorollary01}
 	The rescaled process $ \{ X_{ts}/\boldsymbol{c}(t):s\geq 0 \}$ under $\mathbf{P}_x(\cdot\,|\,\tau_0^->t)$ converges weakly to $\mathcal{C}^{Y,1}$ in $D([0,\infty);\mathbb{R})$ as $t\to\infty$. 
 \end{corollary}

  \begin{figure}	
  	\centering
  	  	\begin{minipage}{0.27\textwidth}
  		\centering
  		\scalebox{0.8}{
  			\begin{tikzpicture}[baseline=(current bounding box.south)]
  				\def\r{1.5}
  				\def\d{1.5}
  				\draw[thick] (-\d/2,0) circle (\r);
  				\draw[thick,blue] (\d/2,0) circle (\r);
  				\node at (-\d/2-0.7, 0) {$\zeta >t$};
  				\node at (\d/2+0.75, 0) {$\overline{\epsilon} >c(t)$};
  		\end{tikzpicture}}
  	\end{minipage}
  	\hfill
  	\raisebox{-0.25cm}{%
  		\begin{minipage}{0.36\textwidth}
  			\centering
  			\scalebox{0.4}{
  				\begin{tikzpicture}[baseline=(current bounding box.south)]
  					\node[anchor=south west, inner sep=0] at (-1.8,-1.99)
  					{\includegraphics[width=12cm]{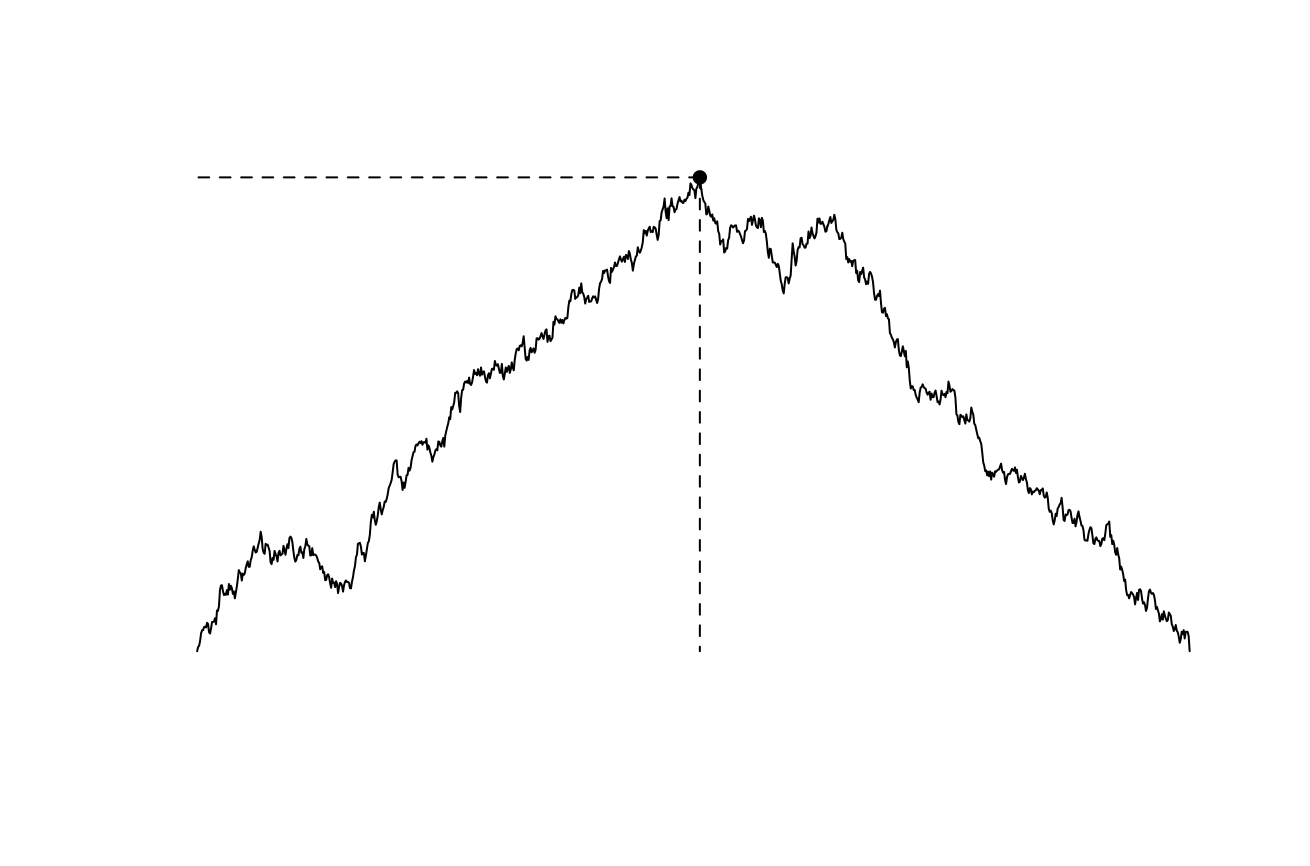}};
  					\node[anchor=south west, inner sep=0] at (-1.4,-1.5)
  					{\includegraphics[width=14cm,height=10cm]{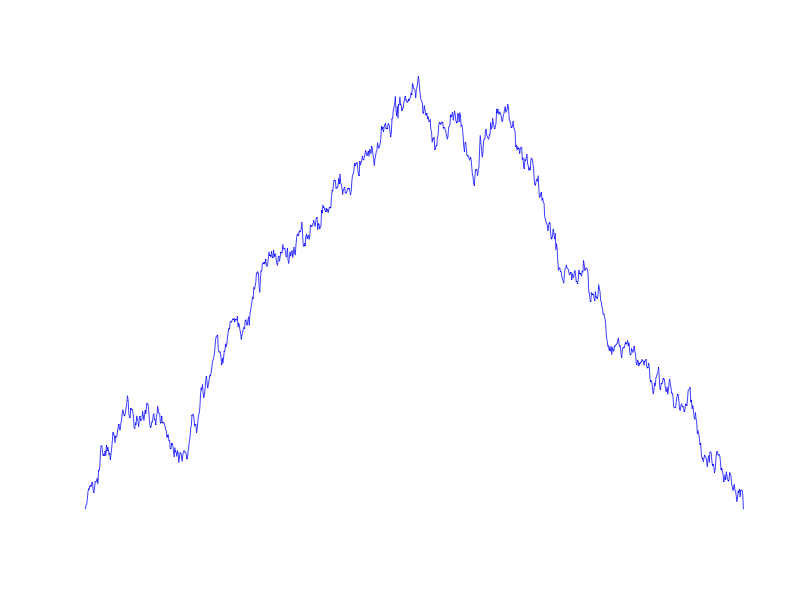}};
  					\node[anchor=south west, inner sep=0] at (-2.1,-0.6)
  					{\includegraphics[width=14cm,height=3.5cm]{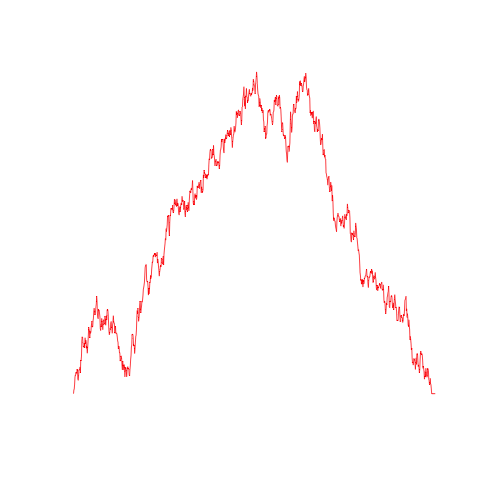}};
  					\node at (6,8) {\huge $\underline{n}(\epsilon_{ts}/\boldsymbol{c}(t)\mid \zeta>t)$};
  					\node at (11,-0.5) {\huge$\zeta$};
  					\draw[->, thick] (0,0) -- (0,8) node[left] {\huge$\overline{\epsilon}$};
  					\draw[->, thick] (0,0) -- (13,0) node[below] {\huge$t$};
  			\end{tikzpicture}}
  	\end{minipage}}
  	\hfill
  	\begin{minipage}{0.35\textwidth}
  		\centering
  		\scalebox{0.4}{
  			\begin{tikzpicture}[baseline=(current bounding box.south)]
  				\node[anchor=south west, inner sep=0] at (-1.8,-1.99)
  				{\includegraphics[width=12cm]{2.png}};
  				\node[anchor=south west, inner sep=0] at (-1.5,-1.8)
  				{\includegraphics[width=10cm,height=10cm]{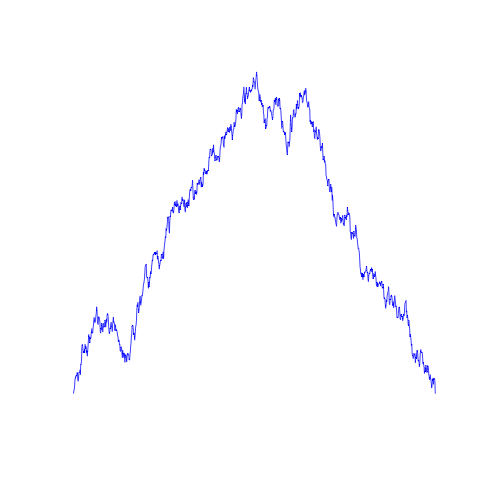}};
  				\node[anchor=south west, inner sep=0] at (-0.75,-2.15)
  				{\includegraphics[width=5cm,height=12cm]{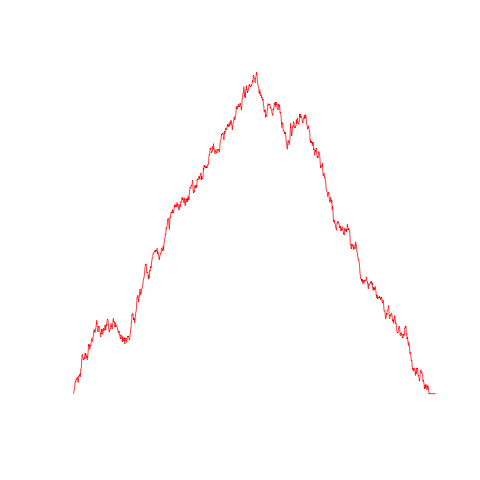}};
  				\node at (6,8) {\huge $\underline{n}(\epsilon_{ts}/\boldsymbol{c}(t)\mid \overline{\epsilon}>\boldsymbol{c}(t))$};
  				\node at (9,-0.5) {\huge$\zeta$};
  				\draw[->, thick] (0,0) -- (0,8) node[left] {\huge$\overline{\epsilon}$};
  				\draw[->, thick] (0,0) -- (11,0) node[below] {\huge$t$};
  		\end{tikzpicture}}
  	\end{minipage}
  	\captionsetup[figure]{labelfont=it,textfont={bf,it}}
  	\caption{The first sub-graph depicts the asymptotic relation between $\zeta>t$ and $\overline{\epsilon}>\boldmath{c}(t)$. The other two sub-graphs then show several typical rescaled excursion paths conditioned on each event.}\label{fig:1}
  \end{figure}

  Theorem~\ref{MainThm01} also implies that an excursion with long lifetime need not to reach a substantial height; see the second sub-graph in Figure~\ref{fig:1}. More precisely, there is a good chance that excursions remaining near the bottom can nonetheless survive for extended periods.
  Specifically,  for any $\delta >0$, 
  \beqnn
  \lim_{t\to\infty}\underline{n}\big(\overline{\epsilon}\leq \delta \cdot \boldsymbol{c}(t) \,|\, \zeta >t\big)=  \underline{n}^Y\big(\overline{\epsilon}\leq \delta  \,|\, \zeta >1\big)>0 .
  \eeqnn
  On the other hand, this limit also indicates that even if $\underline{n} \big(\overline{\epsilon}>\delta_0 \cdot \boldsymbol{c}(t) \big) \sim \underline{n}( \zeta >t)$ for some $\delta_0>0$, the events $\overline{\epsilon}>\delta_0 \cdot \boldsymbol{c}(t)$ and $\zeta >t$ are not asymptotically equivalent; see the first sub-graph in Figure~\ref{fig:1}.  
  This leads naturally to the following questions: Can excursions with large height also persist for a long time? 
  And, conditioned on large height, do scaling limit theorems hold under $\underline{n}(d\epsilon)$?

   In contrast to the extensive literature on the lifetime $\zeta$, the asymptotic behavior of  $\overline{\epsilon}$ has received less attention. 
  However, in the spectrally positive case,  the explicit connection between the scale function and $\underline{n}(\overline{\epsilon}>x)$ (see the identity (8.22) in \cite[p.239]{Kyprianou2014}) offers a pathway to analyze the law of $\overline{\epsilon}$ and address the questions above. 
  In this case, under Assumption~\ref{Assumption01}, we have that $\alpha \in (1,2]$ and $\alpha\rho =1$; see \cite[p.218]{Bertoin1996}.



  \begin{theorem} \label{MainThm03}
 	If $X$ is spectrally positive, then $\underline{n} \big(\overline{\epsilon}>\boldsymbol{c}(t)\big) 
 	\sim (\alpha-1) \cdot \Gamma(1-\rho)  \cdot 	\underline{n} \big(\zeta>t\big)  $, and under the scaling map $\epsilon\mapsto  \{ \epsilon_{ts}/\boldsymbol{c}(t):s\geq 0 \}$, the conditional law $\underline{n} \big(\cdot\,|\,\overline{\epsilon}>\boldsymbol{c}(t)\big)$ converges weakly to $\underline{n}^Y\big(\cdot \,|\,\overline{\epsilon}>1\big)$  on $D\big([0,\infty);\mathbb{R}\big)$ as $t\to\infty$. 
 
 \end{theorem}
 
 The conditional limit theorems in Theorem~\ref{MainThm01} and \ref{MainThm03} show that rescaled excursions with long lifetime or large height can be well approximated by stable excursions with distinct characteristics. 
 Nevertheless, their behavior can differ substantially. 
 For instance,excursions that remain near the bottom may nevertheless persist for a long time, while short-lived excursions can still reach a high level; see the last two sub-graphs in Figure~\ref{fig:1}. 
 Intuitively, as $\rho\to 0+$, excursions with large height tend to become thinner like the red curve in the final sub-graph.

 \begin{remark}
 The proof of Theorem~\ref{MainThm03}  heavily relies on the asymptotic relation that as $x\to\infty$, 
 \beqlb\label{eqn.102}
  \underline{n}(\overline{\epsilon}>x) \sim C\cdot x^{-1},
  \eeqlb
  This not only compares the long-term behavior of $\zeta$ and $\overline{\epsilon}$, but also allows us to transfer conditional limit theorems from $\underline{n}  (\cdot\,|\,\overline{\epsilon}>\boldsymbol{c}(t) )$ to  $\underline{n} (\cdot\,|\,\zeta>t )$. 
   For general L\'evy processes, excursions with large height are unlikely to jump down significantly before attaining their maximum, so their behavior should resemble the spectrally positive case. 
   It is therefore natural to conjecture that \eqref{eqn.102} and the conclusions of Theorem~\ref{MainThm03} remain valid. 
 	\end{remark}

 \textit{\textbf {The negative-drift case.}} 
 We now present asymptotic results for excursions under $\underline{n}(d\epsilon)$ when the L\'evy process $X$ has negative drift 
 \beqnn
 \beta:=-\mathbf{E}\big[X_1\big]\in (0,\infty).
 \eeqnn
 These results are proved under the following heavy-tailed assumption on  the tail-distribution of L\'evy measure $\overline{\nu}(x):= \nu\big([x,\infty)\big)$ for $x>0$. 
 \begin{assumption}\label{Assumption02}
  The state $0$ is regular for $(0,\infty)$ and $\overline{\nu} \in \mathrm{RV}^\infty_{-\theta}$ for some $\theta>1$.
 \end{assumption}

 In contrast to the oscillating case, we show in the next theorem that excursions with long lifetime necessarily reach high levels and vice versa. 
 Since $X$ drifts to $-\infty$,  Theorem~6.9(ii) in \cite{Kyprianou2014} deduces that $\underline{n} ( \zeta )<\infty$. 
 
 \begin{theorem}\label{MainThm04}
 Under $\underline{n}(d\epsilon)$,  the events $\zeta>t$ and $\overline\epsilon >\beta t$ are asymptotically equivalent as $t\to\infty$,  that is, 
 	\beqnn
 	\underline{n}\big( \zeta > t\big)
 	\sim
 	\underline{n}\big(\bar\epsilon >\beta t\big)
 	\sim \underline{n}\big( \zeta\big) \cdot \bar{\nu}\big(\beta t\big) 
 	\quad \mbox{and}\quad 
 	\underline{n}\big(\bar\epsilon >\beta t \,\big|\, \zeta >t \big)
 	\sim
 	\underline{n}\big( \zeta >t \,\big|\, \bar\epsilon >\beta t \big) \to 1.
 	\eeqnn
 \end{theorem}
 
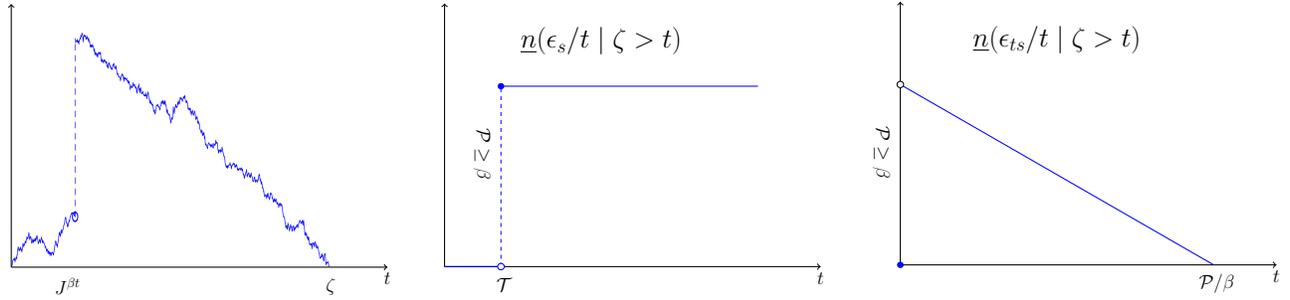
\begin{figure}	
 	\centering 
 		\hspace*{-1cm}
 		\raisebox{-0.14cm}{
 			\begin{minipage}{0.37\textwidth}
 				\centering
 				\scalebox{0.5}{
 					\begin{tikzpicture}[baseline=(current bounding box.south)]
 						\node[anchor=south west, inner sep=0] (image) at (-2.1,-3.15)
 						{\includegraphics[width=12cm, height=12cm]{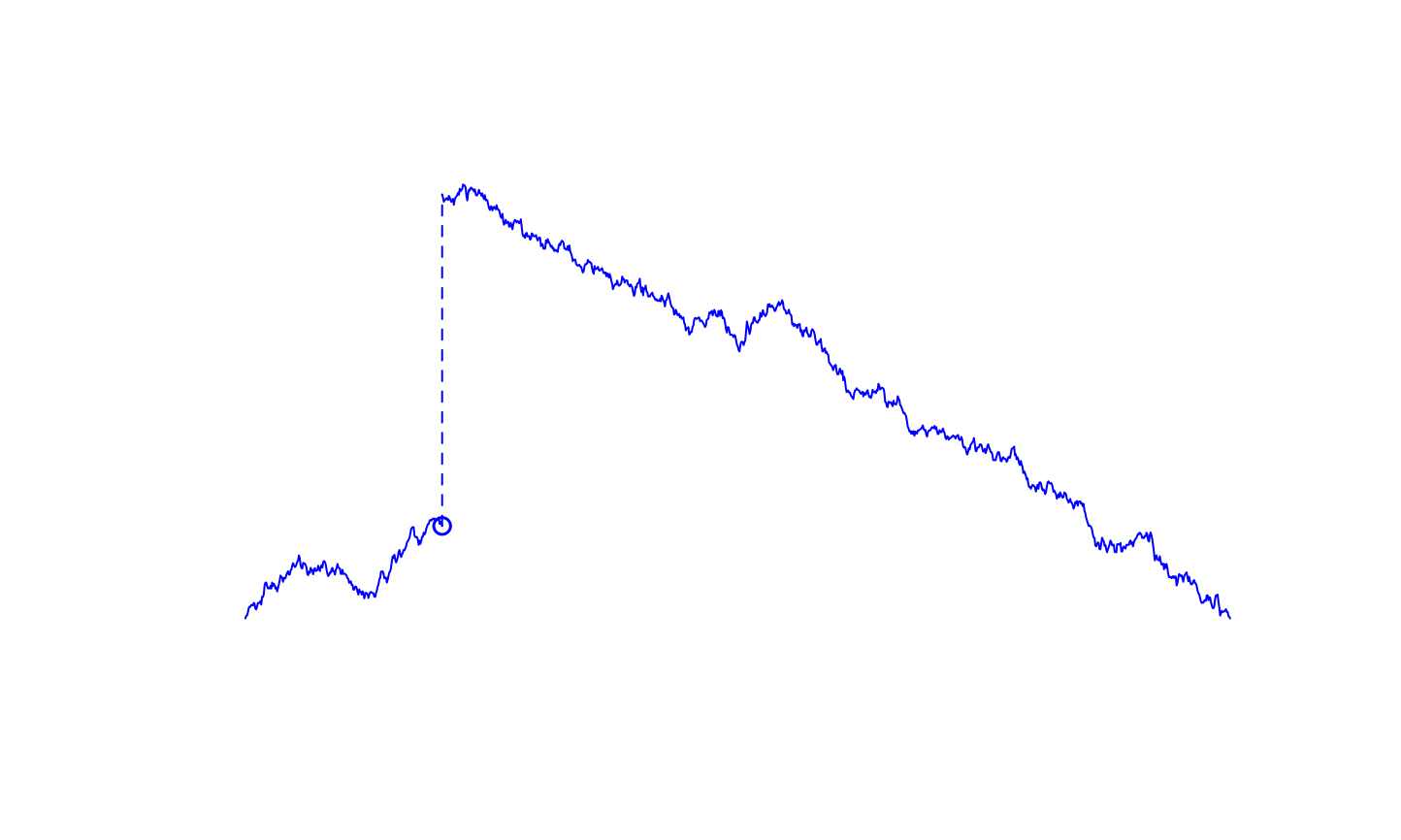}};
 						\node at (1.5,-0.5) {\Large$J^{\beta t}$};
 						\node at (8.5,-0.5) {\Large$\zeta$};
 						\draw[->, thick] (0,0) -- (0,7);
 						\draw[->, thick] (0,0) -- (10,0) node[below] {\Large$t$};
 				\end{tikzpicture}}
 			\end{minipage}
 		}
 		\hfill
 		\begin{minipage}{0.28\textwidth}
 			\centering
 			\scalebox{0.5}{
 				\begin{tikzpicture}[baseline=(current bounding box.south)]
 					\draw[->, thick] (0,0) -- (0,7) node[ rotate=270] at (1,3) {\Large$\mathcal{P}\geq \beta$};
 					\draw[->, thick] (0,0) -- (10,0) node[below] {\Large$t$};
 					
 					\begin{scope}[xscale=8.3333, yscale=8]
 						\node at (0.5,0.75) {\huge $\underline{n}(\epsilon_{s}/t\mid \zeta >t)$};
 						\draw[blue, dashed] (0.18, 0) -- (0.18, 0.6);
 						\draw[blue, thick,domain=0:0.18,smooth,samples=200] plot(\x,0);
 						\draw[blue, thick,domain=0.18:1,smooth,samples=200] plot(\x,0.6);
 						\draw[blue, thick, fill=white] (0.18,0) circle (0.3pt);
 						\node[below right] at (0.15,-0.02) {\Large $\mathcal{T}$};
 						\fill[blue](0.18,0.6) circle (0.3pt);
 					\end{scope}
 			\end{tikzpicture}}
 		\end{minipage}
 		\hfill
 		\begin{minipage}{0.35\textwidth}
 			\centering
 			\scalebox{0.5}{
 				\begin{tikzpicture}[baseline=(current bounding box.south)]
 					\draw[->, thick] (0,0) -- (0,7) node[ rotate=270] at (-0.5,3) {\Large$\mathcal{P}\geq \beta$};
 					\draw[->, thick] (0,0) -- (10,0) node[below] {\Large$t$};
 					
 					\begin{scope}[xscale=8.333, yscale=8] 
 						\node at (0.5,0.75) {\huge $\underline{n}(\epsilon_{ts}/t\mid \zeta >t)$};
 						\draw[blue, thick,domain=0:1,smooth,samples=200] plot(\x,{0.6-0.6*\x});
 						\fill[blue](0,0.6) circle (0.3pt);
 						\node at (1.01,-0.06) {\Large$\mathcal{P}/\beta$};
 						\draw[fill=white] (0,0.6) circle (0.3pt); 
 						\fill[blue](0,0) circle (0.3pt);
 					\end{scope}
 			\end{tikzpicture}}
 		\end{minipage} 
 	\caption{The first sub-graph shows a representative excursion whose long lifetime is sustained by a single large jump. The other two sub-graphs illustrate that, after spatial scaling, each large excursion reduces to a single jump trajectory, while the effect of its negative drift becomes visible only after an additional time scaling.}\label{fig:2}
 \end{figure}

 This phenomenon mainly arises from the interplay between the negative drift and the single-big-jump principle; see the first sub-graph in Figure~\ref{fig:2}. 
 More precisely, the negative drift tends to kill excursions quickly, and only sufficiently large positive jumps can counteract it and allow the excursion to persist. 
 Conversely, if an excursion experiences a significant jump before being extinguished by the negative drift, it will typically survive for a long duration. 
 Under the heavy-tailed condition, the single-big-jump principle asserts that the occurrence of more than one big jump is extremely unlikely. 
 Additionally,  to offset the negative drift, the single large jump is expected to occur as soon as possible. 
 
 We now formulate another main result in this work that not only highlights the key role of the single-big-jump principle, but also precisely describes the arrival time and size of the single big jump.  
 Consider two independent positive random variables $\mathcal{P} $  and $\mathcal{T} $ that are distributed as 
 \beqnn
 \mathbf{P}\big( \mathcal{P}  \leq x \big) = 1- \big( x/\beta \big)^{-\theta}
 \quad \mbox{and}\quad 
 \mathbf{P}\big( \mathcal{T}  \leq t \big) = \frac{\underline{n}( \zeta\wedge t)}{\underline{n}( \zeta)}.
 \eeqnn
 
  \begin{theorem}\label{MainThm07}
 	The following two limit results hold as $t\to\infty$. 
 	\begin{enumerate}
 		\item[(1)]  The pushforward of $\underline{n}( \cdot \,|\, \zeta > t)$ by the scaling map $\epsilon\mapsto  \{ \epsilon_{s}/t :s\geq 0 \}$ converges weakly to the law of $\{ \mathcal{P}  \cdot \mathbf{1}_{\{s\geq \mathcal{T} \}}:  s\geq 0 \}$ on $D\big([0,\infty);\mathbb{R}\big)$.
 		
 		\item[(2)] The pushforward of $\underline{n}( \cdot \,|\, \zeta > t)$ by the scaling map  $\epsilon\mapsto   \{ \epsilon_{ts}/t :s>0 \}$  converges weakly to the law of $\{ (\mathcal{P}  -\beta s)\vee 0: s> 0 \}$  on $D\big((0,\infty);\mathbb{R}\big)$. 
 	\end{enumerate}
 \end{theorem}

 Since the events  $\zeta>t$ and $\overline\epsilon >\beta t$ are asymptotically equivalent,  the two scaling limits in Theorem~\ref{MainThm07} also hold when the conditioning measure $\underline{n}( \cdot \,|\, \zeta > t)$ replaced by $\underline{n}(\cdot\,|\, \overline\epsilon >\beta t)$. 
 The first conclusion reveals that the arrival time and size of the single large jump are  asymptotically independent, which can be attributed to the independent increments of $X$; see the second sub-graph in Figure~\ref{fig:2}.  
 The limiting size-biased distribution of the arrival time arises from offsetting the negative drift early in the excursion. 
 In the heavy-tailed setting, it is unsurprising that the rescaled jump size follows a Pareto distribution. 
 Comparing these two scaling limits shows that the effect of negative drift after the arrival of big jump becomes visible only after long-term accumulation;   see the third sub-graph in Figure~\ref{fig:2}.

 Using the preceding theorem, the next corollary establishes an analogous result for a L\'evy process started from a positive state $x$ and conditioned to stay positive.
 The main difference is that the arrival time of the single-big-jump is asymptotically equal in law to the random variable $\mathcal{T}_x$ with size-biased distribution given by
 \beqnn
  \mathbf{P}\big(\mathcal{T}_x\leq t\big) =\frac{\mathbf{E}_x\big[\tau_0^-\wedge t\big]}{\mathbf{E}_x\big[\tau_0^- \big] } . 
 \eeqnn
 Alternatively, the random variable $\mathcal{T}_x$ also can be constructed from $\mathcal{T}$ and the potential measure $\widehat{U}(dt,dy)$ of the bivariate downward ladder process (see \eqref{eqn.U}) via
 \beqnn
  \mathbf{P}\big(\mathcal{T}_x\leq t\big) = \int_0^t \mathbf{P}\big( \mathcal{T}\leq t-s\big)\, \frac{\widehat{U}\big(ds,[0,x]\big)}{\widehat{U}\big(\mathbb{R}_+,[0,x]\big)} .
 \eeqnn
 A proof is given in  Lemma~\ref{Lemma.408}. 
 
  \begin{corollary}\label{MainThm08}
 	For any $x>0$, the following two limit results hold as $t\to\infty$.
 	\begin{enumerate}
 		\item[(1)]  The spatially scaled process $
 		\{X_s/t : s\geq 0\}$ under $\mathbf{P}_x(\cdot\,|\, \tau_0^->t)$ converges weakly to $\{ \mathcal{P} \cdot \mathbf{1}_{\{s\geq \mathcal{T}_x \}}: s\geq 0 \}$ in $D\big([0,\infty);\mathbb{R}\big)$. 
 		
 		\item[(2)]  The rescaled process $
 		\{X_{ts}/t : s> 0\}$ under $\mathbf{P}_x(\cdot\,|\, \tau_0^->t)$ converges weakly to $\{ \mathcal{P} -\beta s: s> 0 \}$ in $D\big((0,\infty);\mathbb{R}\big)$.  
 	\end{enumerate}
 \end{corollary}
 
 The first conclusion was firstly established by Xu \cite{Xu2021a} using a sophisticated approach based on the fluctuation theory of L\'evy processes. 
 The discrete analogue of our second scaling limit for random walks with finite variance was proved earlier by Durrett \cite{Durrett1980}.
 
  
   \smallskip
  {\it \textbf{Organization of this paper.}}
  In Section~\ref{Sec.Preliminaries}, we introduce some further notation, recall some elementary fluctuation theory for L\'evy process and then present three auxiliary lemmas used in the proof of our main results. 
  Section~\ref{Sec.RecurrentCase} is devoted to the proofs of
  Theorem~\ref{MainThm01}, \ref{MainThm03} and Corollary~\ref{MainCorollary01}. 
  Meanwhile, proofs of Theorem~\ref{MainThm04}, \ref{MainThm07} and Corollary~\ref{MainThm08} can be found in Section~\ref{Sec.TransientCase}.

  \section{Preliminaries} 
 \label{Sec.Preliminaries}
 \setcounter{equation}{0}

  
  This section begins by introducing additional notation and recalling basic fluctuation theory for Lévy processes, for which we refer to \cite{Bertoin1996,Doney2007,Kyprianou2014}. We then present several auxiliary lemmas that will be essential for proving our main results.
 All processes in this work are defined on a complete probability space
 $(\Omega,\mathscr{F},\mathbf{P})$ equipped with a filtration $\{\mathscr{F}_t\}_{t\geq 0}$ satisfying the usual hypotheses. 
 For $x\in\mathbb{R}$, let $\mathbf{P}_x$ and $\mathbf{E}_x$ denote the law and expectation of the  process starting from $x$. 
 For simplicity, we also write $\mathbf{P}= \mathbf{P}_0$ and $\mathbf{E}=\mathbf{E}_0$. 
  Let   $\overset{\rm a.s.}\to$, $\overset{\rm d}\to$  and  $\overset{\rm p}\to$ be the  almost sure convergence,  convergence in distribution and convergence in probability respectively.
 We also use $\overset{\rm a.s.}=$, $\overset{\rm d}=$ and $\overset{\rm p}=$ to denote almost sure equality, equality in distribution and equality in probability respectively.

 \subsection{Fluctuation theory}

  For $t\geq 0$, let $\Delta X_t := X_t-X_{t-}$ denote the jump size of $X$ at time $t$. 
 The point measure of all jumps 
 \beqnn
 \sum_{t>0} \mathbf{1}_{\{ \Delta X_t\neq 0 \}}\cdot \delta_{(t,\Delta X_t)}(ds,dy),
 \eeqnn
 is a Poisson random measure on $(0,\infty)\times \mathbb{R}$ with intensity $ds\, \nu(dy)$.  
 By the compensation formula for Poisson random measures (see \cite[p.7]{Bertoin1996}), we have that for any predictable process $\mathbf{F}=\{ \mathbf{F}_t:t\geq 0 \}$ taking values in the space of non-negative measurable functions on $\mathbb{R}$, 
  \beqlb\label{eqn.CompenstationX}
  \mathbf{E}\bigg[ \sum_{t>0} \mathbf{F}_t(\Delta X_t) \cdot \mathbf{1}_{\{ \Delta X_t\neq 0 \}}  \bigg]= \mathbf{E}\bigg[\int_0^\infty dt \int_{\mathbb{R}\setminus\{0\}}\mathbf             {F}_t(y)\, \nu(dy)\bigg].
  \eeqlb
  
 The reflected process $X-\underline{X}$ is a Markov process with Feller transition semigroup given by $\mathbf{P}_x(X_t\in dy, \underline{X}_t>0)$ for any $x,y>0$; see, e.g., Proposition~1 in \cite[p.156]{Bertoin1996}. 
 Its excursions away from $0$ form a Poisson point process on $\mathcal{E}$ with intensity $\underline{n}(d\epsilon)$.  
 Under $\underline{n}(d\epsilon)$, the excursion process is Markovian and shares the same semigroup of $X$ killed upon entering $(-\infty,0]$, that is, for any $t,s>0$, any measurable functionals $F$ on $D\big([0,t];\mathbb{R}\big)$ and  $G$ on $D\big([0,s];\mathbb{R}\big)$, 
 \beqlb\label{eqn.n-Markov}
 \lefteqn{\underline{n}\Big(F\big(\epsilon_r, r \in[0,t]\big) \cdot G\big(\epsilon_r, r\in[t,t+s]\big),\, \zeta>t+s\Big)}\ar\ar\cr
 \ar\ar\cr
 \ar=\ar \underline{n}\Big(F\big(\epsilon_r, r \leq t\big)\cdot \mathbf{E}_{\epsilon_t}\big[G\big(X_r, r\in[0,s]\big),\,  \underline{X}_s>0\big],\, \zeta>t\Big)  .
 \eeqlb
 An important consequence of \eqref{eqn.CompenstationX} and \eqref{eqn.n-Markov} is that the compensation formula also holds under $\underline{n}(d\epsilon)$, i.e., for any predictable process $\mathbf{F}=\{ \mathbf{F}_t:t\geq 0 \}$ taking values in the space of non-negative measurable functions on $\mathbb{R}$, 
 \beqlb\label{eqn.Compensation}
 \underline{n}\bigg(\sum_{0<s\leq \zeta} \mathbf{F}_s(\Delta \epsilon_s) \cdot 1_{\{\Delta \epsilon_s \neq 0\}}\bigg)
 = \underline{n}\bigg(\int_0^\zeta ds \int_{\mathbb{R}\setminus\{0\}} \mathbf{F}_s(y)\,\nu(dy)\bigg).
 \eeqlb



 Let $\widehat{L} = \{\widehat{L}_t: t\ge 0\}$ be the local time at zero of $X-\underline{X}$ in the sense of \cite[p.109]{Bertoin1996}.
 Without loss of generality, we normalize it so that 
 \beqnn
  \mathbf{E}\bigg[ \int_0^\infty e^{-t} d \widehat{L}_t\bigg] =1. 
 \eeqnn
 Its right inverse $\widehat{L}^{-1}$ is referred to as the \textit{downward ladder time process}, and $\widehat{H}= \{ \widehat{H}_t :t\geq 0 \}$ with $\widehat{H}_t := X_{\widehat{L}^{-1}_t}$ is called the \textit{downward ladder height process}. 
 Lemma~2 in \cite[p.157]{Bertoin1996} states that $(\widehat{L}^{-1},\widehat{H})$ is a L\'evy process (possibly killed at an exponential rate) with \textit{bivariate exponent}  given by
 \beqnn
  \widehat\kappa(\lambda,u) =\log \mathbf{E}\Big[e^{-\lambda \widehat{L}^{-1}_1- u\widehat{H}_1} \Big]
  \ar=\ar \exp\Big\{\int_0^\infty \frac{dt}{t}\int_{[0,\infty)}(e^{-t}-e^{-\lambda  t-ux})\mathbf{P}(X_t\in dx)\Big\},\quad \lambda,u\geq 0.
 \eeqnn
  In particular, both $\widehat L^{-1}$ and $\widehat H$ are (possibly killed) subordinators with drifts $\underline{\tt d},\underline{\tt a}\geq 0$.
  Moreover, $\underline{\tt d}=0$ if and only if the state $0$ is regular for   $(0,\infty)$, for $X$. 
  By (6.16) in \cite[p.170]{Kyprianou2014}, the bivariate exponent $\widehat\kappa(\lambda, u)$ also admits the  representation
  \beqlb\label{eqn.kappa}
   \widehat \kappa(\lambda,u) \ar=\ar \underline{\tt d}\cdot \lambda + \underline{\tt a}\cdot u + \int_0^\infty \int_0^\infty (1- e^{-\lambda t-u x})\,\underline{n}(\zeta \in dt,\epsilon_\zeta \in dx). 
  \eeqlb
 Denote by $\widehat{U}(dt,dy)$ the potential measure of $(\widehat{L}^{-1},\widehat{H})$, defined as the $\sigma$-finite measure on $\mathbb{R}^2_+$ given by
 \beqlb\label{eqn.U}
 \widehat{U}(dt,dy):= \int_0^\infty \mathbf{P}\big( \widehat{L}^{-1}_s \in  dt,\widehat{H}_s\in dy \big)\,ds.
 \eeqlb
 In particular, its marginal cumulative function in space is also known as the \textit{renewal function} associated with the ladder height process $\widehat{H}$, that is, 
 \beqnn
  \widehat{V}(x) := \widehat{U}\big(\mathbb{R}_+,[0,x]\big)= \int_0^\infty \mathbf{P}\big( \widehat{H}_s \leq x \big)\,ds,\quad x\geq 0;
 \eeqnn 
 see e.g.  \cite[Chapter VI]{Bertoin1996} and \cite[Chapter 6]{Kyprianou2014}. 
 By \eqref{eqn.kappa} and \eqref{eqn.U}, the function $\widehat{V}$ also satisfies that
 \beqlb\label{eqn.447}
 \widehat{V}(x) = \bar{\tt d} +  \int_0^\infty  \overline{n}(  \epsilon_s\leq x, \zeta>s)\, ds ,\quad x\geq 0;
 \eeqlb
 see also (16) in \cite{DoneyRivero2013}. 
 Analogously, 
 quantities associated with $\overline{X}$ and $\overline{X}-X$ are denoted without hats, e.g., $\kappa(\lambda,u)$, $U(dt,dy)$ and $V(x)$,  and satisfy analogous properties.

 \subsection{Several auxiliary lemmas}
 
 
 We now present several auxiliary lemmas essential for proving our main results. The first extends Theorem~5 in \cite{Chaumont2013}, which decomposes 
 $X$ at its last passage time  $\underline{g}_t$ into two independent excursions. Its proof follows directly from the monotone class theorem and is sketched below.

 \begin{lemma}\label{Lemma.201}	
 	For any $t\geq 0$ and  any bounded measurable functional $F$ on   $\mathbb{R}_+\times \mathbf{D}\big([0,t];\mathbb{R}\big)$, we have 
 \beqlb\label{eqn.ExcursionRep}
 \mathbf{E}\big[F(\underline{g}_t,X)\big]  
 \ar=\ar  \underline{\tt d}\cdot  \int_{\mathcal{E}} F \big( t, \overleftarrow{\epsilon}^t \big) \, \overline{n}(d\epsilon,\zeta>t) + \bar{\tt d}\cdot  \int_{\mathcal{E}} F \big( 0,\epsilon  \big) \, \underline{n}(d\epsilon ,\zeta>t) \cr
 \ar\ar  + \int_{(0,t)}ds \int_{\mathcal{E}}  \overline{n}(d\epsilon, \zeta>s)  \int_{\mathcal{E}}F\big(s,  (\epsilon,\epsilon^*)^s \big) \, \underline{n}(d\epsilon^*, \zeta>t-s).
 \eeqlb
 where $F(s,\omega):= F\big(s, \{\omega_r:r\in[0,t]\}\big)$, $ \overleftarrow{\omega}^s_r :=\omega_{(s-r)-}- \omega_s\wedge\omega_{s-}$ and  $(\omega,\omega^*)^s_r:= \overleftarrow{\omega}^s_r \cdot \mathbf{1}_{\{ r\leq s \}} + \omega^*_{r-s}\cdot \mathbf{1}_{\{ r\geq s \}}$ for $s,r\geq 0$ and $\omega,\omega^*\in D([0,\infty);\mathbb{R})$. 
 \end{lemma}
 
 Analogously to $\mathcal{C}^{Y,1}$, for each $t>0$,  we define a c\`adl\`ag process $\mathcal{C}^{X,t}:=\{ \mathcal{C}^{X,t}_s:s\geq 0 \}$ as the concatenation of $X$ and an independent copy of its meander, that is, 
 \begin{enumerate}
  \item[(1)] The segment $\big\{ \mathcal{C}^{X,t}_s: s\in [0,t]\big\}$ is distributed as the pushforward of $\underline{n}(\cdot\,|\,\zeta>t)$ by the restriction map  $\omega\mapsto \big\{\omega(s): s\in[0,t]
  \big\}$;
  
  \item[(2)] Given  $\mathcal{C}^{X,t}_t=y$, the shifted process
  $\{\mathcal{C}^{X,t}_{t+s}: s\geq 0\}$ is independent of  $\big\{\mathcal{C}^{X,t}_s :s\in[0,t] \big\}$ and has the same law as $X$ under $\mathbf{P}_y$.  
 \end{enumerate}
 Our second auxiliary lemma shows that the process $\mathcal{C}^{X,t}$ stopped at $\sigma_0^{X,t}:=\inf\{s>0:\mathcal{C}^{X,t}_{s}\leq 0\}$ realizes  the probability law $\underline{n}(\cdot\,|\,\zeta >t)$ on $D([0,\infty);\mathbb{R})$.

 \begin{lemma}\label{Lemma.302}
	For any $T\geq 0$ and any bounded measurable functional $H$  on   $D([0,T];\mathbb{R})$,  we have  
 	\beqlb\label{eqn.221}
 	\mathbf{E}\Big[ H\Big(\mathcal{C}^{X,t}_{s\wedge \sigma_0^{X,t}}\vee 0: s\in [0,T]\Big)\Big]=\underline{n}\big(H\big(\epsilon_s:s\in [0,T]\big) \,\big|\,\zeta >t\big). 
 	\eeqlb
 	
 \end{lemma}
 \proof 
 First, note that the definition of $\mathcal{C}^{X,t}$ implies that $\sigma_0>t $ almost surely, so \eqref{eqn.221} holds for $T\leq t$.  
 For $T>t$,  by the monotone class theorem it suffices to prove that for any  bounded measurable functionals  $F$ on $D([0,t];\mathbb{R})$ and $G$  on  $D([0,T-t];\mathbb{R})$, 
 \beqlb\label{eqn4.44}
 \mathbf{E}\big[F\big(\mathcal{C}^{X,t}_{s}, s\in [0,t]\big)\cdot G\big(\mathcal{C}^{X,t}_{s\wedge \sigma_0^{X,t}}\vee 0: s\in [t,T]\big)\big]
 =\underline{n}\big(F\big(\epsilon_{s}, s\in [0,t]\big)\cdot G\big(\epsilon_s,s\in [t,T]\big)\,\big|\,\zeta >t\big).
 \eeqlb
 We first consider the right-hand side that can be written as
 \beqnn
 \frac{1}{\underline{n} ( \zeta >t )}\int_{t}^{\infty} \int_0^\infty\underline{n}\big(F\big(\epsilon_{s}, s\in [0,t]\big)\cdot G\big(\epsilon_s, s\in [t,T]\big),\zeta \in dr,\epsilon_t\in dy\big),
 \eeqnn
 which, by using the Markov property to the last integral at time $t$, equals to 
 \beqnn 
 \frac{1}{\underline{n} ( \zeta >t )}\int_{0}^{\infty} \underline{n}\big(F\big(\epsilon_{s}, s\in [0,t]\big), \epsilon_{t}\in dy,\zeta >t\big) \int_{t}^{\infty} \mathbf{E}_y\big[G\big(X_{s\wedge (r-t)}\vee 0, s\in [0,T-t]\big)\big), \sigma_0^X\in dr-t\big] ,
 \eeqnn
 where $\sigma_0^X:= \inf\{ s\geq 0: X_s\leq 0 \}$. 
 It is easy to identify that the inner integral equals to $\mathbf{E}_y\big[G\big(X_{s\wedge \sigma_0^X}\vee 0: s\in [0,T-t]\big)\big)\big]$ and then the right-hand side of \eqref{eqn4.44} equals to 
 \beqlb\label{eqn.448}
 \int_{0}^{\infty}\mathbf{E}_y\big[G\big(X_{s\wedge \sigma_0^X}\vee 0, s\in [0,T-t]\big)\big)\big] \cdot \underline{n}\big(F\big(\epsilon_{s}, s\in [0,t]\big), \epsilon_{t}\in dy\,\big|\,\zeta >t\big).
 \eeqlb
 For the left-hand side of \eqref{eqn4.44}, it also can be written as 
 \beqnn
 \lefteqn{\int_{t}^{\infty}\int_{0}^{\infty}\mathbf{E}\big[F\big(\mathcal{C}^{X,t}_{s}, s\in [0,t]\big)\cdot G\big(\mathcal{C}^{X,t}_{t+s}\cdot \mathbf{1}_{\{t+s< r\}} , s\in [0,T-t]\big), \sigma_0^{X,t} \in dr, \mathcal{C}^{X,t}_t\in dy\big]}\ar\ar\cr
 \ar=\ar \int_{t}^{\infty}\int_{0}^{\infty}\mathbf{E}\big[F\big(\mathcal{C}^{X,t}_{s},s\in [0,t]\big)\cdot G\big(\mathcal{C}^{X,t}_{t+s}\cdot \mathbf{1}_{\{t+s< r\}} , s\in [0,T-t]\big),  \sigma_0^{X,t} \in dr\, \big|\, \mathcal{C}^{X,t}_t = y\big]\mathbf{P}\big(\mathcal{C}^{X,t}_t\in dy\big).
 \eeqnn
 By the conditional independence between the two segments of $\mathcal{C}^{X,t}$ before and after time $t$, it further equals to 
 \beqnn
 \int_{0}^{\infty} \mathbf{E}\big[F\big(\mathcal{C}^{X,t}_{s}, s\in [0,t]\big), \mathcal{C}^{X,t}_t\in dy\big]\int_{t}^{\infty}\mathbf{E}\big[ G\big(\mathcal{C}^{X,t}_{t+s} \cdot \mathbf{1}_{\{t+s< r\}} , s\in [0,T-t]\big), \sigma_0^{X,t} \in dr\, \big|\, \mathcal{C}^{X,t}_t= y\big]. 
 \eeqnn
 Since 
 $\{\mathcal{C}^{X,t}_{t+s}: s\geq 0\}$ given $\mathcal{C}^{X,t}_t = y$ equals in law to $X$ under $\mathbf{P}_y$, it also equals to 
 \beqnn
 \ar\ar\int_0^{\infty} \mathbf{E}\big[F\big(\mathcal{C}^{X,t}_{s}, s\in [0,t]\big),  \mathcal{C}^{X,t}_t\in dy\big] \int_{0}^{\infty}\mathbf{E}_y\big[ G\big(X_{s}\cdot \mathbf{1}_{\{s< r\}} ,s\in [0,T-t]\big), \sigma_0^X \in dr\big] .
 \eeqnn
 It is easy to see that the inner integral equals to $\mathbf{E}_y\big[G\big(X_{s\wedge \sigma_0^X}\vee 0,s\in [0,T-t]\big) \big]$. 
 Since $\{ \mathcal{C}^{X,t}_s: s\in [0,t] \}$ is distributed as the pushforward of $\underline{n}(\cdot\,|\,\zeta>t)$ by the mapping $\omega\mapsto \{\omega(s): s\in[0,t]\}$, we also have 
 \beqnn
 \mathbf{E}\big[F\big(\mathcal{C}^{X,t}_{s}, s\in [0,t]\big),  \mathcal{C}^{X,t}_t\in dy\big] =\underline{n}\big(F\big(\epsilon_{s},s\in [0,t]\big), \epsilon_{t}\in dy \,\big|\,\zeta >t\big). 
 \eeqnn
 Putting them together, we see that the left-hand side of \eqref{eqn4.44} also equals to \eqref{eqn.448}.  The proof ends.  
 \qed
 
 As a third auxiliary lemma, we prove that if a sequence of c\`ad\`ag processes converges weakly and the limit process has the property that $0$ is regular for $(-\infty,0)$, then the weak convergence is inherited by the corresponding processes stopped at the first passage into $(-\infty,0]$.
 
 
  \begin{lemma}\label{Lemma.304}
 	Consider processes $\xi^t$, $\xi^* $ such that $\xi^t\overset{\rm d}\to \xi^*$   in $D\big([0,\infty);\mathbb{R}\big)$ as $t\to\infty$. 
 	Let $\sigma^t_0$ and $\sigma^{*}_0$ be their  first passage times into $(-\infty,0]$ respectively. 
 	If $\sigma^{*}_0<\infty$ a.s. and $0$ is regular for $(-\infty,0)$, for $\xi^*$,
 	then $ \xi^t_{\cdot \wedge \sigma^{t}_0} \overset{\rm d}\to  \xi^*_{\cdot \wedge \sigma^{*}_0}$  in $D\big([0,\infty);\mathbb{R}\big)$ as $t\to\infty$. 
 \end{lemma}
 \proof  
 By Skorokhod's representation theorem, we may assume that $\xi^{(t)}\overset{\rm a.s.}\to \xi^*$  in $D([0,\infty);\mathbb{R})$.
 Consequently, for any $T\geq 0$ there exists random time changes $\{ \lambda^t_s:s\geq 0 \}_{t>0}$ with $ \lambda^t_0\overset{\rm a.s.}=0$ and $ \lambda^t_T\overset{\rm a.s.}=T$ such that as $t\to\infty$,
 \beqlb\label{eqn.341}
 \sup_{s\in[0,T]} \big| \lambda^t_s-s \big| + \sup_{s\in[0,T]}\Big| \xi^t_{\lambda^t_s} - \xi^*_s \Big| \overset{\rm a.s.}\to 0,
 \eeqlb
 Hence it suffices to prove that $\sigma^{t}_0\overset{\rm a.s.}\to \sigma^{*}_0$ as $t\to\infty$. 
 For any $s>0$, if $\sigma^{*}_0>s$, we see that $\underline{\xi}^*_s>0$ and hence $\underline{\xi}^t_s>0$ for all $t$ large enough, which yields that $\sigma^t_0\geq s$. 
 Since $s$ is arbitrary, we obtain that
 \beqnn
 \liminf_{t\to\infty}\sigma^t_0 \geq  \sigma^*_0.
 \eeqnn 
 On the other hand, if $\sigma^{*}_0<s$, we have $\underline{\xi}^*_s\leq0$.
 Since $0$ is regular for $(-\infty,0)$ for $\xi^*$, there exists some constant $a<0$ such that $\underline{\xi}^*_s\leq a$.
 By \eqref{eqn.341}, we have $\underline{\xi}^t_s\leq a/2$ for all large $t$ and hence $\sigma^t_0\leq s$. 
 This along with the arbitrariness of $s$ yields that 
 \beqnn
 \limsup_{t\to\infty}\sigma^t_0\leq \sigma^*_0.
 \eeqnn 
 Hence $\sigma^{t}_0\overset{\rm a.s.}\to \sigma^{*}_0$, which completes the proof. 
 \qed

  \section{The oscillating case} 
 \label{Sec.RecurrentCase}
 \setcounter{equation}{0}
 
 In this section, we prove our main asymptotic results for oscillating L\'evy processes under Assumption~\ref{Assumption01}. 
 As the preparation, we first recall several auxiliary properties of $X$.   
 For the quantities introduced to $X$, the corresponding ones for $Y$ are denoted by a superscript $Y$. Many of these admit explicit expressions; e.g. for $t>0 $ and $ x,\lambda\geq 0$, 
 \beqlb\label{eqn.4006}
 \widehat{\kappa}^Y(0,\lambda)=\widehat{\kappa}^Y(0,1) \cdot \lambda^{\alpha\rho},\quad \underline{n}^Y \big( \zeta >t \big) =\frac{t^{-\rho}}{\Gamma(1-\rho)}   \quad\mbox{and}\quad	\widehat{V}^Y(x)= \frac{ x^{\alpha\rho}}{	\Gamma\big(1+\alpha \rho\big) \cdot  \widehat{\kappa}^Y(0,1)},
 \eeqlb
 as described in Chapter~VIII of \cite{Bertoin1996}. 
 From Lemma 14 and its proof in \cite{DoneyRivero2013}, we obtain that as $t\to\infty$,
  \beqlb\label{eqn.4002}
  \widehat{V}(t) \in \mathrm{RV}^\infty_{\alpha\rho} 
  \quad \mbox{and}\quad
  \Gamma(1-\rho)\cdot \widehat\kappa^Y(0,1)\cdot \underline{n}(\zeta>t) \sim \widehat{\kappa}\big(0,1/\boldsymbol{c}(t)\big) \sim \frac{1}{\Gamma(1+\alpha\rho)\widehat{V}(\boldsymbol{c}(t))}.
 \eeqlb
Moreover, for any $x> 0$,  Proposition~6 of \cite{DoneyRivero2013} gives constants $C_1,C_2>0$ such that as $t\to\infty$, 
 \beqlb\label{eqn.334} 
 \overline{n}(\epsilon_t\leq x,\zeta >t) 
 \ar\sim\ar \frac{C_1}{t\cdot \boldsymbol{c}(t)}  
 \quad\mbox{and}\quad 
 \underline{n}(\epsilon_t\leq x,\zeta >t) 
 \sim \frac{C_2}{t\cdot \boldsymbol{c}(t)}. 
 \eeqlb
 
 For $q>0$, let $e(q)$ denote  an exponentially distributed random time with parameter $q>0$, independent of $X$.  
 From the last equality in the proof of Theorem~20 in \cite[p.176]{Bertoin1996} we have as $q\to 0+$,
 \beqlb\label{eqn.202}
 \mathbf{E}_x\big[1-e^{-q\tau_0^-}\big]
 =\mathbf{P}\big(-\underline{X}_{e(q)} \leq x\big)\sim \widehat{\kappa}(q, 0)\cdot \widehat{V}(x).
 \eeqlb 
 Theorem~14(iv) and Proposition 1 in \cite[p.169, 74]{Bertoin1996} give that $\widehat{\kappa}(q, 0)\sim  \Gamma(1-\rho) \cdot  \underline{n}(\zeta>1/q) \in \mathrm{RV}^0_{\rho}$. 
 Applying Karamata's Tauberian theorem (see e.g. Corollary 8.1.7 in \cite[p.334]{BinghamGoldieTeugels1987})  yields that as $t\to\infty$, 
 \beqlb\label{eqn.4009}
 \mathbf{P}_x(\tau_0^{-}>t)\sim \widehat{V}(x)\cdot\frac{\hat{\kappa}(1/t,0)}{\Gamma(1-\rho)}\sim \widehat{V}(x)\cdot\underline{n}(\zeta >t) \in \mathrm{RV}^\infty_{-\rho}.
 \eeqlb
 
 Finally, we recall the absolute continuity relation between the law of 
 $X$ conditioned to stay positive and the excursion measure $\underline{n}(d\epsilon)$, established in \cite[Theorem 3]{Chaumont1996}. Define the positive strong Markov process $X^{\uparrow}$ by
  \beqnn
 \mathbf{E}_x\big[G(X^{\uparrow}_t)] =   \mathbf{E}_x\bigg[ \frac{\widehat{V}(X_t)}{\widehat{V}(x)} \cdot G(X_t), \tau_0^->t \bigg],\quad x>0,\, t\geq 0,
 \eeqnn
 for any measurable function $G$ on $\mathbb{R}$. 
 Then  for any $t\geq 0$ and any measurable functional $F$ on $D\big([0,t];\mathbb{R}\big)$,
 \beqlb\label{eqn.443}
 \underline{n}\big(F(\epsilon), \zeta>t\big)= \mathbf{E} \bigg[\frac{F(X^\uparrow)}{\widehat{V}(X_t^{\uparrow})}\bigg] := \lim_{x\to 0+}\mathbf{E}_x \bigg[\frac{F(X^\uparrow)}{\widehat{V}(X_t^{\uparrow})}\bigg].
 \eeqlb
 Analogously, one defines the process $Y^{\uparrow}$, and the same relation also holds for $(\underline{n}^Y,Y^\uparrow )$. 
  
 With these preparations, we now prove Theorem~\ref{MainThm01}, Theorem~\ref{MainThm03}, and Corollary~\ref{MainCorollary01} in detail.
  For any $t>0$ and $\omega\in D([0,\infty);\mathbb{R})$, define the rescaled path $\omega^{(t)}:=\{\omega_{ts}/\boldsymbol{c}(t):s\geq 0 \}$. 
  Recall the convention that for a measurable functional $F$ on $D([0,t];\mathbb{R})$, we write $F(\omega):= F(\omega_s:s\in[0,t])$. 
 In terms of the process $\mathcal{C}^{X,t}$ from Lemma~\ref{Lemma.302}, we consider its rescaled version
 \beqnn
 \mathcal{C}^{X,(t)}_s:= \frac{\mathcal{C}^{X,t}_{ts}}{\boldsymbol{c}(t)},\quad s\geq 0.
 \eeqnn 
 Recall $\mathcal{C}^{Y,1}$ defined above  Corollary~\ref{MainCorollary01}.  
 The conditional independence of the two segments of $\mathcal{C}^{X,(t)}$ before and after time $1$ prompts us to decompose it at time $1$ and then show that each part converges weakly to the corresponding part of $\mathcal{C}^{Y,1}$.
 In the next proposition, we first prove that rescaled meander under $\underline{n}(\cdot\,|\, \zeta>t)$ converge weakly to the stable meander  by using the absolute continuity relationship \eqref{eqn.443}.

 \begin{proposition}\label{Prop.401}
 	For any  bounded, uniformly continuous functional $G$ on $D([0,1];\mathbb{R})$, we have  as $t\to\infty$,
 	\beqlb\label{eqn.333}
 	\underline{n}\big(G  (\epsilon^{(t)}  )\,\big|\, \zeta>t\big) 
 	\to
 	\underline{n}^Y\big(G (\epsilon  ) \,\big|\, \zeta>1\big).
 	\eeqlb
 	\end{proposition}
 	\proof  
 	For  $\delta>0$,   the left-hand side of \eqref{eqn.333} can be decomposed into the next two terms:
 	\beqnn 
 	\underline{n}\big(G(\epsilon^{(t)}),\epsilon^{(t)}_1<\delta \,|\, \zeta>t\big) 
 	\quad \mbox{and}\quad 
 	\underline{n}\big(G(\epsilon^{(t)}),\epsilon^{(t)}_1\geq  \delta \,|\, \zeta>t\big).
 	\eeqnn 
 	Firstly, the upper bound of $G$ allows us to bound the first term by $C \cdot \underline{n}\big(\epsilon^{(t)}_1< \delta \,|\, \zeta>t\big)$ uniformly in $t,\delta>0$.  
 	Additionally, by \eqref{eqn. StableMeander} we have as $t\to\infty$,
 	\beqnn
 	\underline{n}\big( \epsilon^{(t)}_1< \delta \,\big|\, \zeta>t\big)
 	\to  \underline{n}^Y\big(\epsilon_1< \delta \,\big|\, \zeta>1\big),
 	\eeqnn
 	which  vanishes as $\delta \to 0+$. This follows that 
 	\beqlb\label{eqn.301}
 	\lim_{\delta \to 0+} \limsup_{n\to\infty}  \underline{n}\big(G (\epsilon^{(t)}),\epsilon^{(t)}_1<\delta \,\big|\, \zeta>t\big) =0.
 	\eeqlb
 For the second term, an application of \eqref{eqn.443} with $F(\omega)=G(\omega)\cdot \mathbf{1}_{\{ \omega_1\geq \delta \}}$ shows that
 \beqlb\label{eqn.331}
  \underline{n}\big(G(\epsilon^{(t)}),\epsilon^{(t)}_1\geq  \delta \,|\,  \zeta>t\big)
  \ar=\ar \frac{1}{ \underline{n}(\zeta >t)} \cdot \mathbf{E}\bigg[ \frac{G (X^{\uparrow,(t)})}{\widehat{V}\big(X^{\uparrow,(t)}_1\cdot \boldsymbol{c}(t)\big)} , X^{\uparrow,(t)}_1\geq \delta\bigg] . 
 \eeqlb
  Repeating the proof of Theorem~4 in \cite{ChaumontDoney2010} with the random walk replaced by $X$, we can obtain that $X^{\uparrow,(t)}\to Y^{\uparrow}$ weakly in $D\big([0,\infty);\mathbb{R}\big)$. 
 By Skorokhod's representation theorem, we may assume that this convergence holds almost surely in $D\big([0,\infty);\mathbb{R}\big)$.
 This along with the fact that $\widehat{V} \in \mathrm{RV}^\infty_{\alpha\rho}$ induces that 
 \beqnn
   \widehat{V}\big(X^{\uparrow,(t)}_1\cdot \boldsymbol{c}(t)\big) \sim (Y^{\uparrow}_1)^{\alpha\rho} \cdot \widehat{V}\big( \boldsymbol{c}(t)\big), 
 \eeqnn  
 as $t\to\infty$.
 This allows us to use the dominated convergence theorem to \eqref{eqn.331} and then apply the preceding two asymptotic results to obtain that as  $t\to\infty$, 
 \beqnn
  \underline{n}\big(G(\epsilon^{(t)}),\epsilon^{(t)}_1\geq  \delta \,|\,  \zeta>t\big)
  \ar\sim\ar \frac{1}{ \underline{n}(\zeta >t)} \cdot \mathbf{E}\bigg[ \frac{G (Y^{\uparrow})}{(Y^{\uparrow}_1)^{\alpha\rho} \cdot \widehat{V}\big( \boldsymbol{c}(t)\big)} , Y^{\uparrow}_1\geq \delta\bigg] \cr
  \ar\to\ar \Gamma\big(1+\alpha \rho\big) \cdot \Gamma(1-\rho)\cdot \widehat{\kappa}^Y(0,1)\cdot  \mathbf{E}\bigg[ \frac{G (Y^{\uparrow})}{(Y_1^{\uparrow})^{\alpha\rho}} , Y_1^{\uparrow}\geq \delta\bigg]\cr
  \ar=\ar \Gamma(1-\rho)\cdot \mathbf{E}\bigg[ \frac{G (Y^{\uparrow})}{\widehat{V}^Y(Y_1^{\uparrow})} , Y_1^{\uparrow}\geq \delta\bigg].
 \eeqnn
 Here the limit and the last equality follow from \eqref{eqn.4002} and  \eqref{eqn.4006} respectively. 
 Finally, by using the absolute continuity relationship \eqref{eqn.443} for $(\underline{n}^Y,Y^\uparrow )$ and the fact $  \Gamma(1-\rho) \cdot  \underline{n}^Y \big( \zeta >1 \big)=1$; see \eqref{eqn.4006} with $t=1$,  
 	\beqnn
 	\lim_{t\to\infty}\underline{n}\big(G(\epsilon^{(t)}),\epsilon^{(t)}_1\geq  \delta \,|\,  \zeta>t\big) 
  =\Gamma(1-\rho)\cdot  \underline{n}^Y \big( G(\epsilon),\epsilon_1\geq  \delta, \zeta >1 \big)
  = \underline{n}^Y \big( G(\epsilon),\epsilon_1\geq  \delta\,\big|\, \zeta >1 \big).
 	\eeqnn
 Combining this together with \eqref{eqn.301} and then using the arbitrariness of $\delta$, we can get the desired limit \eqref{eqn.333} immediately.
 \qed

 \textit{\textbf{Proof of Theorem~\ref{MainThm01}.}} 
 Denote by  $\sigma_0^{(t)}$ and $\sigma_0^Y$ the first passage times of $\mathcal{C}^{X,(t)}$ and $\mathcal{C}^{Y,1}$ into $(-\infty,0]$ respectively.  
 A direct corollary of Lemma~\ref{Lemma.302} shows that  the stopped process
 $\mathcal{C}^{X,(t)}_{\cdot \wedge \sigma_0^{(t)}}$ is distributed as the push-forward of $\underline{n}\big(\cdot\,|\,\zeta>t\big)$ by the scaling map $\epsilon\mapsto  \epsilon^{(t)} $; meanwhile, the stopped process $\mathcal{C}^{Y,1}_{\cdot \wedge \sigma_0^{Y}}$ is distributed as  $\underline{n}^Y\big(\cdot\,|\,\zeta>1\big)$. 
 Since $\sigma_0^{Y}>1$ and $\mathcal{C}^{Y,1}$ evolves as $Y$ after time $1$, it is obvious that $0$ is regular for $(-\infty,0)$, for $\mathcal{C}^{Y,1}$. 
 Hence, by Lemma~\ref{Lemma.304} it suffices to prove the weak convergence of $\mathcal{C}^{X,(t)}$ to $\mathcal{C}^{Y,1}$ in $D\big([0,\infty);\mathbb{R}\big)$ as $t\to\infty$. 
 Furthermore, since the two segments of $\mathcal{C}^{X,(t)}$ and $\mathcal{C}^{Y,1}$ before and after time $1$ are conditionally independent, we just need to prove separately that as $t\to\infty$,
 \beqnn
 \big\{\mathcal{C}^{X,(t)}_{s}: s\in[0,1]\big\} \overset{\rm d}\to  \big\{\mathcal{C}^{Y,1}_{s}: s\in[0,1]\big\}
 \quad \mbox{and}\quad 
 \big\{\mathcal{C}^{X,(t)}_{1+s}-\mathcal{C}^{X,(t)}_1: s\geq 0\big\} \overset{\rm d}\to  \big\{\mathcal{C}^{Y,1}_{1+s}-\mathcal{C}^{Y,1}_1: s\geq 0\big\},
 \eeqnn
 respectively in $D\big([0,1];\mathbb{R}\big)$ and $D\big([0,\infty);\mathbb{R}\big)$.
 The second limit is a direct consequence of the two facts that its two sides are equal in law to $X^{(t)}$ and $Y$ respectively, and $X$ is in the domain of attraction of $Y$. 
 For the first one, it is sufficient to prove that  for any  bounded, uniformly continuous non-negative functional $G$ on $D([0,1];\mathbb{R})$,
 \beqnn
 \mathbf{E}\big[ G\big(\mathcal{C}^{X,(t)}\big) \big] \to \mathbf{E}\big[ G\big(\mathcal{C}^{Y,1} \big) \big],
 \eeqnn 
 which follows directly from Proposition~\ref{Prop.401}.
 \qed

 To prove Corollary~\ref{MainCorollary01}, we need the following well-known asymptotic result for the convolution of two regularly varying functions; see e.g. \cite{AsmussenFossKorshunov2003,Cline1986}. 
 \begin{proposition}\label{Prop.Convoluton}
 	Consider three positive functions $f_1$, $f_2$ and $g$ on $\mathbb{R}_+$ such that $f_1(x)\sim C_1\cdot g(x)$ and $f_2(x)\sim C_2\cdot g(x)$ for some constants $C_1,C_2\geq 0$. If $g$ is regularly varying at infinity, we have as $x\to\infty$, 
 	\beqnn
 	\frac{1}{g(x)}\int_0^x f_1(x-y)f_2(y)\,dy \to  C_1 \int_0^\infty f_2(y)\,dy + C_2 \int_0^\infty f_1(y)\,dy.
 	\eeqnn 
 	\end{proposition}

 \textit{\textbf{Proof of Corollary~\ref{MainCorollary01}.}} 
 Firstly, the independent increments of  $X^{(t)}$ allows us to decompose it at time $1$ into the following two independent parts 
 \beqnn 
 X^{(t)}_s=  X^{(t)}_{1\wedge s} + \widetilde{X}^{(t)}_{(s-1)\vee 0} 
 \quad \mbox{with}\quad 
 \widetilde{X}^{(t)}_s:= X^{(t)}_{1+s}- X^{(t)}_1 ,\quad s\geq0.
 \eeqnn
 Similarly as in the proof of Theorem~\ref{MainThm01}, it suffices to prove separately
 that 
 $\widetilde{X}^{(t)} \overset{\rm d}\to Y$ 
 in $D(\mathbb{R}_+;\mathbb{R})$ and  
 $\big\{X^{(t)}_s: s\in[0,1]\big\}$ under $\mathbf{P}_x(\cdot\,|\,\tau_0^->t)$ converge weakly to $\{\mathcal{C}^{Y,1}_s:s\in[0,1]\}$ in $D([0,1]; \mathbb{R}_+)$ as $t\to\infty$.   
 The first limit is obvious  since $\widetilde{X}^{(t)}\overset{\rm d}=X^{(t)}$.  
 For the second one, we just need to  identify that for any bounded, uniformly continuous non-negative functional $G$ on  $D([0,1];\mathbb{R})$,
 \beqlb\label{eqn.323}
 \mathbf{E}_x\big[ G\big( X^{(t)} \big)\,\big|\, \tau_0^->t\big]  
 \to   \underline{n}^Y\big(G(\epsilon)\,|\,\zeta>1\big) ,
 \eeqlb
 as $t\to\infty$. 
 Without loss of generality, we assume that $G(\cdot)\leq 1$. 
 For any constant $K>0$,   we decompose the left-hand side at time $\underline{g}_t$ into the following two terms:
 \beqlb\label{eqn.318}
 \mathbf{E}_x\big[ G\big( X^{(t)}\big); \underline{g}_t> K\,\big|\,\tau_0^->t \big]  
 \quad \mbox{and}\quad   
  \mathbf{E}_x\big[ G\big( X^{(t)}\big); \underline{g}_t\leq K\,\big|\, \tau_0^->t \big] . 
 \eeqlb
 In the next two steps, we prove that they converge to $0$ and $ \underline{n}^Y\big(G(\epsilon)\,|\,\zeta>1\big)$ respectively as $t\to\infty$ and then $K\to\infty$. 
 \medskip
 
 {\bf Step 1.} The first  conditional expectation in \eqref{eqn.318} can be bounded by
 \beqnn
 \mathbf{E}_x\big[  \underline{g}_t>K\,\big|\, \tau_0^->t \big] 
 = \frac{\mathbf{E}_x\big[  \underline{g}_t>K, \tau_0^->t \big] }{\mathbf{P}_x(\tau_0^->t)}
 = \frac{\mathbf{E}\big[  \underline{g}_t>K, \tau_{-x}^->t \big] }{\mathbf{P}_x(\tau_0^->t)}
 =\frac{\mathbf{E}\big[ \underline{g}_t>K, \underline{X}_t\geq -x \big]}{\mathbf{P}_x(\tau_0^->t)},
 \eeqnn
 which, by using  \eqref{eqn.ExcursionRep} with  $F(\underline{g}_t,X)=\mathbf{1}_{\{ \underline{g}_t>K \}}\cdot \mathbf{1}_{\{\underline{X}_t\geq -x\}}= \mathbf{1}_{\{ \underline{g}_t>K \}}\cdot \mathbf{1}_{\{X_{\underline{g}_t} \wedge X_{\underline{g}_t-} \geq -x\}}$, equals to
 \beqlb\label{eqn.420}
 \frac{1}{\mathbf{P}_x(\tau_0^->t)} \cdot \Big( \underline{
 \tt d}\cdot \overline{n}( \epsilon_t\leq x,  \zeta>t) +   
 \int_{(K,t)}  \overline{n}\big( \epsilon_s\leq x,  \zeta>s \big) \cdot  \underline{n}\big(\zeta>t-s\big)  \, ds   \Big) . 
 \eeqlb
 Since $    \overline{n}( \epsilon_t\leq x,  \zeta>t ) \in {\rm RV}^\infty_{-1-1/\alpha} $ and $ \mathbf{P}_x(\tau_0^{-}>t) \sim \widehat{V}(x)\cdot\underline{n}(\zeta >t) \in \mathrm{RV}^\infty_{-\rho}$; see \eqref{eqn.334} and \eqref{eqn.4009}, we have  
 \beqnn
 \int_1^\infty \overline{n}(  \epsilon_t\leq x,  \zeta>t )\, dt <\infty
 \quad \mbox{and}\quad 
  \lim_{t\to\infty}\frac{\overline{n}( \epsilon_t\leq x,  \zeta>t)}{\mathbf{P}_x(\tau_0^->t)}  =0.
 \eeqnn
 These allow us to apply Proposition~\ref{Prop.Convoluton} to the last integral in \eqref{eqn.420} and get that as $t\to\infty$, 
 \beqnn
 \frac{1}{\mathbf{P}_x(\tau_0^->t)}  \int_{(K,t)}  \overline{n}(  \epsilon_s\leq x,  \zeta>s ) \underline{n}(\zeta>t-s)\, ds \to  \frac{1}{ \widehat{V}(x)}\int_{K}^\infty  \overline{n}( \epsilon_s\leq x,  \zeta>s )\, ds ,
 \eeqnn
 which  vanishes as $K\to\infty$.
 Combining the preceding two limits together, we get that 
 \beqlb\label{eqn.322}
 \lim_{K\to\infty} \limsup_{t\to\infty}  \mathbf{E}_x\big[ G\big( X^{(t)}\big); \underline{g}_t> K\,\big|\,\tau_0^->t \big]  
  =0. 
 \eeqlb

 {\bf Step 2.}
 For the second conditional expectation in \eqref{eqn.318}.
 By the spatial homogeneity of $X$ we have 
 \beqnn
  \mathbf{E}_x\big[ G\big( X^{(t)}\big); \underline{g}_t \leq K\,\big|\, \tau_0^->t \big] 
  \ar=\ar \frac{1}{ \mathbf{P}_x\big(  \tau_0^->t \big)}
  \cdot  \mathbf{E}\big[ G\big( X^{(t)}+ x/\boldsymbol{c}(t)\big); \underline{g}_t\leq K, \underline{X}_t\geq x \big]. 
 \eeqnn
 Since $F$ is uniformly continuous and $\boldsymbol{c}(t)\to\infty$,  the  expectation on the right-hand side equals to
 \beqnn 
  \mathbf{E}\big[ G\big( X^{(t)} \big); \underline{g}_t\leq K, \underline{X}_t\geq -x \big] 
 +O\big(x/ \mathbf{c}(t)\big) \cdot \mathbf{P}\big(\underline{g}_t\leq K, \underline{X}_t\geq -x \big) . 
 \eeqnn
 Note that $\mathbf{P}\big(\underline{g}_t\leq K, \underline{X}_t\geq -x \big) \leq \mathbf{P}\big( \underline{X}_t\geq -x \big)= \mathbf{P}_x\big(\tau_0^-> t \big)$, we see that as $t\to\infty$,
 \beqlb\label{eqn.340}
  \mathbf{E}_x\big[ G\big( X^{(t)}\big); \underline{g}_t\leq K\,\big|\, \tau_0^->t \big] 
 \sim \frac{1}{\mathbf{P}_x\big( \tau_0^->t \big)}\cdot 
 \mathbf{E}\big[ G\big( X^{(t)} \big); \underline{g}_t\leq K, \underline{X}_t\geq -x \big]  . 
 \eeqlb
We write $ I_K(t)$ for the last expectation. 
 By using \eqref{eqn.ExcursionRep} with $F(\underline{g}_t,X)=  G\big(X^{(t)}\big) \cdot \mathbf{1}_{\{\underline{g}_t\in[0,K], X_{\underline{g}_t}\wedge X_{\underline{g}_t-}\geq -x \}}$,  
 \beqlb\label{eqn.5001}
  I_K(t)= \bar{\mathtt{d}} \cdot  \underline{n}\big(G(\epsilon^{(t)}), \zeta>t\big) + \int_0^K ds \int_{\mathcal{E}}   \overline{n}(d\epsilon, \epsilon_s\leq x, \zeta>s)  \int_{\mathcal{E}}G\big( \widetilde{\epsilon}_r^{(t),s}\big)\, 
 \, \underline{n}(d\epsilon^*, \zeta>t-s) .
 \eeqlb
 with $\widetilde{\epsilon}_r^{(t),s}:=  (\epsilon,\epsilon^{*})^{s}_{tr}/\boldsymbol{c}(t)$ for $r \in [0,1]$. 
 By Theorem~\ref{MainThm01}  and  \eqref{eqn.4009}, we  first have  as $t\to\infty$,
 \beqlb\label{eqn.339}
 \frac{  \bar{\mathtt{d}}\cdot  \underline{n}\big(G(\epsilon^{(t)}), \zeta>t\big)}{\mathbf{P}_x(\tau_0^->t)} =  \bar{\mathtt{d}} \cdot \frac{  \underline{n}(\zeta>t)}{\mathbf{P}_x(\tau_0^->t)} \cdot    \underline{n}\big(G(\epsilon^{(t)})\,\big|\, \zeta>t\big) \to \frac{ \bar{\mathtt{d}}}{\widehat{V}(x)} \cdot \underline{n}^Y\big(G(\epsilon)\,|\,\zeta>1\big). 
 \eeqlb
 For the integral in \eqref{eqn.5001}, we can write it as the summation of the next two terms:
 \beqnn
   \widetilde{I}_{K}(t) \ar:=\ar \int_0^K ds \int_{\mathcal{E}}  \overline{n}(d\epsilon,  \epsilon_s\leq x, \zeta>s)  \int_{\mathcal{E}}
 G\Big( \mathbf{1}_{\{ r\geq s/t \}}\cdot  \epsilon^{*,(t)}_{r-s/t} ; r\in[0,1]\Big) \, \underline{n}(d\epsilon^*, \zeta>t-s),\cr
  \varepsilon_{K}(t)  \ar:=\ar  \int_0^K ds \int_{\mathcal{E}}   \overline{n}(d\epsilon, \epsilon_s\leq x, \zeta>s)  \int_{\mathcal{E}}\varepsilon^G_s(\epsilon,\epsilon^*)\, 
 \, \underline{n}(d\epsilon^*, \zeta>t-s).
 \eeqnn
 where $\varepsilon^{G,(t)}_s(\epsilon,\epsilon^*):= G\big( \widetilde{\epsilon}_r^{(t),s}\big)-G\big( \mathbf{1}_{\{ r\geq s/t \}}\cdot  \epsilon^{*,(t)}_{r-s/t} ; r\in[0,1]\big)$.
  Since the integrand function in $\widetilde{I}_{K}(t)$ depends only on $\epsilon^*$, we have  
 \beqnn
 \widetilde{I}_{K}(t)= \int_0^K  \overline{n}( \epsilon_s\leq x, \zeta>s) \cdot \underline{n}( \zeta>t-s) \cdot \underline{n}\Big( G\Big( \mathbf{1}_{\{ r\geq s/t \}}\cdot  \epsilon^{*,(t)}_{r-s/t} ; r\in[0,1]\Big)\,\Big|\, \zeta>t-s \Big) \,  ds    .
 \eeqnn 
 Two direct consequences of   \eqref{eqn.4009} and Theorem~\ref{MainThm01} are that as $t\to\infty$,
 \beqnn
 \frac{\underline{n}( \zeta>t-s)}{\mathbf{P}_x\big(\tau_0^-> t \big)}\sim  \frac{1 }{\widehat{V}(x)} 
 \quad \mbox{and}\quad 
 \underline{n}\Big( G\Big( \mathbf{1}_{\{ r\geq s/t \}}\cdot  \epsilon^{*,(t)}_{r-s/t} ; r\in[0,1]\Big)\,\Big|\, \zeta>t-s \Big) \to  \underline{n}^Y \big(  G(\epsilon)\,|\,\zeta >1 \big),
 \eeqnn
 uniformly in $s\in[0,K]$. 
 These yield that
 \beqlb\label{eqn.5002}
 \lim_{t\to\infty} \frac{ \widetilde{I}_{K}(t)}{\mathbf{P}_x\big(\tau_0^-> t \big)} = \frac{\underline{n}^Y\big(G(\epsilon)\,|\,\zeta>1\big)}{\widehat{V}(x)} \cdot \int_0^K \overline{n}(  \epsilon_s\leq  x, \zeta>s) ds.
 \eeqlb
 We now turn to consider $ \varepsilon_{K}(t)$. 
 By \eqref{eqn.447} and $\underline{n}(  \zeta>t-s)/\mathbf{P}_x\big( \tau_0^->t \big) \to 1/\widehat{V}(x)$ uniformly in $s\in[0,K]$; see \eqref{eqn.4009}, we have uniformly in $t>0$, 
 \beqnn
 \int_0^K ds \int_{\mathcal{E}}  \overline{n}\big(d\epsilon,  \epsilon_s\leq x, \zeta>s\big)  \int_{\mathcal{E}}
 \, \frac{\underline{n}(d\epsilon^*, \zeta>t-s)}{\mathbf{P}_x\big( \tau_0^->t \big)}
 \ar=\ar\int_0^K    \overline{n}\big(   \epsilon_s\leq x, \zeta>s\big) \cdot  \frac{\underline{n}(  \zeta>t-s)}{\mathbf{P}_x\big( \tau_0^->t \big)}\,ds<\infty. 
 \eeqnn  
 This allows us to apply the dominated convergence theorem to obtain that 
 \beqnn
 \limsup_{t\to\infty} \frac{\big| \varepsilon_{K}(t) \big|}{\mathbf{P}_x\big( \tau_0^->t \big)} 
 \ar\leq\ar 
 \int_0^K ds \int_{\mathcal{E}}   \overline{n}(d\epsilon, \epsilon_s\leq x, \zeta>s)  \limsup_{t\to\infty} \int_{\mathcal{E}}\big|\varepsilon^{G,(t)}_s(\epsilon,\epsilon^*)\big|\, 
 \, \frac{\underline{n}(d\epsilon^*, \zeta>t-s) }{\mathbf{P}_x\big( \tau_0^->t \big)},
 \eeqnn
 which equals to $0$, since $G$ is uniformly continuous and hence $\sup_{\epsilon^* \in\mathcal{E}}\big|\varepsilon^{G,(t)}_s(\epsilon,\epsilon^*)\big| \to 0$ as $t\to\infty$.  
 Taking this, \eqref{eqn.5002}, \eqref{eqn.339} firstly back into \eqref{eqn.5001}  and then back into \eqref{eqn.340},  we obtain that as $t\to\infty$ and then $K\to\infty$
 \beqnn
  \mathbf{E}_x\Big[ G\big( X^{(t)}\big); \underline{g}_t\leq K\,\big|\, \tau_0^->t \Big]
 	\sim  \frac{I_{K}(t)}{ \mathbf{P}_x\big(\tau_0^-> t \big)}
 \to \frac{\underline{n}^Y\big(G(\epsilon)\,|\,\zeta>1\big)}{\widehat{V}(x)} \cdot \Big( \bar{\tt d} +  \int_0^\infty  \overline{n}(  \epsilon_s\leq x, \zeta>s) ds \Big) ,
 \eeqnn
  which equals to $\underline{n}^Y\big(G(\epsilon)\,|\,\zeta>1\big) $ because of the identity \eqref{eqn.447}. 
  The proof ends. 
 \qed

 Before proving Theorem~\ref{MainThm03}, we need to recall several well-known properties of spectrally positive L\'evy processes.  
 By (2) in \cite[p.191]{Bertoin1996}, the absence of negative jumps ensures that
 \beqlb
 \widehat{\kappa}(0,1)= \widehat{\kappa}^Y(0,1)=1
 \quad \mbox{and}\quad 
 \widehat{V}(x)=\widehat{V}^Y(x)=x,\quad x\geq 0.
 \eeqlb  
 For all $\lambda \geq 0$, let $\psi(\lambda):=\mathbf{E}[\exp\{ -\lambda X_1 \}]$ be the Laplace exponent of $X$, which is zero at zero and increases strictly to infinity at infinity. It is strictly convex and infinitely differentiable on $(0,\infty)$. 
 The \textit{scale function} $W=\{ W(x): x\in\mathbb{R} \}$ associated to $\psi$ is a non-negative function, which is identically zero on $(-\infty,0)$ and characterized on $[0,\infty)$ as a continuous, strictly increasing function whose Laplace transform is given by 
 \beqlb\label{eqn.ScaleF}
   \int_0^\infty e^{-\lambda x}W(x)\, dx = \frac{1}{\psi(\lambda)},\quad \lambda >0. 
 \eeqlb 
 Recall the identity (8.22) in Kyprianou \cite[p.219]{Kyprianou2014} that provides an alternate presentation for $W$:
 \beqlb\label{eqn.500}
 W(x)=  W(1)\exp\bigg\{ \int_1^x \underline{n}(\bar\epsilon>r)\, dr \bigg\} ,\quad x\geq 0.
 \eeqlb 
 Similarly, the scale function of $Y$ has an exact expression $W^Y (x)= x^{\alpha-1}/\Gamma(\alpha)$ for $x\geq 0$.
 This along with the preceding identity  deduces that 
 \beqlb\label{eqn.446}
 \underline{n}^Y(\overline{\epsilon}>x) = \frac{\alpha-1}{x},\quad x>0. 
 \eeqlb

 \textit{\textbf{Proof of Theorem~\ref{MainThm03}.}} 
 For the first claim, applying integration by parts to \eqref{eqn.ScaleF} shows that  
 \beqnn
 \int_0^\infty  e^{-\lambda x}\, dW(x) 
 =  \int_0^\infty \lambda e^{-\lambda x} \big(W(x)-W(0)\big)\, dx 
 =\frac{\lambda }{\psi(\lambda )}-W(0),
 \eeqnn
 which goes to $\infty$ as $\lambda \to 0$ and also belongs to $\mathrm{RV}^0_{1-\alpha}$.
 By the Tauberian theorem; see \cite[p.10]{Bertoin1996}, we have $W \in \mathrm{RV}^\infty_{\alpha-1}$.  
 By the identity \ref{eqn.500} and the Karamata representation theorem; see Theorem~1.3.1 in \cite[p.12]{BinghamGoldieTeugels1987},  we have  as $t\to\infty$,
 \beqnn
 \underline{n}\big(\overline{\epsilon}>t\big) \sim \frac{\alpha-1}{t}.
 \eeqnn
 This along with \eqref{eqn.4002} induces the asymptotic equivalence \eqref{eqn.102} immediately.

 We now turn to prove the second claim, which follows if for any bounded functional $F$ on $\mathcal{E}$,
 \beqlb\label{eqn.321}
 \underline{n}\big(F(\epsilon^{(t)})\,|\,\bar{\epsilon}>\boldsymbol{c}(t)\big) \to  \underline{n}^Y \big(F(\epsilon)\,|\,\bar{\epsilon}>1\big),
 \eeqlb
 as $t\to\infty$. 
 For any $t>1$ and $\delta\in(0,1)$, we can split the left-hand side into the following two terms:
 \beqlb\label{eqn.348}
 \underline{n}\big(F(\epsilon^{(t)}), \zeta >\delta t\,|\,\overline{\epsilon}>\boldsymbol{c}(t)\big)
 \quad \mbox{and}\quad 
 \underline{n}\big(F(\epsilon^{(t)}),\zeta \leq \delta t\,|\,\overline{\epsilon}>\boldsymbol{c}(t)\big).
 \eeqlb
 Firstly, we rewrite the first term as
\beqlb\label{eqn.311}
 \underline{n}\big(F(\epsilon^{(t)}), \zeta >\delta t\,|\,\overline{\epsilon}>\boldsymbol{c}(t)\big) 
 \ar=\ar \frac{ \underline{n}\big(F(\epsilon^{(t)}), \zeta >\delta t,\overline{\epsilon}>\boldsymbol{c}(t)\big)}{ \underline{n}\big(\overline{\epsilon}>\boldsymbol{c}(t)\big)} \cr
 \ar=\ar  \underline{n}\bigg(F(\epsilon^{(t)}), \sup_{s\geq 0}\epsilon_{ts}>\boldsymbol{c}(t)\,\Big|\, \zeta >\delta t\bigg) \cdot \frac{\underline{n}\big(\zeta >\delta  t \big)}{\underline{n}\big(\overline{\epsilon}>\boldsymbol{c}(t)\big)}\cr
 \ar=\ar  \underline{n}\bigg(F(\epsilon^{(t)}), \sup_{s\geq 0}\epsilon^{(t)}_{s}>1\,\Big|\, \zeta >\delta t\bigg) \cdot \frac{\underline{n}\big(\zeta >\delta  t \big)}{\underline{n}\big(\overline{\epsilon}>\boldsymbol{c}(t)\big)}.
\eeqlb
Note that $\rho=1/\alpha$ and $\epsilon^{(t)}_s =\frac{\boldsymbol{c}(\delta t)}{\boldsymbol{c}(t)}\cdot  \epsilon^{(\delta t)}_{s/\delta}$ for all $s\geq 0$.
By Theorem~\ref{MainThm01} and $\boldsymbol{c} \in \mathrm{RV}^\infty_{1/\alpha}$, we have $\boldsymbol{c}(\delta t)/\boldsymbol{c}( t)\to \delta^{\rho}$ and 
\beqnn
 \lim_{t\to\infty} \underline{n}\bigg(F(\epsilon^{(t)}), \sup_{0\leq s<\zeta/t}\epsilon^{(t)}_{s}>1\,\Big|\, \zeta >\delta t\bigg)   
\ar=\ar \underline{n}^Y\bigg( F(\delta^{\rho}\cdot \epsilon_{\cdot/\delta}), \sup_{s\geq 0}\,\delta^{\rho}\cdot \epsilon_{s/\delta}>1\,\Big|\, \zeta >1 \bigg),
\eeqnn
 which equals to $ \underline{n}^Y \big( F( \epsilon),   \overline{\epsilon}>1\,|\, \zeta >\delta \big)$ because of the scaling property of stable excursions; see Lemma~14 in \cite[p.232]{Bertoin1996}. 
Moreover, the first claim \eqref{eqn.102} induces that as $t\to\infty$, 
\beqnn
\frac{\underline{n}\big(\zeta >\delta  t \big)}{ \underline{n}\big(\overline{\epsilon}>\boldsymbol{c}(t)\big)} \to \frac{\delta^{-\rho}}{(\alpha -1)\cdot \Gamma(1-\rho) }. 
\eeqnn
Plugging these two estimates back into the right-hand side of the last equality in \eqref{eqn.311} and then using the identity $\underline{n}^Y\big(  \zeta >\delta \big)= \delta^{-\rho}/\Gamma(1-\rho)$; see \eqref{eqn.4006},   
\beqnn
 \lim_{t\to\infty} \underline{n}\big(F(\epsilon^{(t)}), \zeta >\delta t\,\big|\,\overline{\epsilon}>\boldsymbol{c}(t)\big) 
\ar\to \ar \frac{\delta^{-\rho}\cdot \underline{n}^Y \big( F( \epsilon),   \overline{\epsilon}>1\,\big|\, \zeta >\delta \big)}{(\alpha -1)\cdot \Gamma(1-\rho) }
=  \frac{ \underline{n}^Y \big( F( \epsilon),   \overline{\epsilon}>1,  \zeta >\delta \big)}{\alpha -1 }  .
\eeqnn
By using the monotone convergence theorem and then \eqref{eqn.446} with $x=1$, we have 
\beqlb\label{eqn.317}
 \lim_{\delta \to 0+} \lim_{t\to\infty} \underline{n}\big(F(\epsilon^{(t)}), \zeta >\delta t\,\big|\,\overline{\epsilon}>\boldsymbol{c}(t)\big)  
 = \frac{ \underline{n}^Y \big( F( \epsilon),   \overline{\epsilon}>1  \big)}{\alpha -1 }
 =  \underline{n}^Y \big( F( \epsilon) \,|\,  \overline{\epsilon}>1\big) . 
\eeqlb
For the second term in \eqref{eqn.348}, the upper bound of $F$ yields that  uniformly in $t> 1$ and $\delta \in(0,1)$, 
\beqnn
\underline{n}\big(F(\epsilon^{(t)}),\zeta \leq \delta t\,|\,\overline{\epsilon}>\boldsymbol{c}(t)\big)
\leq C\cdot \underline{n}\big(\zeta \leq \delta t\,|\,\overline{\epsilon}>\boldsymbol{c}(t)\big)
\ar=\ar C\cdot \big(1-\underline{n}\big(\zeta >\delta t\,|\, \overline{\epsilon}>\boldsymbol{c}(t)\big) \big),
\eeqnn
which vanishes as $\delta  \to 0+$ because of \eqref{eqn.317} with $F\equiv1 $.  
 Combining this together with \eqref{eqn.317} induces the desired limit \eqref{eqn.321}. 
 The whole proof ends. 
 \qed

  \section{The negative-drift case} 
 \label{Sec.TransientCase}
 \setcounter{equation}{0}
 
 This section is devoted to proving Theorem~\ref{MainThm04}, \ref{MainThm07} and Corollary~\ref{MainThm08}. 
 Since $0$ is assumed to be regular for $(0, \infty)$, for $X$, we have $\underline{\mathrm{d}}=0$.
 Applying the change of variables to \eqref{eqn.kappa} with $u=0$ gives that 
 \beqlb\label{eqn.551}
 	\widehat \kappa(\lambda,0) =  \int_0^\infty (1- e^{-\lambda t})\,\underline{n}(\zeta \in dt)
 	= \lambda \int_0^\infty  e^{-\lambda t}\underline{n}(\zeta > t)\, dt,\quad \lambda \geq 0.
 \eeqlb
 By Theorem 6.9 (ii) in \cite{Kyprianou2014},  we have that $\kappa(0,0) >0$ since $X$ drift to  $-\infty$ and hence 
 \beqlb\label{eqn.501}
 \underline{n}(\zeta )= \int_0^\infty  \underline{n}(\zeta > t)dt=\lim_{\lambda\to 0}\frac{\widehat \kappa(\lambda,0)}{\lambda} =\frac{1}{ \kappa(0,0)} < \infty.
 \eeqlb
 Here the last equality follows from the identity $\widehat \kappa(\lambda,0) \cdot \kappa(\lambda,0)=\lambda$ for all $\lambda \geq 0$. 
 Additionally, for any $x>0$, by Proposition 17 (ii) in \cite[p.172]{Bertoin1996} and Theorem 2.2  in \cite{DenisovShneer2013} we have 
 \beqlb  \label{Asytau} 
    \mathbf{E}_{x}[\tau_0^-]=\frac{\widehat{V}(x)}{\kappa(0,0)} = \underline{n}(\zeta ) \cdot \widehat{V}(x)
    \quad \mbox{and}\quad 
    \mathbf{P}_x\big(\tau_0^->t\big)  
    \sim  \underline{n}(\zeta ) \cdot \widehat{V}(x) \cdot \overline{\nu}(\beta t) \in \mathrm{RV}^\infty_{-\theta} . 
 \eeqlb
 For $x>0$, let $J^x:=\inf\{s> 0: \Delta \omega_s>x  \}$ be the arrival time of the first jump of $\omega$ larger than $x$ and
 $N^{x}_t :=\# \{ s\in(0,t]: \Delta \omega_s>x\}$ be the number of jumps larger than $x$ before time $t$.
 The compensation formula \eqref{eqn.Compensation} deduces that for all $ T,x\geq 0$, 
 \beqlb\label{eqn.502}
 \underline{n}\big(J^x\leq T\wedge \zeta  \big)
 \leq \underline{n}\big(N^x_{T\wedge \zeta}  \big)
 = \underline{n}\bigg(\sum_{0\leq s\leq T\wedge\zeta }\mathbf{1}_{\{\Delta \epsilon_s>x\}}   \bigg)=\underline{n}(T\wedge\zeta) \cdot \overline\nu(x).
 \eeqlb
 We now start to prove Theorem~\ref{MainThm03} with the help of the next two propositions, which establish a uniform upper bound for the first passage time $\tau_0^-$ and a uniform linear approximation for $X$. 
   \begin{proposition}\label{Prop.406}
 	For any $0<\delta<\beta $, there exists a constant $C>0$ such that for all  $t>0 $ and $x \in[0, \delta t]$,
 	\beqlb
 	 \mathbf{P}_x\big(\tau_0^- >t\big)\leq C \cdot (1+x)\cdot  \bar{\nu}(\beta t).
 	\eeqlb
 \end{proposition}
 \proof
 Let $\{\tau_{-1,i}^{-}\}_{i \geq 1}$ be a sequence of i.i.d. copies of $\tau_{-1}^{-}$. 
 By the strong Markov property of $X$,
 \beqnn
  \mathbf{P}_x\big(\tau_{0}^{-} > t\big) 
  =\mathbf{P}\big(\tau_{-x}^{-} > t\big)   
  \ar\leq\ar \mathbf{P}\bigg(\sum_{i=1}^{[x]+1} \tau_{-1,i}^{-}>t\bigg)  . 
 \eeqnn
 For any constant $\gamma>1$ such that $\gamma \cdot \delta \cdot \mathbf{E}[\tau_{-1}^{-}]<1$, we have for all $t>0$ and $x \in[0, \delta t]$, 
 \beqnn
 \mathbf{P}\bigg(\sum_{i=1}^{[x]+1} \tau_{-1,i}^{-} >t\bigg) 
 \ar=\ar \mathbf{P}\bigg(\sum_{i=1}^{[x]+1}\big(\tau_{-1,i}^{-}-\gamma \cdot  \mathbf{E}[\tau_{-1}^{-}]\big) > t- ([x]+1) \cdot \gamma \cdot  \mathbf{E}[\tau_{-1}^{-}]  \bigg) \cr
 \ar\leq\ar  \mathbf{P}\bigg(\sum_{i=1}^{[x]+1}\big(\tau_{-1,i}^{-}- \gamma \cdot  \mathbf{E}[\tau_{-1}^{-}]\big) > t- (\delta t+1) \cdot \gamma \cdot  \mathbf{E}[\tau_{-1}^{-}]  \bigg).
 \eeqnn 
 Since $t- (\delta t+1) \cdot \gamma \cdot  \mathbf{E}[\tau_{-1}^{-}] \sim \big(1- \delta \cdot \gamma\cdot  \mathbf{E}[\tau_{-1}^{-}] \big)\cdot t \to \infty $, 
 by Theorem 2 in \cite{DenisovFossKorshunov2010} there exists a constant $C>0$ such that for large $t$, the last probability is bounded uniformly in $x \geq 0$  by 
 \beqnn
 C\cdot (x+1) \cdot \mathbf{P}\big(\tau_{-1}^{-} \geq \big(1- \delta \cdot  \gamma\cdot  \mathbf{E}[\tau_{-1}^{-}] \big)\cdot t \big)  ,
 \eeqnn
 which, by \eqref{Asytau} and Assumption~\ref{Assumption02}, is asymptotically equivalent as $t\to\infty$ to 
 \beqnn
  \frac{C}{\beta}\cdot (x+1) \cdot \overline{\nu}\big( \big(1- \delta \cdot \theta\cdot  \mathbf{E}[\tau_{-1}^{-}] \big)\cdot t \big)
  \sim  \frac{C}{\beta}\cdot (x+1) \cdot  \Big(\frac{1-\delta \gamma /\beta }{\beta}\Big)^{-\theta} \cdot \bar{\nu}(\beta t).
  \eeqnn
 Consequently, the desired uniformly upper bound holds.
 \qed
 
 \begin{proposition}\label{Prop.402}
 	For any $T\geq 0$ , we have 
 	$\sup_{s\in[0,T]} \big| X_{ts}/t+\beta s\big|  \to 0$ in probability as $t\to\infty$.  
 \end{proposition}
 \proof By the L\'evy-It\^o decomposition theorem, the process $X$ admits the following representation
 \beqlb\label{eqn.445}
 X_t:= - \beta t + \sigma B_t + \int_0^t \int_\mathbb{R} y \widetilde{N}(dr,dy) ,\quad t\geq 0,
 \eeqlb
 for some constant $\sigma\geq 0$ and compensated Poisson random measure $\widetilde{N}(ds,dy)$ on $(0,\infty)\times \mathbb{R}$ with intensity $ds\,\nu(dy)$. 
 For any $K>0$, we can decompose $X_t+\beta t$ into the following two terms:
 \beqnn
 M^{\leq K}_t:=  \sigma B_t + \int_0^t \int_{|y|\leq K} y \, \widetilde{N}(dr,dy)
 \quad \mbox{and}\quad 
 M^{> K}_t:=   \int_0^t \int_{|y|> K} y\, \widetilde{N}(dr,dy). 
 \eeqnn
 Hence it suffices to prove that for any $\varepsilon>0$, 
 \beqnn
 \lim_{K\to\infty} \limsup_{t\to\infty} \mathbf{P}\Big( \sup_{s\in[0,Tt]} \big| M^{> K}_s\big| \geq \varepsilon t
 \Big) =0
 \quad \mbox{and}\quad 
 \lim_{t\to\infty} \mathbf{P}\Big( \sup_{s\in[0,Tt]} \big| M^{\leq K}_s\big| \geq \varepsilon t
 \Big) =0
 \eeqnn
 Firstly, by Chebyshev's inequality and the fact $|M^{> K}_s|\leq   \int_0^s \int_{|y|> K} |y| N(dr,dy) + s\cdot  \int_{|y|> K} |y| \, \nu(dy) $, 
 \beqnn
 \sup_{t> 0}\mathbf{P}\Big( \sup_{s\in[0,Tt]} \big| M^{> K}_s\big| \geq \varepsilon t
 \Big) 
 \ar\leq\ar \sup_{t> 0} \frac{1}{ \varepsilon t}\cdot \mathbf{E}\bigg[ \sup_{s\in[0,Tt]} \big| M^{> K}_s\big|
 \bigg]
 \leq \frac{2T}{ \varepsilon }\cdot \int_{|y|> K} |y| \, \nu(dy), 
 \eeqnn
 which goes to $0$ as $K\to\infty$. 
 Moreover, by using Chebyshev's inequality again and then the Burkholder-Davis-Gundy inequalities to the first probability, there exists some constant $C>0$ independent of $t$, $T$ and $K$ such that 
 \beqnn
 \mathbf{P}\Big( \sup_{s\in[0,Tt]} \big| M^{\leq K}_s\big| \geq \varepsilon t
 \Big)  
 \leq \frac{1}{|\varepsilon t|^2} \mathbf{E}\bigg[ \sup_{s\in[0,Tt]} \big| M^{\leq K}_s\big|^2
 \bigg]
 \ar\leq\ar  \frac{C\cdot T}{|\varepsilon |^2\cdot t} \cdot \Big[\sigma^2 + \int_{|y|\leq K} |y|^2 \nu(dy)\Big],
 \eeqnn
 which also vanishes as $t\to \infty$. 
 \qed
 	
 \begin{lemma}\label{Lemma.401}
  As $t\to\infty$, we have $\underline{n}\big(J^{\beta t}\leq T,\zeta>t \big) \sim \underline{n}(T\wedge \zeta) \cdot \overline{\nu}(\beta t)$ uniformly in $T\geq0$; equivalently, for all $T_0\geq 0$,
  \beqlb\label{eqn.401}
  \limsup_{t\to\infty} \sup_{T\in[0,T_0]} \Big| \frac{\underline{n}\big(J^{\beta t}\leq T,\zeta>t \big)}{ \overline{\nu}(\beta t)} - \underline{n}(T\wedge \zeta) \Big|
  +   \lim_{T\to\infty} \limsup_{t\to\infty} \frac{\underline{n}\big(J^{\beta t}> T,\zeta>t \big)}{ \overline{\nu}(\beta t)}
   =0. 
  \eeqlb
 \end{lemma}
 \proof The equivalence of these two claims can be easily identified by using the fact $ \underline{n}(\zeta)<\infty$;  see \eqref{eqn.501}. 
 We now prove the two limits in \eqref{eqn.401} separately in the next two steps.
 
 {\textbf{Step 1.}} For the first limit, it is easy to see that for all $t\geq T_0$,
 \beqlb\label{eqn.503}
 \underline{ n}\big(J^{\beta t} \leq T,\zeta >t\big) \leq \underline{n}\big(J^{\beta t} \leq T\wedge \zeta \big) = \underline{n}\big(N^{\beta t}_{T\wedge \zeta}\geq 1 \big)\leq \underline{n}\big(N^{\beta t}_{T\wedge \zeta}  \big) .
 \eeqlb
 This along with \eqref{eqn.502} induces that uniformly in $T\in[0,T_0]$, 
 \beqnn
 \frac{\underline{n}\big(J^{\beta t} \leq T,\zeta >t\big) }{\bar{\nu}(\beta t)}
 \leq  \underline{n}(T\wedge\zeta)   . 
 \eeqnn
 It remains to prove that the asymptotic lower bound also holds, i.e., uniformly in $T\in [0,T_0]$,
 \beqlb\label{eqn.402}
  \liminf_{t\to\infty} \frac{\underline{n}\big(J^{\beta t} \leq T,\zeta >t\big) }{\bar{\nu}(\beta t)}
  \geq \underline{n}(T\wedge\zeta)  . 
 \eeqlb
 For any $z>\beta$, we see that $ \underline{n}\big(J^{\beta t} \leq T,\zeta >t\big)\geq \underline{n}\big(J^{z t} \leq T, \zeta>t \big)$.
 By the Markov property,  
 \beqnn
 \underline{n}\big(J^{z t} \leq T, \zeta>t \big)
 \ar=\ar \int_{0}^{T}\int_{zt}^{\infty}    \underline{n}\big(J^{z t} \in ds, \epsilon_s\in dx, \zeta>t\big)\cr
  \ar=\ar \int_{0}^{T}\int_{zt}^{\infty}\mathbf{P}_{x}\big( \underline{X}_{t-s}>0\big) \, \underline{n}\big(J^{z t}\in ds,\epsilon_s\in dx,\zeta >s\big),
 \eeqnn  
 which, by the monotone of $\underline{X}$, can be bounded from below by $\mathbf{P}_{zt}\big( \underline{X}_{t}>0\big)\cdot \underline{n}\big(J^{z t}\leq T\wedge \zeta \big)$. 
 Since $z>\beta=-\mathbf{E}[X_1]$, the L\'evy process $W:=\{W_t:=\frac{z+\beta}{2}\cdot t+ X_t, t\geq 0 \}$ drifts to infinity  and hence  as $t\to\infty$,
 \beqnn
 \mathbf{P}_{zt}\big( \underline{X}_{t}>0\big)
 = \mathbf{P}\big( \underline{X}_{t} +zt>0\big) =  \mathbf{P}\bigg(  \underline{X}_t +\frac{z+\beta}{2}\cdot t >-\frac{z-\beta}{2}\cdot t \bigg) 
 \geq  \mathbf{P}\bigg(  \underline{W}_t >-\frac{z-\beta}{2}\cdot t \bigg) 
 \to 1.
 \eeqnn
 Consequently, the lower bound \eqref{eqn.402} follows if we can identify that uniformly in $T\in[0,T_0]$,
 \beqlb\label{eqn.421}
 \lim_{z\to \beta+} \liminf_{t\to\infty}
 \frac{\underline{n}\big(J^{z t}\leq T\wedge \zeta \big)}{\overline{\nu}(\beta t)} \geq \underline{n}\big(T\wedge \zeta\big).
 \eeqlb
 To prove this inequality, we first observe that $\{J^{z t}\leq T\wedge \zeta \} = \{N^{z t}_{T\wedge \zeta} \geq 1\}$, which implies that
 \beqnn
 \underline{n}\big(J^{z t}\leq T\wedge \zeta \big) 
 = \underline{n}\big(  N^{z t}_{T\wedge \zeta} \big)-\underline{n}\big( N^{z t}_{T\wedge \zeta}-1 ;N_{T\wedge\zeta}^{z t}\geq 2 \big) .
 \eeqnn
 By \eqref{eqn.502} and Assumption~\ref{Assumption02}, we have 
 $ \underline{n}\big(  N^{z t}_{T\wedge \zeta} \big) \sim (z/\beta)^{-\theta}\cdot \underline{n}\big(T\wedge \zeta\big) \cdot \overline{\nu}(\beta t) $ 
 as $t\to\infty$, 
 and then  uniformly in $T\in[0,T_0]$,
 \beqlb\label{eqn.414}
  \lim_{z\to \beta+} \liminf_{t\to\infty}
  \frac{\underline{n}\big(J^{z t}\leq T\wedge \zeta \big)}{\overline{\nu}(\beta t)} \geq \underline{n}\big(T\wedge \zeta\big) - \lim_{z\to \beta+} \limsup_{t\to\infty} \sup_{T\in[0,T_0]}
  \frac{ \underline{n}\big( N^{z t}_{T\wedge \zeta}-1 ;N_{T\wedge\zeta}^{z t}\geq 2 \big) }{\overline{\nu}(\beta t)}.
 \eeqlb
 Hence the lower bound \eqref{eqn.421} follows if the last limit equals to $0$. 
 Summation by parts deduces that 
 \beqlb\label{eqn.404}
  \underline{n}\big( N^{z t}_{T\wedge \zeta} -1;N_{T\wedge\zeta}^{z t}\geq 2 \big) 
  = \sum_{k=2}^\infty \underline{n}\big( N_{T\wedge\zeta}^{z t}\geq k \big).
 \eeqlb
Using the Markov property and then the independent increments to summand, 
 \beqnn
 \underline{n}\big( N_{T\wedge\zeta}^{z t}\geq k \big)
 \ar=\ar \int_{0}^{T}\int_{zt}^{\infty}\underline{n}\big( N_{T\wedge\zeta}^{z t}\geq k,  J^{z t}\in ds,\epsilon_s\in dy \big)\cr
 \ar=\ar 
 \int_{0}^{T}\int_{zt}^{\infty}\mathbf{P}_{y}\Big( N_{(T-s)\wedge\tau_0^-}^{zt}\geq k-1 \Big)\, \underline{n}\big(  J^{z t}\in ds,\epsilon_s\in dy, \zeta>s\big) , 
 \eeqnn
 which can be bounded by $ \mathbf{P}\big( N_{T}^{zt}\geq k-1 \big)\cdot  \underline{n}\big( J^{z t}\leq T \wedge \zeta\big)$ since $ N_s^{zt}$  is  non-decreasing in $s$. 
 By \eqref{eqn.502}, \eqref{eqn.503} and the fact that $N^{zt}_T$ is Poisson distributed with rate $T\cdot \overline{\nu}(zt)$,
 \beqnn
  \underline{n}\big( N_{T\wedge\zeta}^{z t}\geq k \big) 
  \leq \mathbf{P}\big( N_{T}^{zt}\geq k-1 \big)\cdot  \underline{n}(T\wedge\zeta) \cdot \bar{\nu}(z t)
  = \underline{n}\big( T\wedge\zeta \big) \cdot \overline{\nu}(z t)\cdot \sum_{j=k-1}^{\infty}\frac{(T\cdot \overline{\nu}(z t))^{j}}{j!} \cdot e^{-T\cdot \overline{\nu}(z t) } .
 \eeqnn 
 Plugging this back into the sum in \eqref{eqn.404}, then using the change of variables and summation by parts,  
 \beqlb\label{eqn.42}
 \underline{n}\big( N^{z t}_{T\wedge \zeta} -1 ;N_{T\wedge\zeta}^{z t}\geq 2 \big) 
  \ar\leq\ar \underline{n}\big( T\wedge\zeta \big) \cdot \overline{\nu}(z t)\cdot \sum_{k=2}^\infty  \sum_{j=k-1}^{\infty}\frac{(T\cdot \overline{\nu}(z t))^{j}}{j!} \cdot e^{-T\cdot \overline{\nu}(z t) }\cr
  \ar=\ar  \underline{n}\big( T \wedge\zeta \big) \cdot \overline{\nu}(z t)\cdot \sum_{k=1}^\infty  \sum_{j=k}^{\infty}\frac{(T\cdot \overline{\nu}(z t))^{j}}{j!} \cdot e^{-T \cdot \overline{\nu}(z t) } \cr
  \ar=\ar \underline{n}\big( T \wedge\zeta \big) \cdot \overline{\nu}(z t)\cdot \sum_{k=1}^\infty k\cdot \frac{(T\cdot \overline{\nu}(z t))^k}{k!} \cdot e^{-T\cdot\overline{\nu}(z t) }, 
 \eeqlb
 which equals to $\underline{n}\big( T\wedge\zeta \big) \cdot T\cdot \big|\bar{\nu}(z t)\big|^2$ and hence the last term in \eqref{eqn.414} equals to $0$. 
 Consequently, the lower bound \eqref{eqn.402} holds and then the first limit in \eqref{eqn.401} follows. 

 {\textbf{Step 2.}}  For the second limit in \eqref{eqn.401}, 
 we continue to use the notation $\tau_x^+ =\inf\{ s>0: \epsilon_s>x \}$ and  then 
 decompose $\underline{n}\big(J^{\beta t}> T,\zeta>t \big)$ into the following two terms
 \beqlb\label{Ra114}
 \underline{n}\big(J^{\beta t} > T,\tau_{1}^{+}> T/2, \zeta >t\big)
 \quad \mbox{and}\quad 
 \underline{n}\big(J^{\beta t} > T,\tau_{1}^{+}\leq T/2, \zeta >t\big).
 \eeqlb 
An application of the Markov property to the first term at time  $T/2$ shows that 
\begin{align*}
 \underline{n}\big(J^{\beta t} > T,\tau_{1}^{+}> T/2, \zeta >t\big) = \int_{0}^{1}\mathbf{P}_x\big(J^{\beta t}>T/2,\tau_0^{-}>t-T/2 \big)\, \underline{ n} \big(\epsilon_{T/2}\in dx,\zeta\wedge \tau_1^+>T/2 \big),
\end{align*}
which is smaller than $\mathbf{P}_1\big(\tau_0^{-}>t-T/2 \big)
\cdot \underline{ n}\big(\epsilon_{T/2}\leq 1,\zeta\wedge \tau_1^+>T/2 \big)$. 
By Proposition~\ref{Prop.406},  
\beqlb\label{eqn.418}
	\limsup_{t\to \infty}\frac{1}{\bar{\nu}(\beta t)} \cdot \underline{n}\big(J^{\beta t} > T,\tau_{1}^{+}> T/2, \zeta >t\big)
	\leq C\cdot  \underline{n} \big(\epsilon_{T/2}\leq 1,T/2<\zeta\wedge \tau_1^+ \big),
\eeqlb
 which goes to $ 0 $ as $ T \to \infty $.
 
 We now consider the second term in \eqref{Ra114}. For any $ \delta\in(0, \beta)$, it can be further decomposed into the following two terms by the overshoot size at $ \tau_{1}^+ $,
 \beqlb\label{eqn.416}
 	\underline{n}\big( J^{\beta t} > T,\tau_{1}^{+}\leq  T/2, 1<\epsilon_{\tau_{1}^{+}}\leq  \delta t, \zeta >t\big)
   \quad \mbox{and}\quad 
	\underline{n}\big( J^{\beta t} > T,\tau_{1}^{+}\leq  T/2, \epsilon_{\tau_{1}^{+}}>\delta t, \zeta >t\big). 
\eeqlb
 For $t\geq 2T$, we use the Markov property to the first term and get that  it equals to
 \beqlb\label{eqn.412}
  \lefteqn{ \int_0^{T/2} \int_1^{\delta t} \underline{n}\big( \tau_{1}^{+}\in ds, \epsilon_s \in dx,  J^{\beta t} > T, \zeta >t\big)}\ar\ar\cr
  \ar=\ar \int_0^{T/2} \int_1^{\delta t} \mathbf{P}_x\big(J^{\beta t} > T-s, \tau_0^->t-s\big) \, \underline{n}\big( \tau_{1}^{+}\in ds, \epsilon_s \in dx, \zeta >s\big)\cr
  \ar\leq\ar \int_0^{T/2} \int_1^{\delta t} \mathbf{P}_x\big(J^{\beta t} > T/2, \tau_0^->t/2\big) \, \underline{n}\big( \tau_{1}^{+}\in ds, \epsilon_s \in dx, \zeta >s\big)\cr
  \ar\leq\ar \int_1^{\delta t} \mathbf{P}_x\big(J^{\beta t} > T/2, \tau_0^->t/2\big) \, \underline{n}\big(   \epsilon_{\tau_{1}^{+}} \in dx, \zeta >\tau_{1}^{+}\big). 
 \eeqlb
By Proposition~\ref{Prop.406} and Theorem~3.4 in \cite{Xu2021a}, we have uniformly in $t\geq 2T$ and $x\in [1,\delta t]$,
 \beqnn
  \frac{1 }{\overline{\nu}(\beta t)} \cdot \mathbf{P}_x\big(J^{\beta t} > T/2, \tau_0^->t/2\big)
  \ar=\ar \frac{\mathbf{P}_x\big(\tau_0^->t/2\big)}{\overline{\nu}(\beta t)} \cdot \mathbf{P}_x\big(J^{\beta t} > T/2\,\big|\, \tau_0^->t/2\big) \cr
  \ar\leq \ar C\cdot x \cdot \Big(1- \frac{\mathbf{E}_x\big[\tau_0^- \wedge \frac{T}{2}\big]}{\mathbf{E}_x\big[\tau_0^- \big]} \Big).
 \eeqnn 
 Taking this back into the last integral in \eqref{eqn.412},  we obtain that
 \beqlb\label{eqn.415}
  \lefteqn{\sup_{t\geq 2T}\frac{1}{\overline{\nu}(\beta t)} \cdot \underline{n}\big( J^{\beta t} > T,\tau_{1}^{+}\leq  T/2, 1<\epsilon_{\tau_{1}^{+}}\leq  \delta t, \zeta >t\big) }\quad \ar\ar\cr
  \ar\leq\ar C  \int_1^\infty x \cdot \Big(1- \frac{\mathbf{E}_x\big[\tau_0^- \wedge  \frac{T}{2}\big]}{\mathbf{E}_x\big[\tau_0^- \big]} \Big) \,   \underline{n}\big(   \epsilon_{\tau_{1}^{+}} \in dx, \zeta >\tau_{1}^{+}\big) .
 \eeqlb
By using integration by parts, we have
\beqnn
\int_{1}^{\infty} x \, \underline{n}(\epsilon_{\tau_1^+}\in dx, \zeta> \tau_1^+ )
\ar=\ar \int_{0}^{\infty}  \underline{n}(\epsilon_{\tau_1^+}>1\vee x, \zeta> \tau_1^+ )\, dx\cr
\ar\leq\ar 2 \cdot \underline{n}(\overline{\epsilon}>1 )+\int_2^{\infty}  \underline{n}( \epsilon_{\tau_1^+}> x, \zeta> \tau_1^+)\, dx .
\eeqnn
Note that $\{ \epsilon_{\tau_1^+}> x, \zeta> \tau_1^+\} \subset\{ \Delta \epsilon_{\tau_1^+}>x-1, \zeta> \tau_1^+\} 
\subset\{ N^{x-1}_\zeta\geq 1  \}$. 
By \eqref{eqn.502} with $T=\infty$, 
\beqnn
\int_{1}^{\infty} x \, \underline{n}(\epsilon_{\tau_1^+}\in dx, \zeta> \tau_1^+ )
\ar\leq\ar 2 \cdot \underline{n}(\overline{\epsilon}>1 )+ \int_{2}^{\infty}  \underline{n}\big( N^{x-1}_{\zeta}  \big)\, dx \cr
\ar=\ar 2 \cdot \underline{n}(\overline{\epsilon}>1 )+ \underline{n}(\zeta) \cdot  \int_{2}^{\infty} \bar{\nu}(x-1)dx <\infty.
\eeqnn
This allows us to apply the dominated convergence theorem to   \eqref{eqn.415} and obtain that
\beqlb\label{eqn.417}
 \lim_{T\to\infty} \sup_{t\geq 2T}\frac{1}{\overline{\nu}(\beta t)} \cdot \underline{n}\big( J^{\beta t} > T,\tau_{1}^{+}\leq  T/2, 1<\epsilon_{\tau_{1}^{+}}\leq  \delta t, \zeta >t\big) =0.
\eeqlb
We now analyze the second term in \eqref{eqn.416}. 
Firstly, our observation 
$$\{ J^{\beta t} > T,\tau_{1}^{+}\leq  T/2, \epsilon_{\tau_{1}^{+}}>\delta t, \zeta >t \} \subset  \{  \tau_{1}^{+}\leq    \zeta, \delta t -1 <\Delta\epsilon_{\tau_{1}^{+}}< \beta t  \},\quad \forall\, t\geq 2T$$  
tells that it can be bounded by
\beqnn
\underline{n} \big( \tau_{1}^{+}\leq   \zeta, \delta t -1 <\Delta\epsilon_{\tau_{1}^{+}}< \beta t  \big)
\ar\leq\ar  \underline{n}\bigg(\sum_{0< s\leq  \zeta }\mathbf{1}_{\{\delta t-1 \leq \Delta\epsilon_s\leq \beta t\}}\bigg)
= \underline{n}\big(  N^{\delta t-1}_{ \zeta} \big)  -  \underline{n}\big(  N^{\beta t}_{\zeta} \big) ,
\eeqnn
which, by \eqref{eqn.502} with $T=\infty$ again, equals to $ \underline{n}(\zeta) \cdot \big(\overline{\nu}(\delta t-1)- \overline{\nu}(\beta t) \big)$. 
Consequently, we have
 \beqnn
 	\limsup_{t\to \infty}\frac{1}{\bar\nu(\beta t)} \cdot \underline{n}\big(\mathcal{J}^{\beta t} > T,\tau_{1}^{+}\leq  T/2, \epsilon_{\tau_{1}^{+}}\geq (\beta-\delta) t, \zeta >t\big)
 \leq   \underline{n}\big(\zeta \big) \cdot  \big[\big(\delta/\beta\big)^{-\theta}-1\big], 
 \eeqnn
 which goes to $0$ as $\delta \to \beta-$. 
 Put this together with \eqref{eqn.417}, we see that 
 \beqnn
 \lim_{t\to\infty} \frac{1}{\bar\nu(\beta t)} \cdot  \underline{n}\big(J^{\beta t} > T,\tau_{1}^{+}\leq T/2, \zeta >t\big) =0.
 \eeqnn
 This along with \eqref{eqn.418} induces that the second  limit in \eqref{eqn.401} equals to $0$. The whole proof ends. 
\qed

 \textit{\textbf{ Proof for Theorem~\ref{MainThm04}.}}
 A direct consequence of Lemma~\ref{Lemma.401} and \eqref{eqn.501} shows that
 \beqlb\label{eqn.5021}
	\lim_{t \to \infty}\frac{	\underline{n}\big(\zeta >t\big)}{\bar{\nu}(\beta t)}=\lim_{T\to \infty}\lim_{t \to \infty}\frac{	\underline{n}\big(J^{\beta t} \leq  T,\zeta >t\big)}{\bar{\nu}(\beta t)}+\lim_{T\to \infty}\lim_{t \to \infty}\frac{\underline{ n}\big(J^{\beta t} > T,\zeta >t\big)}{\bar{\nu}(\beta t)}
  = 	\underline{n}(\zeta)<\infty. 
 \eeqlb
 We now prove that $\underline{n}\big(\overline{\epsilon}>\beta t\big) \sim \underline{n}\big(\zeta >t\big)$ as $t\to\infty$.  
 Firstly, it is obvious that
 \beqlb\label{eqn.419}
  \underline{n}\big(\overline{\epsilon}>\beta t\big) \geq \underline{n}\big(\overline{\epsilon}>\beta t,\zeta >t\big) 
 = \underline{n}\big(\zeta >t\big) - \underline{n}\big(\overline{\epsilon}\leq \beta t,\zeta >t\big) 
 \eeqlb 
 Recall the notation $\tau_x^+ =\inf\{ s>0: \epsilon_s>x \}$. 
 Note that $\{\tau_{\beta t}^+>t, \zeta> t\}\subset \{J^{\beta t}>t,\zeta> t\}$, by the second limit in \eqref{eqn.401} we have as $t\to\infty$,
 \beqnn
  \underline{n}\big(\bar{\epsilon}\leq \beta t, \zeta>t \big) 
  \leq \underline{n}\big( \tau_{\beta t}^+>t, \zeta>t \big)
  \leq  \underline{n}\big( J^{\beta t}>t, \zeta>t \big) = o\big(   \underline{n}\big(  \zeta>t \big)\big).
 \eeqnn
 Taking this back into \eqref{eqn.419} gives that 
 \beqlb\label{eqn.422}
  \liminf_{t\to\infty} \frac{\underline{n}\big(\overline{\epsilon}>\beta t\big) }{ \underline{n}\big(  \zeta>t \big)} \geq 1.
 \eeqlb
 On the other hand,   we have 
 $\underline{n}\big(\overline{\epsilon}>\beta t\big) = \underline{n}\big(\overline{\epsilon}>\beta t, \zeta\leq  \delta t \big) + \underline{n}\big(\overline{\epsilon}>\beta t, \zeta >\delta t\big)$ for any $\delta\in(0,1)$. 
 By the Markov property, 
 \beqnn
 	\underline{n}\big(\overline{\epsilon}> \beta t,\zeta\leq \delta t \big)  =  \underline{n}\big( \tau_{\beta t}^+ <\zeta\leq\delta t\big)
 	\ar=\ar \int_0^{\delta t} \int_{\beta t}^\infty \underline{n}\big(  \tau_{\beta t}^+ \in ds, \epsilon_{s}\in dx , s<\zeta\leq\delta t\big)\cr
 	\ar=\ar  \int_0^{\delta t} \int_{\beta t}^\infty \mathbf{P}_x\big( \tau_0^-\leq\delta t-s \big)\, \underline{n}\big(  \tau_{\beta t}^+ \in ds, \epsilon_{s}\in dx ,  \zeta>s\big)\cr
 	\ar\leq\ar  \int_0^{\delta t} \mathbf{P}_{\beta t}\big( \tau_0^-\leq\delta t-s \big) \,  \underline{n}\big(  \tau_{\beta t}^+ \in ds,   \zeta>s\big),
 \eeqnn 
 which can be bounded by $
  \mathbf{P}_{\beta t}\big( \tau_0^-\leq\delta t \big) \cdot \underline{n}\big(  \overline{\epsilon}> \beta t\big)$.
  Since $\mathbf{E}[X_1]=-\beta$ and $\delta\in (0,1)$, the L\'evy process $\eta_t:= X_t +(2-\delta) \cdot \beta t$ drifts to infinity a.s.
  Note that $\underline{\eta}_t \leq \underline{X}_t +(2-\delta) \cdot \beta t$, we have that as $t\to\infty$,
  \beqnn
  \mathbf{P}_{\beta t}\big( \tau_0^-\leq\delta t \big) 
  \ar=\ar \mathbf{P}_{\beta t}\big(\underline{X}_{ \delta t}\leq0 \big) 
  =  \mathbf{P} \big(  \underline{X}_{ \delta t}\leq- \beta t\big) 
  \leq \mathbf{P} \big( \underline{\eta}_{ \delta t}\leq- (1-\delta )\beta t \big)  
  \to 0,
  \eeqnn  
 and hence $	\underline{n}\big(\overline{\epsilon}> \beta t,\zeta\leq\delta t \big) = o\big(\underline{n}\big(  \overline{\epsilon}> \beta t\big) \big)$.
 This induces that 
 \beqlb\label{eqn.423}
 \underline{n}\big(\overline{\epsilon}>\beta t\big) \sim  \underline{n}\big(\overline{\epsilon}>\beta t, \zeta >\delta t\big),
 \eeqlb 
  and then by \eqref{eqn.5021},
 \beqnn
  \limsup_{t\to\infty} \frac{\underline{n}\big(\overline{\epsilon}>\beta t\big) }{\underline{n}\big(  \zeta >  t \big) }
  \leq   \limsup_{t\to\infty} \frac{\underline{n}\big(\overline{\epsilon}>\beta t, \zeta >\delta t\big)}{\underline{n}\big(\zeta > t\big)} 
  \leq  \lim_{t\to\infty} \frac{\underline{n}\big(\zeta >\delta t\big)}{\underline{n}\big(\zeta > t\big)} = \delta^{-\theta} \to 1,
 \eeqnn
 as $\delta \to 1-$. This along with the lower bound \eqref{eqn.422} induces that $\underline{n}\big(\overline{\epsilon}>\beta t\big) \sim \underline{n}\big(\zeta >t\big)$ as $t\to\infty$. 
 
 We now prove $\underline{n}\big(\overline\epsilon >\beta t \,\big|\, \zeta >t \big)
 \sim
 \underline{n}\big( \zeta >t \,\big|\, \overline\epsilon >\beta t \big) \to 1$. Our preceding result tells that it holds if 
 \beqnn
  \underline{n}\big(\overline\epsilon >\beta t, \zeta \leq t \big) = o\big(\underline{n}(\overline\epsilon >\beta t)\big)
  \quad \mbox{and}\quad 
  \underline{n}\big(\overline\epsilon \leq \beta t, \zeta >t \big) = o\big(\underline{n}( \zeta >t)\big).
 \eeqnn 
 We just need to the first one and the second one follows immediately. 
 For any $\delta\in(0,1)$, we have 
 \beqnn
  \underline{n}\big(\overline\epsilon >\beta t, \zeta \leq t \big)= \underline{n}\big(\overline\epsilon >\beta t  \big) -\underline{n}\big(\overline\epsilon >\beta t, \zeta > \delta t \big)+ \underline{n}\big(\overline\epsilon >\beta t,\delta t< \zeta \leq t \big) .
 \eeqnn
 For the second term on the right-hand side, by \eqref{eqn.423} we have 
 \beqnn
  \limsup_{t\to\infty}\frac{\underline{n}\big(\overline\epsilon >\beta t, \zeta > \delta t \big)}{\underline{n}( \overline\epsilon >\beta t)} =1 . 
  \eeqnn
 Moreover, by \eqref{eqn.5021} we also have as $t\to\infty$,
 \beqnn
  \frac{ \underline{n}\big(\overline\epsilon >\beta t,\delta t< \zeta \leq t \big)}{\underline{n}(\overline\epsilon >\beta t)}
 \leq  
  \frac{ \underline{n}\big( \delta t< \zeta \leq t \big)}{\underline{n}(\overline\epsilon >\beta t)} 
  =  \frac{ \underline{n}\big(  \zeta >\delta  t \big)}{\underline{n}(\overline\epsilon >\beta t)} -  \frac{ \underline{n}\big(\zeta>t  \big)}{\underline{n}(\overline\epsilon >\beta t)}
  \to  \delta^{-\theta}-1,
 \eeqnn
 which goes to $0$ as $\delta \to 1-$. 
 Putting these estimates together, we see that $\underline{n}\big(\overline\epsilon >\beta t, \zeta \leq t \big) = o\big(\underline{n}(\overline\epsilon >\beta t)\big)$ and then this proof ends. 
 \qed

 \begin{lemma}\label{Lemma.405}
 Under $\underline{n}(d\epsilon)$,  the events $\zeta>J^{\beta t}$ and  $\zeta>t$ are asymptotically equivalent as $t \to \infty$, i.e., 
 \beqnn
 	\underline{n} \big(\zeta>J^{\beta t} \mid \zeta>t\big) \sim \underline{n} \big(\zeta>t \mid \zeta>J^{\beta t}\big)\to  1.
 \eeqnn
 \end{lemma}
 \proof Let  $A\Delta B$ be the symmetric difference of the two sets $A$ and $B$. 
 It suffices to prove that
 \beqlb\label{eqn.5022}
 \underline{n} \big(\{\zeta>J^{\beta t}) \Delta\{\zeta>t\}\big)=o(\underline{n} (\zeta>t)),
 \eeqlb
  as $t\to\infty$. 
 From the properties of symmetric difference, we have
 \beqlb\label{eqn4.21}
 	\underline{n} \big(\{\zeta>J^{\beta t}) \Delta\{\zeta>t\}\big)=\underline{n} \big(\zeta>t\big)-\underline{n} \big(\zeta>J^{\beta t}\big)+2 \cdot \underline{n} \big(J^{\beta t}<\zeta \leq t \big).
 \eeqlb
 Firstly, since $\{\zeta>J^{\beta t} \} = \{  N^{\beta t}_\zeta \geq 1 \}$, by \eqref{eqn.502} and Theorem~\ref{MainThm04} we have 
 \beqnn
 \underline{n} \big(\zeta>J^{\beta t}\big) =\underline{n} \big(N^{\beta t}_\zeta \geq 1\big) 
 \leq \underline{n} \big(N^{\beta t}_\zeta \big)  = \underline{n} \big( \zeta  \big)  \cdot \overline{\nu}(\beta t) 
 \quad\mbox{and then}\quad 
 \limsup_{t\to \infty} \frac{\underline{n} (\zeta>J^{\beta t})}{\underline{n} (\zeta>t)}\leq 1.
 \eeqnn
 On the other hand,   it follows from  \eqref{eqn.421} that  for any $T>0$
 \beqnn
 	\liminf_{t\to\infty}	\frac{\underline{n} (\zeta>J^{\beta t})}{\underline{n} (\zeta>t)}
 	\geq \lim_{z\to \beta+}\liminf_{t\to\infty} \frac{\underline{n} (J^{z t}\leq \zeta \wedge T)}{\underline{n} (\zeta>t)}\geq \frac{ \underline{n}(\zeta\wedge T)}{\underline{n}(\zeta)},
 \eeqnn
 which goes to 1 as $T \to \infty$. 
Combining these two estimates together,  we obtain that as $t\to\infty$,
 \beqlb\label{eqn4.22}
 	\underline{n} \big(\zeta>J^{\beta t}\big)\sim \underline{n} (\zeta>t).
 \eeqlb
 For the last term in \eqref{eqn4.21}, we see that 
 \beqnn
 \underline{n} \big(J^{\beta t}<\zeta \leq t\big)=\underline{n} \big(  \zeta>J^{\beta t}, J^{\beta t} \leq t\big)-\underline{n}\big(J^{\beta t} \leq t<\zeta\big).
 \eeqnn
For any $T>0$ and any large $t$, by using  Lemma~\ref{Lemma.401} we get 
 \beqnn
 \liminf _{t \to \infty} \frac{\underline{n} (J^{\beta t} \leq t<\zeta)}{\underline{n} (\zeta>t)} \geq
  \liminf _{t \to \infty} \frac{\underline{n} (J^{\beta t} \leq T,\zeta>t)}{\underline{n} (\zeta>t)} 
=  \lim _{t \to \infty} \underline{n} (J^{\beta t} \leq T \mid \zeta>t)=\frac{\underline{n} (\zeta \wedge T)}{\underline{n} (\zeta)},
 \eeqnn
 which goes to 1 as $T \to \infty$. 
 Moreover, by \eqref{eqn4.22}, 
 \beqnn
 \limsup _{t \to \infty} \frac{\underline{n} (\zeta>J^{\beta t}, J^{\beta t} \leq t)}{\underline{n} (\zeta>t)}\leq \lim _{t \to \infty} \frac{\underline{n} (\zeta>J^{\beta t})}{\underline{n} (\zeta>t)}=1.
 \eeqnn
 Putting these two estimates together, we conclude that $\underline{n} (J^{\beta t}<\zeta \leq t)=o(\underline{n} (\zeta >t))$ as $t\to\infty$.
 Finally,  taking this and \eqref{eqn4.22} back into the right-hand sides of \eqref{eqn4.21} immediately yields \eqref{eqn.5022}.
 \qed
 \begin{corollary}\label{Lemma.406}
 For any $ T \geq 0$ and $z \geq \beta$, we have as $t \to \infty$,
 \beqnn
  \underline{n} \big(\Delta \epsilon_{J^{\beta t}}>zt, J^{\beta t} \leq T \mid \zeta>t\big)\to (z / \beta)^{-\theta} \cdot \frac{	\underline{n} (\zeta\wedge T)}{	\underline{n} (\zeta)} .
 \eeqnn
 \end{corollary}
 \proof By Lemma \ref{Lemma.405}, we have that as $t\to\infty$,
 \beqnn
 \underline{n} \big(\Delta \epsilon_{J^{\beta t}}>zt, J^{\beta t} \leq T, \zeta>t\big) 
 \sim \underline{n} \big(\Delta \epsilon_{J^{\beta t}}>zt, J^{\beta t} \leq T, \zeta>J^{\beta t}\big) , 
 \eeqnn
 which can be bounded by $\underline{n}\big(J^{z t}\leq T\wedge \zeta \big) \leq  \underline{n}\big(  N^{z t}_{T\wedge \zeta} \big)$. 
 Then by \eqref{eqn.502} and Theorem~\ref{MainThm04},
 \beqnn
 \limsup_{t\to \infty}\underline{n} \big(\Delta \epsilon_{J^{\beta t}}>zt, J^{\beta t} \leq T \mid \zeta>t\big)
 \ar\leq\ar \limsup_{n\to\infty} \frac{\underline{n}\big(  N^{z t}_{T\wedge \zeta} \big)}{\underline{n} \big(\zeta >t\big)} = (z / \beta)^{-\theta } \cdot \frac{	\underline{n} (\zeta\wedge T)}{	\underline{n} (\zeta)} .
 \eeqnn 
 On the other hand, it is easy to identify that 
 $\{ \Delta \epsilon_{J^{\beta t}}>zt, J^{\beta t} \leq T, \zeta>J^{\beta t} \} \supset  \{N^{zt}_{T \wedge \zeta}\geq 1\}$ and hence 
  \beqnn
 \underline{n}\big(\Delta \epsilon_{\mathcal{J} \beta t}>z t, \mathcal{J}^{\beta t} \leq T, \zeta>\mathcal{J}^{\beta t}\big) \geq 
 \underline{n}\big(  N^{z t}_{T\wedge \zeta} \big)-\underline{n}\big(  N^{\beta  t}_{T\wedge \zeta}; \mathcal{N}_{u \wedge \zeta}^{\beta t}\geq 2\big).
 \eeqnn
 By \eqref{eqn.502} and \eqref{eqn.42}, we have $\underline{n}\big( N^{\beta  t}_{T\wedge \zeta} ;N_{T\wedge\zeta}^{\beta t}\geq 2 \big)=o\big( \underline{n}\big(  N^{z t}_{T\wedge \zeta} \big)\big)$ and then 
 \beqnn
  \liminf_{t\to \infty}\underline{n} \big(\Delta \epsilon_{J^{\beta t}}>zt, J^{\beta t} \leq T \mid \zeta>t\big)
 \ar\geq\ar \liminf_{n\to\infty} \frac{\underline{n}\big(  N^{z t}_{T\wedge \zeta} \big)}{\underline{n} \big(\zeta >t\big)} = (z / \beta)^{-\theta } \cdot \frac{	\underline{n} (\zeta\wedge T)}{	\underline{n} (\zeta)} . 
 \eeqnn
 The desired limit holds.
 \qed

 \textit{\textbf{Proof of Theorem~\ref{MainThm07}(1).}}
 For convention, we define $\xi :=\{ \mathcal{P}  \cdot \mathbf{1}_{\{s\geq \mathcal{T} \}},  s\geq 0 \}$. 
 Note that 
 \beqnn
 \frac{\epsilon_{s}}{t}= \Big(\frac{\epsilon_{s}}{t} - \frac{\Delta\epsilon_{J^{\beta t}}}{t} \cdot\mathbf{1}_{\{s\geq J^{\beta t}\}} \Big)+\frac{ \Delta\epsilon_{J^{\beta t}}}{t} \cdot\mathbf{1}_{\{s\geq J^{\beta t}\}}  ,\quad t>0 , \, s \geq 0. 
 \eeqnn 
  By Theorem~16.4 in \cite[p.170]{Billingsley1999}, it suffices to prove that for any $T>0$ and $\varepsilon>0$,   
  \beqlb\label{eqn.424}
 		\lim_{t\to\infty}\underline{n}\bigg( \sup_{s\in[0,T]}\Big|\frac{\epsilon_{s}}{t}- \frac{\Delta\epsilon_{J^{\beta t}}}{t} \cdot\mathbf{1}_{\{s\geq J^{\beta t}\}}\Big|\geq\varepsilon \ \Big|\  \zeta>t \bigg)=0;
 \eeqlb
 and for any non-negative bounded, uniformly continuous non-negative functional $F$ on  $D([0,T];\mathbb{R})$, 
 \beqlb\label{eqn.425}
 		\lim_{t\to \infty}\underline{n} \bigg(F\Big( \frac{\Delta\epsilon_{J^{\beta t}}}{t} \cdot\mathbf{1}_{\{s\geq J^{\beta t}\}},\,  s\in[0,T] \Big) \,\Big| \, \zeta >t\bigg)=\mathbf{E}\big[F\big(\xi_s, s \in[0,T] \big)\big].
 \eeqlb
 
 {\bf Step 1.} We first prove the limit \eqref{eqn.424}. 
 Decompose the excursion $\epsilon$ at the first big jump $J^{\beta t}$, i.e.,
 \beqlb\label{eqn.427}
 	\epsilon_{s}=\epsilon_{(s\wedge J^{\beta t})-}+\Delta\epsilon_{J^{\beta t}}\cdot \mathbf{1}_{\{s\geq J^{\beta t}\}}+\epsilon_{s\vee J^{\beta t}}-\epsilon_{J^{\beta t}},\quad s\geq 0.
 \eeqlb
 Plugging this into the left-hand side of \eqref{eqn.424} and then using  Lemma~\ref{Lemma.405}, we see that \eqref{eqn.424} holds if 
 \beqlb \label{eqn.426}
  \lim_{t\to\infty}	\underline{n}\bigg( \sup_{s\in[0,T]}\big|\epsilon_{s\vee J^{\beta t}}-\epsilon_{J^{\beta t}}\big|\geq\varepsilon t \,\Big|\,  \zeta>J^{\beta t} \bigg) 
 +  \lim_{t\to\infty}		\underline{n}\bigg( \sup_{s\in[0,T]}\epsilon_{(s\wedge J^{\beta t})-}\geq\varepsilon t \, \Big|\, \zeta>J^{\beta t} \bigg) = 0.
 \eeqlb
For the first limit, using the identity \eqref{eqn.n-Markov} to the third equality, 
 \beqnn
 \underline{n}\bigg(\sup_{s\in[0,T]}\big|\epsilon_{s\vee J^{\beta t}}-\epsilon_{J^{\beta t}}\big|\geq\varepsilon t , \zeta>J^{\beta t} \bigg)
 \ar=\ar \int_0^\infty   \underline{n}\bigg(\sup_{s\in[0,T]}\big|\epsilon_{s\vee r}-\epsilon_r\big|\geq\varepsilon t , \zeta>r, J^{\beta t}\in dr  \bigg)\cr
 \ar=\ar \int_0^T  \underline{n}\bigg(\sup_{s\in[r,T]}\big|\epsilon_{s}-\epsilon_r\big|\geq\varepsilon t , \zeta>r, J^{\beta t}\in dr \bigg)\cr
 \ar=\ar \int_0^T  \mathbf{P}\Big(\sup_{s\leq T-r}\big|X_{s}\big|\geq\varepsilon t   \Big)   \, \underline{n}\big(\zeta>r, J^{\beta t}\in dr  \big)\cr
 \ar\leq\ar     \mathbf{P}\Big(\sup_{s\leq T}\big|X_{s}\big|\geq\varepsilon t   \Big) \cdot \underline{n}\big(\zeta> J^{\beta t} \big),
 \eeqnn
 which follows that as $t\to\infty$,  
 \beqnn
  \underline{n}\bigg(\sup_{s\in[0,T]}\big|\epsilon_{s\vee J^{\beta t}}-\epsilon_{J^{\beta t}}\big|\geq\varepsilon t \,\Big|\,  \zeta>J^{\beta t} \bigg) \leq   \mathbf{P}\Big(\sup_{s\leq T}\big|X_{s}\big|\geq\varepsilon t   \Big) \to 0.
 \eeqnn
 For the second limit,   the Poissonian structure of jumps (see \eqref{eqn.Compensation}) gives that
 \beqnn
 	\underline{n}\bigg( \sup_{s\in[0,T]}\epsilon_{(s\wedge J^{\beta t})-}\geq\varepsilon t,  \zeta>J^{\beta t} \bigg) 
 	\ar\leq\ar \underline{n}\Big( \sum_{ 0\leq s< \zeta}\mathbf{1}_{ \{\Delta\epsilon_{s}>\beta t, \overline{\epsilon}_{s-}>\varepsilon t\}} \Big) 
 	= \underline{n}\Big(\int_{0}^{\zeta}\int_{\beta t}^{\infty}  \mathbf{1}_{\{\overline{\epsilon}_{s}>\varepsilon t\}} \cdot \nu(dx)\, ds \Big) , 
 \eeqnn
 which can be bounded by $\overline{\nu}(\beta t)\cdot \underline{n}\big(\zeta\cdot \mathbf{1}_{\{\overline{\epsilon}>\varepsilon t\}} \big) $. 
 Recall that $\underline{n}(\zeta>J^{\beta t}) \sim \underline{n}(\zeta>t) \sim \bar{\nu}(\beta t)/\beta$; see Lemma~\ref{Lemma.405} and Theorem~\ref{MainThm04}. 
 This along with $\underline{n}(\zeta)<\infty$ induces that as $t\to\infty$,
 \beqnn
 \underline{n}\bigg( \sup_{s\in[0,T]}\epsilon_{(s\wedge J^{\beta t})-}\geq\varepsilon t \, \Big|\, \zeta>J^{\beta t} \bigg)
 \leq \frac{\overline{\nu}(\beta t)}{ \underline{n}(\zeta>J^{\beta t})} \cdot \underline{n}\big(\zeta\cdot \mathbf{1}_{\{\overline{\epsilon}>\varepsilon t\}} \big) 
 \leq C\cdot \underline{n}\big(\zeta \cdot \mathbf{1}_{\{\overline{\epsilon}>\varepsilon t\}} \big) \to0.
 \eeqnn
 Here we have proved the limit \eqref{eqn.424}.
 
  {\bf Step 2.} To prove the limit  \eqref{eqn.425},  we define a one-step process 
  \beqnn
  M^{(t)}_s:= t^{-1}\cdot \Delta\epsilon_{J^{\beta t}}\cdot\mathbf{1}_{\{s\geq J^{\beta t}\}},\quad s\geq0.
  \eeqnn
 Obviously, it holds if and only if $M^{(t)}\to \xi $ weakly in $D([0,T];\mathbb{R})$. 
 Firstly, the finite-dimensional convergence follows directly from  Corollary~\ref{Lemma.406}. It remains to identify the tightness of the sequence $\{ M^{(t)} \}_{t>0}$, which is obvious by using Theorem~13.5 in \cite[p.142]{Billingsley1999} as well as the two facts that 
 \beqnn
 \big|M^{(t)}_{s_2}-M^{(t)}_{s_1}\big| \wedge 
 \big|M^{(t)}_{s_3}-M^{(t)}_{s_2}\big| \overset{\rm a.s.}= 0
 \eeqnn
 for all $0\leq s_1<s_2<s_3\leq T$ and
 $\xi^{\alpha,\beta}$ is stochastically continuous. 
 \qed

  \textit{\textbf{Proof of Theorem~\ref{MainThm07}(2).}} 
  Recall the process $\mathcal{C}^{X,t}$ defined before Lemma~\ref{Lemma.302} for $t>0$. We consider the rescaled process $\xi^{(t)}$ with 
  \beqnn
  \xi^{(t)}_s= \frac{\mathcal{C}^{X,t}_{ts}}{t},\quad s\geq 0,
  \eeqnn
  and also make the convention that $\xi^*:=\{ (\mathcal{P} -\beta s)\vee 0, s> 0 \}$.
  Notice that $0$ is regular for $(-\infty,0)$, for $\xi^*$. 
  By Lemma~\ref{Lemma.302} and \ref{Lemma.304}, it suffices to prove that $\xi^{(t)} \overset{\rm d}\to \xi^*$ in $D\big((0,\infty);\mathbb{R}\big)$ as $t\to\infty$. 
 Since 
 \beqnn
 \{ \xi^{(t)}_{s+1}- \xi^{(t)}_{1} :s\geq 0 \}\overset{\rm d}=\{X_{ts}/t:s\geq 0\}
 \eeqnn 
 and is independent of $\{ \xi^{(t)}_s:s\in[0,1] \}$, we just need to prove that as $t\to\infty$,
 \beqnn
 \xi^{(t)}  \overset{\rm d}\to \xi^*,
 \eeqnn
 in $D((0,1];\mathbb{R})$ and 
 \beqnn
 \{X_{ts}/t:s\geq 0\}\overset{\rm p}\to \{ -\beta s:s \geq 0 \},
 \eeqnn  
 uniformly on compacts. 
 The second limit follows directly from Proposition~\ref{Prop.402}. 
 It remains to prove that  $\xi^{(t)}  \overset{\rm d}\to \xi^*$ in $D([\delta,1];\mathbb{R})$ for any $\delta \in(0,1)$, which follows immediately if  for any $\varepsilon>0$, 
 \beqlb\label{eqn.429}
 \lim_{t\to\infty}\underline{n}\bigg( \sup_{s\in[\delta,1]}\Big|\frac{\epsilon_{ts}}{t}- \frac{\Delta\epsilon_{J^{\beta t}}-\beta ts}{t} \Big|\geq\varepsilon \ \Big|\  \zeta>t \bigg)  =0,
\eeqlb
and for any non-negative, bounded, uniformly continuous function $F$ on $D([\delta,1];\mathbb{R})$,  
\beqnn
 \lim_{t\to\infty}\underline{n} \bigg(F\Big( \frac{\Delta\epsilon_{J^{\beta t}}}{t} -\beta s ,\,  s\in[\delta ,1] \Big) \,\Big| \, \zeta >t\bigg)
 = \mathbf{E}\big[ F\big(  \xi^*_s, s\in[\delta,1] \big) \big].
 \eeqnn
 The second limit follows directly from  Corollary~\ref{Lemma.406}.
 We now prove the limit \eqref{eqn.429}. 
 Since $\Delta\epsilon_{J^{\beta t}}=  \epsilon_{J^{\beta t}}- \epsilon_{J^{\beta t}-}$, the left-hand side of \eqref{eqn.429} can be bounded by 
 \beqlb\label{eqn.432}
    \lim_{t\to\infty}\underline{n}\bigg( \sup_{s\in[\delta,1]}\big|\epsilon_{ts}- \epsilon_{J^{\beta t}} +\beta ts \big|\geq \frac{\varepsilon t}{2}\,\Big|\,  \zeta>t \bigg) 
    +   \lim_{t\to\infty} \underline{n}\Big(    \frac{\epsilon_{J^{\beta t}-} }{t} \geq \frac{\varepsilon }{2} \,\Big|\,  \zeta>t \bigg).
 \eeqlb
 Theorem~\ref{MainThm07}(1) tells that the second term equals to $0$. 
 For the first term, by Lemma~\ref{Lemma.401} we have $\underline{n}(J^{\beta t}> \delta t\,|\,\zeta >t)\to 1$ and hence it remains to prove that
 \beqlb\label{eqn.441}
   \lim_{t\to\infty} \underline{n}\bigg( \sup_{s\in[\delta,1]}\big|\epsilon_{ts}- \epsilon_{J^{\beta t}} +\beta ts \big|\geq \frac{\varepsilon t}{2} ,J^{\beta t}\leq \delta t \,\Big|\,  \zeta>t \bigg)=0.
 \eeqlb
 By  the  Markov property and the independent increments of $\epsilon$ under $\underline{n}(\cdot,\zeta>t)$; see \eqref{eqn.n-Markov},  
 \beqnn
  \lefteqn{ \underline{n}\bigg( \sup_{s\in[\delta,1]}\big|\epsilon_{ts}- \epsilon_{J^{\beta t}} +\beta ts \big|\geq \frac{\varepsilon t}{2} ,J^{\beta t}\leq \delta t ,  \zeta>t \bigg)}\ar\ar\cr
  \ar=\ar  \int_0^{\delta t} \int_{\beta t}^\infty  \underline{n}\bigg( \sup_{s\in[\delta,1]}\big|\epsilon_{ts}- \epsilon_{J^{\beta t}} +\beta ts \big|\geq \frac{\varepsilon t}{2} ,J^{\beta t}\leq \delta t ,  \zeta>t ,J^{\beta t} \in dr, \epsilon_{r}\in dx \bigg) \cr
   \ar=\ar \int_0^{\delta t} \int_{\beta t}^\infty  \mathbf{P}_x\bigg(  \sup_{s\in[\delta-r/t,1-r/t]}\big|X_{ts}-x+\beta (ts+r)\big| \geq \frac{\varepsilon t}{2} , \tau_0^->t-r \bigg) \,\underline{n}\big( J^{\beta t} \in dr, \epsilon_{r}\in dx ,  \zeta>r\big)\cr
   \ar\leq \ar \int_0^{\delta t} \int_{\beta t}^\infty  \mathbf{P}\bigg(  \sup_{s\in[0,1]}\big|X_{ts}+\beta ts\big| + \beta r \geq \frac{\varepsilon t}{2} \bigg) \,\underline{n}\big( J^{\beta t} \in  dr, \epsilon_{r}\in dx,  \zeta>r\big) \cr
   \ar\leq\ar  \int_0^{\frac{\varepsilon t}{4\beta}} \int_{\beta t}^\infty  \mathbf{P}\bigg(  \sup_{s\in[0,1]}\big|X_{ts}+\beta ts\big|   \geq \frac{\varepsilon t}{4} \bigg) \,\underline{n}\big( J^{\beta t} \in  dr, \epsilon_{r}\in dx,  \zeta>r\big)   +\underline{n}\Big( \frac{\varepsilon t}{4\beta}< J^{\beta t} \leq \delta t,  \zeta> J^{\beta t}\Big)\cr
   \ar\leq \ar \mathbf{P}\bigg(  \sup_{s\in[0,1]}\big|X_{ts}+\beta ts\big|  \geq  \frac{\varepsilon t}{4} \bigg) \cdot  \underline{n}\big( J^{\beta t} \leq \delta t,  \zeta> J^{\beta t}\big) + 
  \underline{n}\Big(  J^{\beta t} > \frac{\varepsilon t}{4\beta},  \zeta> J^{\beta t}\Big).
 \eeqnn 
 By Proposition~\ref{Prop.402}, the first probability on the right-hand side of the last inequality vanishes as $t\to\infty$. 
 Moreover, by using Lemma~\ref{Lemma.405} and Lemma~\ref{Lemma.401}, we have $ \underline{n}\big( J^{\beta t} \leq \delta t,  \zeta> J^{\beta t}\big) \sim \underline{n} ( \zeta> t)$ and 
 \beqnn
 \frac{\underline{n}\big(  J^{\beta t} > \frac{\varepsilon t}{4\beta},  \zeta> J^{\beta t}\big)}{\underline{n} ( \zeta> t)} 
 \sim  \frac{\underline{n}\big(  J^{\beta t}>  \frac{\varepsilon t}{4\beta},  \zeta> t\big)}{\underline{n} ( \zeta> t)}
 = \underline{n}\Big(  J^{\beta t}> \frac{\varepsilon t}{4\beta}\,\Big|\, \zeta> t\Big)  ,
 \eeqnn
  which vanishes as $t\to\infty$. 
  Consequently, the limit \eqref{eqn.441} holds and the proof ends.
 \qed

 \begin{lemma}\label{Lemma.403}
  For any $u>0$, we have  $\underline{n} (\epsilon_u>t, \zeta >u )
  \sim \underline{n}(u\wedge\zeta) \cdot\bar{\nu}(t)$ as $t\to\infty$. 
 \end{lemma}	
 \proof Equivalently, it suffices to prove that 
 $	\underline{n}\big(\epsilon_u>2\beta  t, \zeta >u	\big)
 \sim \underline{\tt n}(u\wedge\zeta) \cdot\bar{\nu}\big(2\beta  t\big)$. 
 For large $t$, the left-hand side can be decomposed into the next two terms
 \beqlb\label{eqn.438}
 \underline{n}\big(\epsilon_u>2\beta  t, u<\zeta \leq t	\big)
  \quad\mbox{and}\quad
 \underline{n}\big(\epsilon_u>2\beta  t, \zeta >t	\big) . 
 \eeqlb
 Applying the Markov property to the first term, we have 
 \beqnn
  \underline{n}\big(\epsilon_u>2\beta  t, u<\zeta \leq t	\big)
  \ar=\ar \int_{2\beta  t}^{\infty}\mathbf{P}_y\big(\tau_{0}^{-}\leq t-u\big)\, \underline{  n}\big(\epsilon_u\in dy ,\zeta >u\big),
 \eeqnn
 which, obviously, can be bounded by $ \mathbf{P}_{2\beta  t}\big(\tau_{0}^{-}\leq  t \big) \cdot \underline{n}\big( \epsilon_u>2\beta t, \zeta >u	\big)$. 
 Note that $\frac{3\beta}{2}  t +X_t \to \infty$ as $t\to\infty$, we see that
 \beqnn
 \mathbf{P}_{2\beta  t}\big(\tau_{0}^{-}\leq  t \big)
 \ar=\ar \mathbf{P}\big(2\beta  t + \underline{X}_t \leq 0\big)
 \leq \mathbf{P}\Big(\inf_{s\in[0,t]} \Big(\frac{3\beta}{2}  s + X_s\Big)  \leq -\frac{\beta t}{2}\Big) \to0,
 \eeqnn
 and then $ \underline{n}\big(\epsilon_u>2\beta t, u<\zeta \leq t	\big) = o\big( \underline{n}  (\epsilon_u>2\beta t,\zeta >u	 ) \big)$. 
 For the second term in \eqref{eqn.438}, Theorem~\ref{MainThm07}(1) and Theorem~\ref{MainThm04} ensures that as $t\to\infty$, 
 \beqnn
 \underline{n}\big(\epsilon_u>2\beta  t, \zeta >t	\big) 
 \ar=\ar \underline{n}\big(\epsilon_u>2\beta  t\,\big|\, \zeta >t	\big) \cdot  \underline{n}\big(  \zeta >t	\big) 
 \sim \frac{\underline{ n}(u\wedge\zeta)}{\underline{ n}(\zeta)}\cdot \frac{  \underline{n}\big(  \zeta >t	\big)}{2^{\alpha}} 
 \sim \underline{ n}(u\wedge\zeta) \cdot \frac{\bar{\nu}(\beta t)}{2^{\theta}}    ,
 \eeqnn
 which is asymptotically equivalent to $  \underline{ n}(u\wedge\zeta) \cdot\bar{\nu}(2\beta t)$ since $\bar{\nu} \in \mathrm{RV}^\infty_{-\theta}$.
 The proof ends.  
 \qed
 
 \begin{lemma}\label{Lemma.404}
 	For any $x> 0$, we have $\mathbf{P}_x(X_u>t,\,\tau_0^{-}>u)\sim \mathbf{E}_x[\tau_0^{-}\wedge u]\cdot\bar{\nu}(t).$ as $t\to\infty$. 
 	\end{lemma}  
 \proof  We first decompose the left-hand side at the stopping time $J^t$ into  the next two parts:
 \begin{align}\label{8.1.1}
 \mathbf{P}_x\big(X_u>t, J^{ t}\leq u< \tau_0^- \big)
 \quad \mbox{and}\quad 
 \mathbf{P}_x\big(X_u>t, J^{ t}> u, \tau_0^- >u\big).
 \end{align}
 To obtain the desired limit, we prove in the next two steps that after dividing by $\bar{\nu}(t)$, these two probabilities converge to $\mathbf{E}_x[\tau_0^{-}\wedge u]$ and $0$ respectively. 
 
 {\bf Step 1.} For the first probability, we see that $X_u=X_{J^{ t}-}+\Delta  X_{J^{ t}}+X_u-X_{J^{ t}}$  and then for any $\varepsilon>0$,  
 \beqlb\label{eqn.812}
 \lefteqn{\mathbf{P}_x\bigg( \Big|\frac{X_u}{t}-\frac{\Delta X_{J^{ t}}}{t}\Big|>\varepsilon,J^{ t}\leq u< \tau_0^-  \bigg)  }\ar\ar\cr
 \ar\ar\cr
  \ar\leq\ar \mathbf{P}_x\Big( X_{J^{ t}-}> \frac{\varepsilon t}{2},J^{ t}\leq u< \tau_0^- \Big) 
  + \mathbf{P}_x\Big( \big|X_u-X_{J^{ t}}\big|>\frac{\varepsilon t}{2},J^{ t}\leq u< \tau_0^-   \Big). 
 \eeqlb
 Note that the event $\{X_{J^{ t}-}>\varepsilon t/2,J^{ t}\leq u< \tau_0^-\}$ occurs if and only if 
 \beqnn
  \#\big\{s\in[0,u]: \Delta X_{s}> t,X_{s-}>\varepsilon t/2, \underline{X}_{s}\geq 0 \big\} \geq 1.
 \eeqnn
 By using the compensation formula \eqref{eqn.CompenstationX} with $\mathbf{F}_s(\Delta X_{s}):= \mathbf{1}_{\{\Delta X_{s}> t,X_{s-}>\varepsilon t/2, \underline{X}_{s}\geq 0\}}$, we have  
 \beqnn
  \mathbf{P}_x\big( X_{J^{ t}-}>\varepsilon t/2,J^{ t}\leq u< \tau_0^-\big) 
  \ar\leq\ar \mathbf{E}_x\bigg[\int_{0}^{u}\int_{ t}^{\infty}  \mathbf{1}_{ \{ X_{s}>\varepsilon t/2, \underline{X}_{s}\geq 0\}} \nu(dy)ds \bigg]\cr
  \ar=\ar \bar{\nu}(t)\cdot \int_{0}^{u} \mathbf{P}_x\big(  X_{s}>\varepsilon t/2, \underline{X}_{s}\geq 0 \big)\,ds\cr
  \ar\leq\ar u\cdot \bar{\nu}(t)\cdot \mathbf{P}_x\bigg( \sup_{0\leq s\leq \tau_0^-} X_{s}>\frac{\varepsilon t}{2} \bigg),
 \eeqnn
 which equals to $o\big(\bar{\nu}(t)\big)$ as $t\to\infty$,  since  the last probability vanishes. 
 Moreover, by the Markov property, the second probability on the right-hand side of \eqref{eqn.812} equals to
 \beqnn
 \lefteqn{ \int_0^u\int_t^\infty  \mathbf{P}_x\big( \big|X_u-X_{J^{ t}}\big|>\varepsilon t/2, \tau_0^- >u, J^{ t}\in ds,  X_s\in dy \big)}\ar\ar\cr
 \ar=\ar \int_{0}^{u} \int_{ t}^{\infty} \mathbf{P}_y\big( |X_{u-s}-y|>\varepsilon t/2,\tau_0^- >u-s \big)\cdot \mathbf{P}_x\big(J^{ t}\in ds, X_s \in dy, \tau_0^- >s\big)\cr
 \ar\leq\ar \int_{0}^{u} \int_{ t}^{\infty} \mathbf{P}\big( |X_{u-s}|>\varepsilon t/2 \big)\cdot \mathbf{P}_x\big(J^{ t}\in ds, X_s \in dy, \tau_0^- >s\big) , 
 \eeqnn
 which, along with the two facts that $\sup_{s\in[0,u]} \mathbf{P}\big( |X_{u-s}|>\varepsilon t/2 \big)\leq  \mathbf{P}\big( \sup_{s\in [0,u]}|X_{s}|>\varepsilon t/2 \big) \to 0$ as $t\to\infty$ and $ \mathbf{P}_x\big(J^{ t}\leq u \big) = \mathbf{P}\big(J^{ t}\leq u \big)= u\cdot \bar{\nu}( t)$, induces that 
 \beqnn
 \mathbf{P}_x\Big( \big|X_u-X_{J^{ t}}\big|>\varepsilon t/2,J^{ t}\leq u< \tau_0^-  \Big) 
 = o\Big( \mathbf{P}_x\big(J^{ t}\leq u, \tau_0^- > J^{ t}\big) \Big) 
 \leq o\Big( \mathbf{P}_x\big(J^{ t}\leq u \big) \Big) 
 =  o\big(u\cdot \bar{\nu}( t)\big).
 \eeqnn
 Combining the preceding two estimates together, we obtain that  as $t\to\infty$,
 \begin{align}\label{8.1.5}
 	\mathbf{P}_x\big( \big|X_u/t-\Delta X_{J^{ t}}/t\big|>\varepsilon,J^{ t}\leq u< \tau_0^-   \big) = o\big(\bar{\nu}( t)\big).
 \end{align}
 and hence the first probability in \eqref{8.1.1} can be well approximated by  
 \beqnn
 \mathbf{P}_x\big(\Delta X_{J^{ t}}>t, J^{ t}\leq u<\tau_0^- \big)
 =\mathbf{P}_x\big(J^{ t}\leq u<\tau_0^- \big)
 =\mathbf{P}_x\big(J^{ t}\leq u, J^{ t} <\tau_0^- \big)-\mathbf{P}_x\big( J^{ t}<\tau_0^-\leq u\big)  .
 \eeqnn
 By Theorem 3.4 and Lemma 3.5 in \cite{Xu2021a}, we have $\mathbf{P}_x\big(J^{ t}\leq u, J^{ t} <\tau_0^- \big) \sim \mathbf{E}_x[\tau_0^{-}\wedge u]\cdot\bar{\nu}(t)$ as $t\to\infty$.
 Meanwhile, by using  the strong Markov property,
 \beqnn
 \mathbf{P}_x\big( J^{ t}<\tau_0^-\leq u\big) 
  \ar\leq\ar \mathbf{P}_x\Big(  \sup_{s\in[J^{ t},u]}\big|X_s-X_{J^{ t}}\big|> t,J^{ t}\leq u   \Big)\cr
  \ar\leq\ar \mathbf{P}_x\Big(  \sup_{s\in[J^{ t},J^{ t}+u]}\big|X_s-X_{J^{ t}}\big|> t,J^{ t}\leq u  \Big)\cr
  \ar=\ar \mathbf{P}_x\big( J^{ t}\leq u \big)  \cdot \mathbf{P}\Big( \sup_{s\in[0,u]}\big|X_s\big|> t\Big) ,
 \eeqnn
 which equals to $o\big(\bar{\nu}( t)\big)$ as $t\to\infty$ since $ \mathbf{P}_x\big(J^{ t}\leq u \big) = u\cdot \bar{\nu}( t)$.
 In conclusion, we have as $t\to\infty$, 
 \beqnn
 \mathbf{P}_x\big(X_u>t, J^{ t}\leq u< \tau_0^- \big)
 \sim \mathbf{P}_x\big(\Delta X_{J^{ t}}>t, J^{ t}\leq u<\tau_0^- \big)
 \sim \mathbf{E}_x[\tau_0^{-}\wedge u]\cdot\bar{\nu}(t).
 \eeqnn
 
 {\bf Step 2.} We now consider the second probability in  \eqref{8.1.1}, which can be bounded by
 \begin{align}\label{8.1.7}
 \mathbf{P}_x\big(X_u>t, J^{ t}> u \big)
 = \mathbf{P}_{x}\big( X_{u}>t,N^{t}_{u}=0 \big)
  =\mathbf{P}_{x}\big( X_{u}>t \big) - \mathbf{P}_{x}\big( X_{u}>t,N^{t}_{u}\geq 1 \big).
 \end{align} 
 Let $\widetilde{X}^1$ be the process defined by moving all jumps larger than $1$ away from $X$, i.e.,
 \beqnn
 \widetilde{X}_t^{1} = X_t - \sum_{s\leq t} \Delta X_s \cdot \mathbf{1}_{\{\Delta X_s>1\}},\quad t\geq 0. 
 \eeqnn
 For any $\delta >0$, conditionally on $N_{u}^{(1+\delta)t}=1  $ we have $X_{u}\geq \widetilde{X}_{u}^{1}+\Delta X_{J^{(1+\delta)t}}\geq \widetilde{X}_{u}^{1}+(1+\delta)t$ and then 
 \beqnn
  \big\{X_{u}>t, N_{u}^{t}=1\big\} \supseteq \big\{ \widetilde{X}_{u}^{1}+(1+\delta)t>t \big\}\cap \big\{ N_{u}^{t}-N_{u}^{(1+\delta)t}=0 \big\} \cap \big\{ N_{u}^{(1+\delta)t}=1 \big\}.
 \eeqnn
 It is easy to identify that the three events on the right-hand side are mutually independent and then 
 \beqnn
 	\mathbf{P}_x\big(X_{u}>t, N_{u}^{t}=1\big)\geq \mathbf{P}_x \big( \widetilde{X}_{u}^{1}>-\delta t \big)\cdot \mathbf{P}_x \big(N_{u}^{t}-N_{u}^{(1+\delta)t}=0 \big)\cdot \mathbf{P}_x \big( N_{u}^{(1+\delta)t}=1 \big) .
 \eeqnn
 Since $\mathbf{E}[\widetilde{X}_1^{1}]<0$, we have  $ \mathbf{P}_{x}\big(  \widetilde{X}_{u}^{1}>-\delta t \big) \to 1$ as $t\to\infty$.
 Moreover, since $N_{u}^{t}-N_{u}^{(1+\delta)t}$ and $N_{u}^{(1+\delta)t}$ are Poisson distributed with rate $u\cdot \big(\bar{\nu}(t)-\bar{\nu}((1+\delta)t)\big)$ and $u\cdot \bar{\nu}((1+\delta)t)$ respectively, we also have 
 \beqnn
 \mathbf{P}_{x}\big( N_{u}^{t}-N_{u}^{(1+\delta)t}=0 \big) \to 1
 \quad \mbox{and}\quad 
 \mathbf{P}_{x}\big( N_{u}^{(1+\delta)t}=1 \big) \sim  u\cdot\bar{\nu}\big((1+\delta)t\big),
 \eeqnn 
 as $t\to\infty$. 
 Combining these estimates together and then using the fact that $\bar{\nu} \in \mathrm{RV}^\infty_{-\theta}$,
 \beqnn
 \liminf_{t\to\infty} \frac{\mathbf{P}_{x}\big( X_{u}>t,N^{t}_{u}\geq 1 \big)}{\bar{\nu}(t)} \geq  u\cdot \lim_{t\to\infty} \frac{\bar{\nu}\big((1+\delta)t\big)}{\bar{\nu}(t)} = \frac{u}{ (1+\delta)^{\theta}}.
 \eeqnn 
 Additionally, repeating the proof of (3.2) in \cite{Xu2021a} gives that  
 $\mathbf{P}_{x}\big( X_{u}>t \big) \sim u\cdot \bar{\nu}(t)$ as $t\to\infty$.
 Taking these two estimates back into \eqref{8.1.7}, we have 
 \begin{align}\label{8.1.10}
 	 \limsup_{t\to\infty}	\frac{1}{\bar{\nu}(t)}\cdot \mathbf{P}_{x}\big( X_{u}>t,J^{ t}>u,\tau_0^->u \big)
 	\leq  u\cdot\big(1-(1+\delta)^{-\theta}\big),
 \end{align}
 which goes to $0$ as $\delta \to 0+$.  The whole proof ends.
 \qed
 
 \begin{lemma}\label{Lemma.408}
 For any $x>0$, we have the following identity
 \beqlb
  \mathbf{P}\big(\mathcal{T}_x \leq u\big)
  =\frac{1}{\widehat{V}(x)}\int_0^{\infty} \mathbf{P}\big(\mathcal{T}\leq u-s\big)\cdot \bar{n}(\epsilon_s\le x,\,\zeta>s)\,ds 
  + \frac{\bar{\tt d}}{\widehat{V}(x)}\,\mathbf{P}\big(\mathcal{T}\leq u\big),\quad u\geq 0.
 \eeqlb
 \end{lemma}
 	\proof  For any $x, t,u>0$,  Lemma \ref{Lemma.201} allows us to write the probability $\mathbf{P}_x(X_u>t,\,\tau_0^{-}>u)$ as
 	\beqlb\label{eqn4.30}  
    \bar{\tt d}\cdot \underline{n}\big(\epsilon_u>t-x,\,\zeta>u \big) 
    + \int_0^u \int_0^x \underline{ n}\big(\epsilon_{u-s}>t+z-x,\,\zeta>u-s\big)\cdot \overline{n}\big(\epsilon_s\in dz,\,\zeta>s\big)\,ds.  
 	\eeqlb
 From Lemma \ref{Lemma.403}, we see that  
 $\underline{n}(\epsilon_u>t-y,\,\zeta>u)\sim\underline{n}(\zeta\wedge u)\cdot\bar{\nu}(t)$
  as $t\to\infty$ uniformly in $y\in[0,x]$.
  Moreover, by Theorem~\ref{MainThm04}, 
  \beqnn
 \limsup_{t\to\infty} \sup_{s\in[0,u]}\sup_{z\in[0,x]} \frac{\underline{ n}\big(\epsilon_{u-s}>t+z-x,\,\zeta>u-s\big)}{\bar{\nu}(t)} 
 \leq   \limsup_{t\to\infty} \frac{\underline{ n}(\overline{\epsilon}>t-x)}{\bar{\nu}(t)} 
 = \underline{ n}(\zeta)<\infty. 
  \eeqnn
  Based on these two asymptotic estimates, we divide   \eqref{eqn4.30} by $\bar{\nu}(t)$ and then apply the dominated convergence theorem to the last integral to obtain that 
 	\begin{align}\label{recall} 
 	\lim_{t\to\infty}\frac{1}{\bar{\nu}(t)}\cdot \mathbf{P}_x(X_u>t,\,\tau_0^{-}>u)
 	=  \bar{\tt d}\cdot\underline{n}(\zeta\wedge u)+\int_0^{\infty}
 	\underline{n}\big(\zeta\wedge((u-s)\vee 0)\big)\cdot
 	\bar{n}(\epsilon_s\leq x,\,\zeta>s)\,ds.
 	\end{align} 
 	On the other hand, by Lemma \ref{Lemma.404} we also have the left-hand side also equals to $\mathbf{E}_x[\tau_0^{-}\wedge u]$ and then
  \beqnn
  \mathbf{E}_x[\tau_0^{-}\wedge u]=\bar{\tt d}\cdot\underline{n}(\zeta\wedge u)+\int_0^{\infty}
  \underline{n}\big(\zeta\wedge((u-s)\vee 0)\big)\cdot
  \bar{n}(\epsilon_s\leq x,\,\zeta>s)\,ds.
  \eeqnn
  Dividing both sides by $\mathbf{E}_{x}[\tau_0^-]$, we see that the desired identity follows from the distributions of $\mathcal{T}_x$ and $\mathcal{T} $ as well as the equality $ \mathbf{E}_{x}[\tau_0^-]=\widehat{V}(x)\cdot \underline{n}(\zeta)$; see \eqref{Asytau}.
 \qed

 \textit{\textbf{ Proof of Corollary~\ref{MainThm08}(1).}}
 It suffices to prove that for any $T>0$ and any bounded, uniformly continuous non-negative functional $F$ on $D([0,T];\mathbb{R})$,
 \beqnn
 \mathbf{E}_x \big[F\big( X_s/t,\,  s\in [ 0,T] \big)\,\big|\, \tau_{0}^- >t\big]\to\mathbf{E}\big[F\big(\mathcal{P}  \cdot \mathbf{1}_{\{\mathcal{T}_x \leq s\}}, s\in [ 0,T] \big)\big],
 \eeqnn
 as $t\to\infty$. 
 For any $0<K< t$, we decompose the left-hand side into the following two parts:
 \beqlb\label{eqn.442}
 	  \mathbf{E}_x\big[ F( X_{s}/t,  s\in[0,T] ), \underline{g}_t>K\,\big|\,\tau_0^->t \big]
  \quad	\mbox{and}\quad  
  \mathbf{E}_x\big[ F( X_{s}/t,  s\in[0,T] ), \underline{g}_t\leq K\,\big|\,\tau_0^->t \big] .
 \eeqlb
 Without loss of generality, we may assume that $F(\cdot)\leq 1$. 
 The first conditional expectation can be bounded by 
 \beqlb\label{eqn.4341}
 \mathbf{E}_x\big[ \underline{g}_t> K \,\big|\,  \tau_0^->t \big] 
 \ar=\ar \frac{\mathbf{P}_x\big(  \underline{g}_t> K,  \tau_0^->t \big) }{\mathbf{P}_x\big( \tau_0^->t \big) }
 = 1- \frac{\mathbf{P}(\underline{g}_t\leq K, \underline{X}_t\geq -x )}{\mathbf{P}_x(\tau_0^->t)}.
 \eeqlb
 By Lemma~\ref{Lemma.201}, the last fraction can be represented as
 \beqnn
 \frac{1}{\mathbf{P}_x(\tau_0^->t)} \cdot \Big( \bar{d}\cdot \underline{n} (\zeta >t)+
 \int_{0}^{K}  \bar{n}( \epsilon_s\leq x,  \zeta>s ) \underline{n}(\zeta>t-s)\, ds \Big) . 
 \eeqnn
 By \eqref{Asytau} and Theorem~\ref{MainThm04}, we have $\mathbf{P}_x(\tau_0^->t) \sim \widehat{V}(x) \cdot \underline{n} (\zeta >t-s)$ uniformly in $s\in[0,K]$ and hence  
 \beqlb\label{eqn.4351}
 \lim_{t\to\infty}\frac{\mathbf{P}(\underline{g}_t\leq K, \underline{X}_t\geq -x )}{\mathbf{P}_x(\tau_0^->t)} 
 = \frac{1}{\widehat{V}(x)} \Big( \bar{d} + \int_{0}^K  \overline{n}( \epsilon_s\leq x,  \zeta>s) \, ds \Big) ,
 \eeqlb 
 which goes to $1$ as $K\to\infty$ by \eqref{eqn.447}.  
 Taking this back into \eqref{eqn.4341}, we have 
 \beqnn
  \lim_{K\to\infty} \lim_{t\to\infty}\mathbf{E}_x\big[ F( X_{s}/t,  s\in[0,T] ), \underline{g}_t>K\,\big|\,\tau_0^->t \big] =0.
 \eeqnn
 We now consider the second conditional expectation in \eqref{eqn.442}, which equals to 
 \beqnn
 	\frac{1}{\mathbf{P}_x( \tau_0^->t )} \cdot \mathbf{E} \bigg[ F\Big( \frac{x}{t} +  \frac{X_{s}}{t}:s\in[0,T] \Big), \underline{g}_t\leq K, \underline{X}_t\geq -x \bigg].
 \eeqnn
 By the uniform continuity of $F$, it further equals to
 \beqlb \label{eqn.444}
	\frac{1}{\mathbf{P}_x( \tau_0^->t )}  \cdot  \Big( \mathbf{E} \big[ F\big(   X_{s}/t:s\in[0,T] \big), \underline{g}_t\leq K, \underline{X}_t\geq -x \big] + O(x/t)\cdot  \mathbf{P} \big( \underline{g}_t\leq K, \underline{X}_t\geq -x \big) \Big) . 
 \eeqlb 
 This along with \eqref{eqn.4351} induces that as $t\to\infty$,
 \beqnn
 \mathbf{E}_x\big[ F( X_{s}/t,  s\in[0,T] ), \underline{g}_t\leq K\,\big|\,\tau_0^->t \big]
 \sim \frac{1}{\mathbf{P}_x( \tau_0^->t )}  \cdot    \mathbf{E} \big[ F\big(   X_{s}/t:s\in[0,T] \big), \underline{g}_t\leq K, \underline{X}_t\geq -x \big].
 \eeqnn
 By Lemma~\ref{Lemma.201}, the last expectation can be represented as
 \beqnn
 \bar{\tt d}\cdot\underline{n}\Big(F\Big(\frac{\epsilon_s}{t},s\in[0,T]\Big),\zeta>t\Big) \ar+\ar\int_{0}^{K}ds \int_{\mathcal{E}} \overline{n}( d\epsilon,\epsilon_s\leq x,  \zeta>s )  \int_\mathcal{E}F\Big(\frac{(\epsilon,\epsilon^*)^s_r}{t}, r\in[0,T] \Big)\, \underline{n}(d\epsilon^*, \zeta>t-s). 
 \eeqnn
 From \eqref{Asytau} and Theorem~\ref{MainThm07}(1), it follows that  as $t\to\infty$, 
 \beqlb\label{eqn.4371}
  \frac{\bar{\tt d}\cdot\underline{n}\big(F\big(\epsilon_{s}/t,s\in[0,T]\big),\zeta>t\big)}{\mathbf{P}_x( \tau_0^->t )} 
  \ar=\ar  \frac{\bar{\tt d}\cdot\underline{n}\big( \zeta>t\big)}{\mathbf{P}_x( \tau_0^->t )} \cdot \underline{n}\big(F\big(\epsilon_{s}/t,s\in[0,T]\big)\,\big|\, \zeta>t\big)\cr
  \ar\sim\ar \frac{\bar{\tt d} }{\widehat{V}(x)}\cdot \mathbf{E}\big[F( \mathcal{P}\cdot\mathbf{1}_{\{\mathcal{T}\leq s \}}, s\in[0,T])\big].
  \eeqlb
 Similarly as in Step 2 of the proof of Corollary~\ref{MainCorollary01},  by the uniform continuity of $F$ we have  as $t\to\infty$,
  \beqnn 
  \lefteqn{ \frac{1}{ \mathbf{P}_x(\tau_0^->t)}\int_{0}^{K}ds \int_{\mathcal{E}}    \overline{n}( d\epsilon,\epsilon_s\leq x,  \zeta>s )  \int_\mathcal{E}F\Big( \frac{(\epsilon,\epsilon^*)^s_r}{t}, r\in[0,T] \Big)\, \underline{n}(d\epsilon^*, \zeta>t-s)}\ar\ar\cr
  \ar\sim\ar  \int_{0}^{K} \overline{n}(  \epsilon_s\leq x,  \zeta>s )   \cdot \underline{n}\Big(F\Big(\frac{ \epsilon_{r-s}}{t}\cdot \mathbf{1}_{\{r\geq s \}}, r\in[0,T] \Big)\,\Big|\, \zeta>t-s \Big)\cdot \frac{ \underline{n}(\zeta >t-s)}{ \mathbf{P}_x(\tau_0^->t)}\, ds .
  \eeqnn
 By \eqref{Asytau} and Theorem~\ref{MainThm07}(1), it converges as $t\to\infty$ to
  \beqnn 
  \frac{1}{\widehat{V}(x)} \int_{0}^{K} \mathbf{E}\Big[F\Big( \mathcal{P}  \cdot\mathbf{1}_{\{r\geq \mathcal{T} +s \}} , r\in[0,T] \Big)\Big]\cdot \bar{n}(  \epsilon_s\leq x,  \zeta>s )  \, ds,
  \eeqnn
  which further converges as $K\to\infty$ to
  \beqnn 
   \frac{1}{\widehat{V}(x)} \int_{0}^\infty \mathbf{E}\Big[F\Big( \mathcal{P}  \cdot\mathbf{1}_{\{r\geq \mathcal{T} +s \}} , r\in[0,T] \Big)\Big]\cdot \bar{n}(  \epsilon_s\leq x,  \zeta>s )  \, ds.
  \eeqnn
  Combining this together with \eqref{eqn.4371} and then taking them back into \eqref{eqn.444}, 
  \beqnn
  \lefteqn{\lim_{K\to\infty}\lim_{t\to\infty}  \mathbf{E}_x\big[ F( X_{s}/t,  s\in[0,T] ), \underline{g}_t\leq K\,\big|\,\tau_0^->t \big]}\ar\ar\cr
  \ar=\ar \frac{\bar{\tt d} }{\widehat{V}(x)}\cdot \mathbf{E}\big[F( \mathcal{P} , s\in[0,T])\big] + \frac{1}{\widehat{V}(x)} \int_{0}^\infty \mathbf{E}\Big[F\Big( \mathcal{P}  \cdot\mathbf{1}_{\{r\geq \mathcal{T} +s \}} , r\in[0,T] \Big)\Big]\cdot \bar{n}(  \epsilon_s\leq x,  \zeta>s )  \, ds.
  \eeqnn
  which equals to $\mathbf{E}\big[F\big(\mathcal{P}\cdot \mathbf{1}_{\{\mathcal{T}_x\leq s\}}, s\in [ 0,T] \big)\big]$ by Lemma~\ref{Lemma.408}. The proof ends.
  \qed

 \textit{\textbf{ Proof of Corollary~\ref{MainThm08}(2).}}
 By the independent increments of $X$, it suffices to prove that for any $\delta \in (0,1)$,
 \beqnn
  \Big\{\frac{X_{ts}}{t}:s\in[\delta,1] \Big\}\overset{\rm d}
  \to \big\{ \mathcal{P} -\beta s: s\in[\delta,1] \big\}
  \quad\mbox{and}\quad
  \Big\{\frac{X_{t(1+s)}-X_{t}}{t}:s\geq 0\Big\}\overset{\rm d}\to \big\{-\beta s: s\in[\delta,1] \big\},
 \eeqnn
 in $D\big([\delta ,1];\mathbb{R}\big)$ and $D\big([0,\infty);\mathbb{R}\big)$ as $t\to\infty$.
 The second limit follows directly from Proposition~\ref{Prop.402}. 
 We now prove the first one, which holds if and only if for any bounded, uniformly continuous non-negative functional $F$ on $D \big([\delta,1];\mathbb{R}\big)$,  
 \beqlb\label{eqn.433}
 \mathbf{E}_x\big[F\big(X_{ts}/t, s\in[\delta,1]\big)  \,\big|\,  \tau_0^->t \big] \to \mathbf{E}\big[F\big( \mathcal{P}^{\alpha,\beta}-\beta s, s\in[\delta,1]\big)\big],
 \eeqlb
 as $t\to\infty$. 
 For any  $K>0$, the left-hand side of \eqref{eqn.433} can be decomposed into the next two terms
 \beqlb\label{q1651}
 \mathbf{E}_x\big[F\big(X_{ts}/t, s\in[\delta,1]\big) , \underline{g}_t> K \,\big|\,  \tau_0^->t \big] 
 \quad \mbox{and}\quad 
 \mathbf{E}_x\big[F\big(X_{ts}/t, s\in[\delta,1]\big) ,\underline{g}_t\leq K \,\big|\,  \tau_0^->t \big] . 
 \eeqlb
 Without loss of generality, we still assume that $F(\cdot)\leq 1$. 
 By \eqref{eqn.4351}, 
 \beqlb\label{eqn.434}
  \lim_{K\to\infty}\lim_{t\to\infty}\mathbf{E}_x\big[F\big(X_{ts}/t, s\in[\delta,1]\big) , \underline{g}_t> K \,\big|\,  \tau_0^->t \big] 
  \leq  \lim_{K\to\infty}\lim_{t\to\infty}\mathbf{E}_x\big[ \underline{g}_t> K \,\big|\,  \tau_0^->t \big]   
  =0. 
 \eeqlb 
 Repeating the proof of   Corollary~\ref{MainThm08}(1) to the second conditional expectation in \eqref{q1651} induces that 
 \beqnn
  \mathbf{E}_x\big[F\big(X_{ts}/t, s\in[\delta,1]\big) ,\underline{g}_t\leq K \,\big|\,  \tau_0^->t \big]
  \sim \frac{\mathbf{E} \big[ F\big(   X_{ts}/t:s\in[\delta,1] \big), \underline{g}_t\leq K, \underline{X}_t\geq -x \big]}{\mathbf{P}_x( \tau_0^->t )}  ,
 \eeqnn
  as $t\to\infty$. 
 The last expectation can be represented as
 \beqnn
 \bar{\tt d}\cdot\underline{n}\Big(F\Big(\frac{\epsilon_{ts}}{t},s\in[\delta,1]\Big),\zeta>t\Big) \ar+\ar\int_{0}^{K}ds \int_{\mathcal{E}}    \overline{n}( d\epsilon,\epsilon_s\leq x,  \zeta>s )  \int_\mathcal{E}F\Big(\frac{(\epsilon,\epsilon^*)_{tr}^s}{t}, r\in[\delta,1] \Big)\, \underline{n}(d\epsilon^*, \zeta>t-s). 
 \eeqnn  
 Dividing both terms by $\mathbf{P}_x( \tau_0^->t )$, then applying  \eqref{Asytau} and Theorem~\ref{MainThm07}(2), they converge as $t\to\infty$ to 
 \beqnn
   \mathbf{E}\big[F\big( \mathcal{P} -\beta s, s\in[\delta,1] \big)\big]\cdot   \frac{1}{\widehat{V}(x)} 
    \Big( \bar{\tt d}  +
 \int_{0}^{K}\overline{n}(  \epsilon_s\leq x,  \zeta>s )  \, ds \Big).
 \eeqnn
 which further converges to $ \mathbf{E}\big[F\big( \mathcal{P} -\beta s, s\in[\delta,1] \big)\big]$ as $K\to\infty$. The proof ends.
 \qed

 \bigskip
 \textbf{Acknowledgment.} We would like to thank Jesus Contreras for comments on the intersection between lifetime and height of excursions in the oscillating case.

 \bibliographystyle{plain}
 
 \bibliography{Reference}

 \end{document}